\documentclass[10pt,reqno]{amsart}
\usepackage[T1]{fontenc}
\usepackage{lmodern,microtype}
\usepackage[margin=1in]{geometry}
\usepackage{amsmath,amssymb,mathtools,mathrsfs,bm,esint,graphicx,pgfplots}
\usepgfplotslibrary{groupplots,fillbetween}
\pgfplotsset{compat=1.18}
\usepackage{enumitem,booktabs,array}
\usepackage{tikz}
\usetikzlibrary{arrows.meta,calc}
\usepackage[colorlinks=true,linkcolor=blue!55!black,citecolor=blue!55!black,urlcolor=blue!55!black]{hyperref}
\usepackage[nameinlink,capitalize,noabbrev]{cleveref}
\numberwithin{equation}{section}
\allowdisplaybreaks[1]
\theoremstyle{plain}
\newtheorem{theorem}{Theorem}[section]
\newtheorem{proposition}[theorem]{Proposition}
\newtheorem{lemma}[theorem]{Lemma}
\newtheorem{corollary}[theorem]{Corollary}
\theoremstyle{definition}
\newtheorem{definition}[theorem]{Definition}
\theoremstyle{remark}
\newtheorem{remark}[theorem]{Remark}
\crefname{theorem}{Theorem}{Theorems}
\crefname{proposition}{Proposition}{Propositions}
\crefname{lemma}{Lemma}{Lemmas}
\crefname{corollary}{Corollary}{Corollaries}
\crefname{definition}{Definition}{Definitions}
\crefname{remark}{Remark}{Remarks}
\newcommand{\R}{\mathbb R}
\renewcommand{\S}{\mathbb S}

\renewcommand{\O}{\Omega}
\newcommand{\A}{\mathcal A}

\newcommand{\J}{\mathcal J}
\newcommand{\K}{\mathcal K}

\newcommand{\s}{\mathcal S}

\newcommand{\W}{\mathcal W}
\newcommand{\ep}{\varepsilon}

\newcommand{\p}{\partial}

\newcommand{\defeq}{\vcentcolon=}

\DeclareMathOperator{\Div}{div}

\DeclareMathOperator{\supp}{supp}

\DeclarePairedDelimiter{\norm}{\lVert}{\rVert}
\DeclarePairedDelimiter{\seminorm}{[}{]}
\DeclarePairedDelimiterX{\inner}[2]{\langle}{\rangle}{#1,#2}
\title{The Alt--Phillips Multiple Membrane Free Boundary Problem}
\author{Kunyi (Mark) Ma}
\address{Department of Mathematics, Columbia University, New York, NY 10027}
\email{km4046@columbia.edu}
\date{}
\subjclass[2020]{Primary 35R35; Secondary 35B65, 35J60, 49J40}
\keywords{Free boundary problems,
Alt--Phillips, Vectorial free boundary problems, Membranes problem}
\begin{document}
\begin{abstract}
Motivated by diffusion--reaction processes in porous catalysts,
we study an Alt--Phillips model for ordered membranes with
weighted gap interactions. The model connects adhesive harmonic junctions \cite{DSS24} and forced
multiple membranes \cite{SY19} through the exponent $0<\gamma<1$. For arbitrary positive weights, we prove optimal
$C^{1,\beta-1}$ regularity, where $\beta=2/(2-\gamma)$,
and obtain nontrivial $\beta$-homogeneous blowup limits
through Weiss almost monotonicity.
Under strict discrete concavity of the weights, we establish
a gluing principle, minimality of flat single-gap cones,
and $C^{1,\alpha}$ free boundary regularity at points
admitting such blowups. We also classify all one-dimensional vertically critical cones. For three membranes, we identify cone families that are
non-minimizing for small $\gamma$ and minimizing near
$\gamma=1$.
\end{abstract}
\maketitle
\tableofcontents

\section{Introduction}\label{sec:introduction}

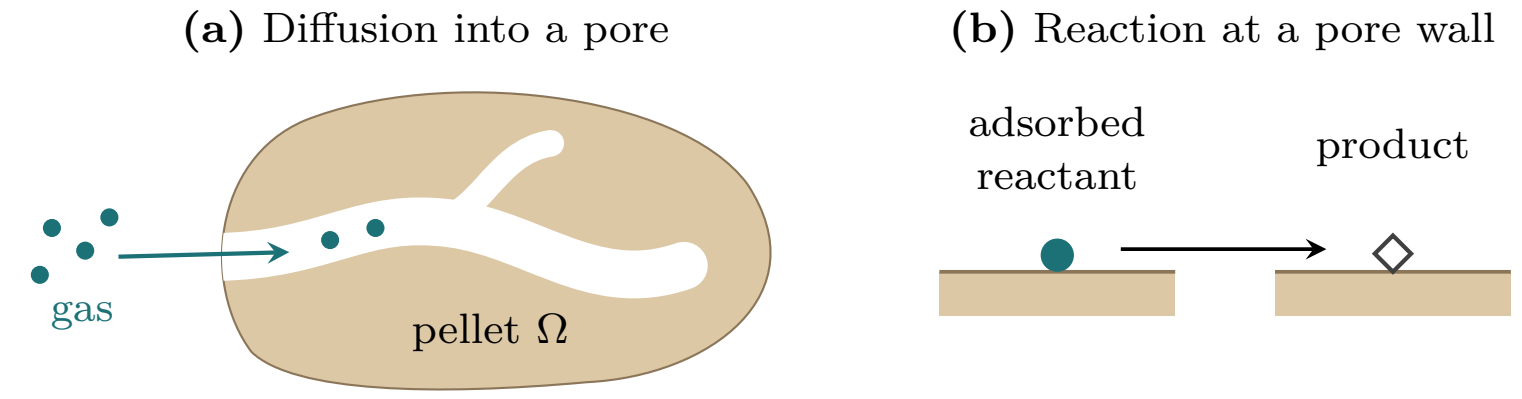
\begin{figure}[htbp]
  \centering
  \setlength{\abovecaptionskip}{1.2pt}
  \begin{tikzpicture}[x=1cm,y=1cm,font=\footnotesize,
      flow/.style={->,>=stealth,draw=pelletGas,line width=.8pt}]
    \definecolor{pelletSolid}{RGB}{221,200,166}
    \definecolor{pelletRim}{RGB}{139,118,90}
    \definecolor{pelletGas}{RGB}{28,113,119}

    % (a) Reactant diffuses from the surrounding gas into a pore.
    \node at (2.15,2.48) {\textbf{(a)} Diffusion into a pore};

    \def\pellet{(1.00,.36)
      .. controls (.67,.80) and (.74,1.63) .. (1.35,1.89)
      .. controls (2.25,2.25) and (3.83,2.08) .. (4.28,1.47)
      .. controls (4.76,.74) and (4.03,.20) .. (3.22,.16)
      .. controls (2.22,.07) and (1.28,.09) .. cycle}
    \filldraw[fill=pelletSolid,draw=pelletRim] \pellet;

    \begin{scope}
      \clip \pellet;
      \draw[white,line width=9pt,line cap=round]
        (.62,.99) .. controls (1.42,.96) and (1.66,1.26)
        .. (2.21,1.22) .. controls (2.82,1.19) and (3.23,.71)
        .. (3.86,.93);
      \draw[white,line width=5pt,line cap=round]
        (2.16,1.22) .. controls (2.52,1.22) and (2.59,1.68)
        .. (2.98,1.74);
    \end{scope}

    \foreach \x/\y in {-.40/.87,-.32/1.18,-.10/1.03,.06/1.25,
                         1.52/1.10,1.82/1.18}
      \fill[pelletGas] (\x,\y) circle (1.7pt);
    \draw[flow] (.12,.99) -- (1.25,1.02);
    \node[text=pelletGas] at (-.12,.60) {gas};
    \node at (2.58,.49) {pellet $\O$};

    % (b) The same wall site before and after the reaction.
    \begin{scope}[shift={(5.30,0)}]
      \node at (2.12,2.48) {\textbf{(b)} Reaction at a pore wall};

      \foreach \x in {1.03,3.25}{
        \fill[pelletSolid] (\x-.78,.60) rectangle (\x+.78,.89);
        \draw[pelletRim,line width=.7pt]
          (\x-.78,.89) -- (\x+.78,.89);
      }

      \fill[pelletGas] (1.03,1.00) circle (.11);
      \node[align=center] at (1.03,1.71) {adsorbed\\reactant};

      \draw[->,>=stealth,line width=.7pt]
        (1.45,1.04) -- (2.81,1.04);

      \draw[black!75,line width=.8pt,fill=white]
        (3.25,.89) -- (3.37,1.01) -- (3.25,1.13)
        -- (3.13,1.01) -- cycle;
      \node at (3.25,1.71) {product};
    \end{scope}
  \end{tikzpicture}
  \caption{Reactant molecules diffuse into the pellet, attach to a pore wall
    and react.}
  \label{fig:intro-catalyst-pellet}
\end{figure}

Consider a reactant gas that diffuses into the pores of a catalyst pellet $\O$. When gas molecules attach to the pore walls, chemical reaction happens. One is interested in the steady state, where the net influx of gas by diffusion
balances the amount consumed by the reaction per unit time; see \Cref{fig:intro-catalyst-pellet}.

At steady state, the normalized gas concentration $u_\varepsilon(x) \in [0, 1]$ satisfies
\[
  \Delta u_\varepsilon=R_\varepsilon(u_\varepsilon)
  \quad\text{in }\O,
  \qquad
  u_\varepsilon=1 \quad\text{on }\p\O,
\]
with consumption rate $R_\varepsilon$, which is the Euler--Lagrange equation of
\[
  \int_\O\left(\frac12|\nabla v|^2+F_\varepsilon(v)\right)\,dx,
  \qquad
  F_\varepsilon(s):=\int_0^s R_\varepsilon(t)\,dt.
\]
The Langmuir--Hinshelwood law in \cite{BN80LNM} admits a limit of the form
\begin{equation}\label{eq:Alt-Phillips-motivation}
  R_\varepsilon(s)\to\lambda s^{\gamma-1},
  \qquad
 F_\varepsilon(s)
  \to\frac{\lambda}{\gamma}s^\gamma,
  \qquad \lambda>0,\quad 0<\gamma<1,
\end{equation}
as $\varepsilon\downarrow0$, for each $s>0$, and recovers the Alt--Phillips problem \cite{Phi83,AP86}
\begin{equation*}
  \min_{\substack{v\in H^1(\O),\ v\geq0\\
                  v=1\ \mathrm{on}\ \p\O}}
  \int_\O\left(\frac12|\nabla v|^2
               +\frac{\lambda}{\gamma}v^\gamma\right)\,dx.
\end{equation*}
A minimizer $u$ represents a limiting steady state gas concentration,
satisfying $\Delta u=\lambda u^{\gamma-1}$ in $\{u>0\}$, and
its free boundary $\O\cap\p\{u>0\}$ is the depletion front of the reactant gas.

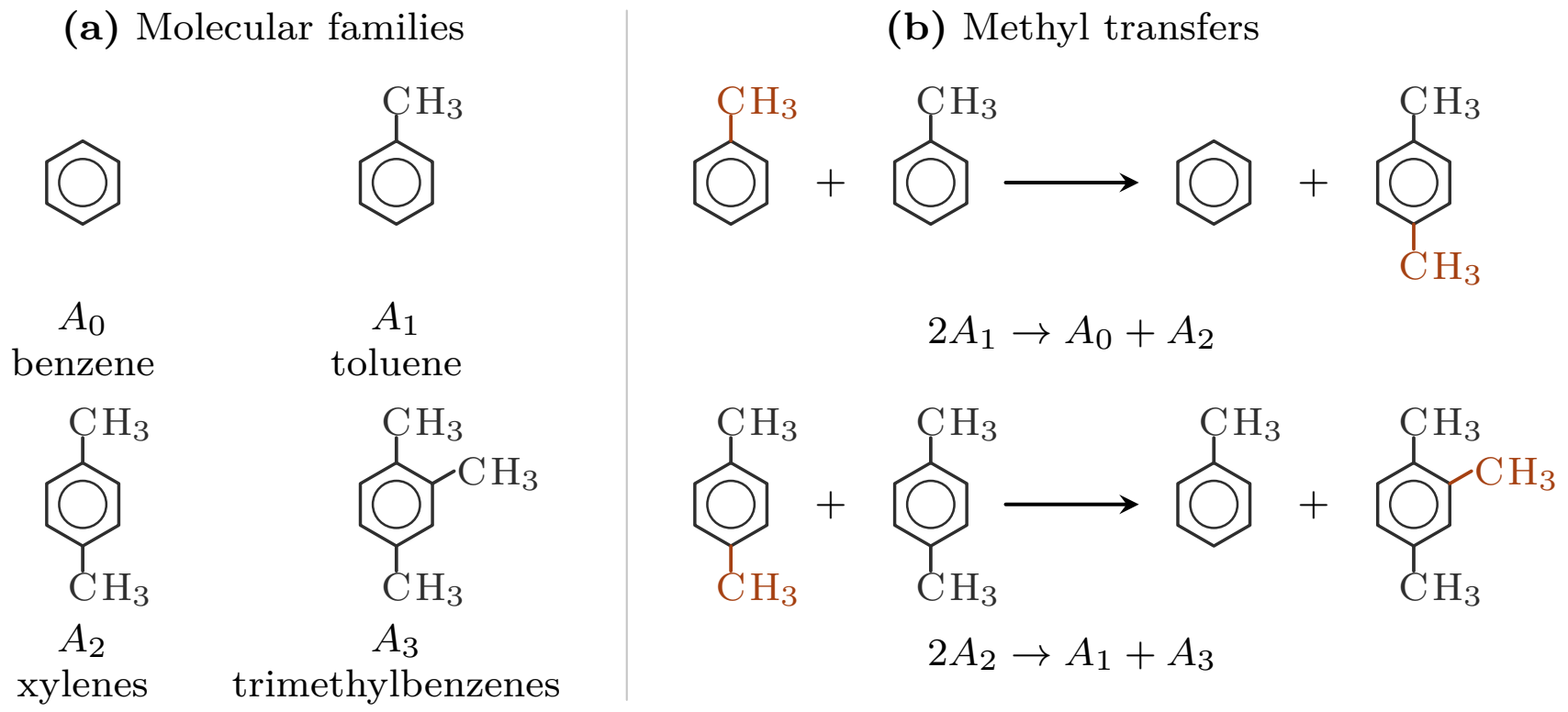
\begin{figure}[htbp]
  \centering
  \vspace{-3pt} 
\setlength{\abovecaptionskip}{1pt}
  \definecolor{transferMethyl}{HTML}{A54011}
  \begin{tikzpicture}[x=1cm,y=1cm,line cap=round,line join=round,
      font=\footnotesize, every node/.style={inner sep=0pt},
      ring/.style={draw=black!82,line width=.65pt}]
    % Center, number of methyl groups, index of the highlighted group (0 = none).
    \def\Aromatic#1#2#3{%
      \begin{scope}[shift={#1}]
        \draw[ring] (90:.30) -- (30:.30) -- (-30:.30) -- (-90:.30)
          -- (-150:.30) -- (150:.30) -- cycle;
        \draw[ring,line width=.48pt] (0,0) circle[radius=.17];
        \foreach \slot/\ang/\anch in {1/90/south,2/-90/north,3/30/west}{%
          \ifnum\slot>#2\relax\else
            \def\bondcolor{black!82}%
            \ifnum\slot=#3\relax\def\bondcolor{transferMethyl}\fi
            \draw[draw=\bondcolor,line width=.7pt] (\ang:.30) -- (\ang:.48);
            \node[anchor=\anch,text=\bondcolor,inner sep=.35pt]
              (methylC) at (\ang:.48) {$\mathrm{C}$};
            \node[anchor=base west,text=\bondcolor] at (methylC.base east)
              {$\mathrm{H}_3$};
          \fi}
      \end{scope}}

    \node at (1.30,1.10) {\textbf{(a)} Molecular families};
    \foreach \x/\y/\j/\family in
      {0/0/0/benzene,2.26/0/1/toluene,
       0/-2.32/2/xylenes,2.26/-2.32/3/trimethylbenzenes}{%
      \Aromatic{(\x,\y)}{\j}{0}
      \node[anchor=north,align=center] at (\x,\y-.86)
        {$A_{\j}$\\[-.2ex]\family};
    }

    \draw[black!22,line width=.4pt] (3.92,1.24) -- (3.92,-3.73);
    \begin{scope}[xshift=4.67cm]
      \node at (2.46,1.10) {\textbf{(b)} Methyl transfers};
      \foreach \y/\j/\lo/\hi in {0/1/0/2,-2.32/2/1/3}{%
        \Aromatic{(0,\y)}{\j}{\j}
        \node at (.72,\y) {$+$};
        \Aromatic{(1.44,\y)}{\j}{0}
        \draw[->,>=stealth,line width=.75pt] (1.98,\y) -- (2.94,\y);
        \Aromatic{(3.48,\y)}{\lo}{0}
        \node at (4.20,\y) {$+$};
        \Aromatic{(4.92,\y)}{\hi}{\hi}
        \node at (2.46,\y-1.08) {$2A_{\j}\to A_{\lo}+A_{\hi}$};
      }
    \end{scope}
  \end{tikzpicture}
  \caption{Molecule structures for $A_0,\ldots,A_3$ and the two transfers.}
  \label{fig:intro-methyl-transfers}
  \label{fig:intro-aromatics}
\end{figure}

In this paper, we study a vectorial Alt--Phillips problem
for $N\geq2$ ordered membranes, minimizing
\begin{equation}
    \label{eq:Alt-Phillips}
  \int_\O\left(
    \sum_{i=1}^N|\nabla u_i|^2
    +\sum_{i=1}^{N-1}w_i(u_i-u_{i+1})^\gamma
  \right)\,dx,
  \qquad w_i>0,\quad 0<\gamma<1,
\end{equation}
over $U=(u_1,\ldots,u_N)$ with $u_1\geq\cdots\geq u_N$
and prescribed boundary data. The `ordering' can represent cumulative concentrations
of several interacting chemicals, as we illustrate for $N=3$ using an experiment by Beltrame et al.~\cite{Beltrame85}. In this experiment, toluene and xylenes were passed
separately over the porous zeolite HZSM-5, undergoing reactions
that redistribute methyl groups ($\mathrm{CH}_3$) between benzene rings. Writing $A_j$ for molecules with
$j$ attached $\mathrm{CH}_3$, $j=0,1,2,3$, we retain the two
transfers in \Cref{fig:intro-methyl-transfers}.

Let $c_j$ be the molar concentration of $A_j$ and $u_i$ be those carrying at least $i$
methyl groups. Then
\begin{equation*}
     u_1=c_1+c_2+c_3\ \geq \ u_2=c_2+c_3\ \geq \ u_3=c_3.
\end{equation*}
A first transfer in \Cref{fig:intro-methyl-transfers} changes  $U=(u_1,u_2,u_3)$ by $(-1,1,0)$, and a second by $(0,-1,1)$. Writing $q_1,q_2$ for the two molar transfer rates, with diffusivity normalized to $2$, the steady state balance writes
\[
    2\Delta U+q_1(-1,1,0)+q_2(0,-1,1)=0,
  \qquad\text{or equivalently}\qquad
  2\Delta U=(q_1,\,q_2-q_1,\,-q_2).
\]
Motivated by \eqref{eq:Alt-Phillips-motivation}, we assume
the transfer rates
\[
  q_i=\lambda_i c_i^{\gamma-1}
     =\gamma w_i(u_i-u_{i+1})^{\gamma-1},
  \qquad c_i=u_i-u_{i+1}>0,\quad i=1,2,
\]
with $\lambda_i=\gamma w_i$.
The steady balances are then precisely the Euler--Lagrange
equations of \eqref{eq:Alt-Phillips} for $N=3$ on $\{u_1>u_2>u_3\}$. A minimizer $U$ models steady cumulative concentrations, and 
the free boundaries
$\Gamma_i(U)=\O\cap\p\{u_i>u_{i+1}\}$, $i=1,2$,
mark the depletion fronts of toluene and xylenes, respectively.

\subsection{The membranes model \texorpdfstring{\(\J_{N,d}^{\gamma,G}\)}{J}}
We now formulate the model over a bounded Lipschitz domain $\O\subset\R^d$. Denoting $U=(u_1,\ldots,u_N): \O \to \R^N$ as the vertical positions of the membranes, we consider
\begin{equation*}
 \J_{N,d}^{\gamma,G}(U,\mathbf{w}_\gamma,\O)
 =\int_\O\left(|\nabla U|^2+W_N^\gamma(U,\mathbf{w}_\gamma)
      +G(x,U)\right)dx.
\end{equation*}
with
\[|\nabla U|^2=\sum_{i=1}^N|\nabla u_i|^2, \qquad  W_N^\gamma(U,\mathbf{w}_\gamma)=\sum_{i=1}^{N-1}
       w_{i}(u_i-u_{i+1})^\gamma, \qquad G \in \mathbf{G}_K, \]
where 
\begin{enumerate}
    \item $\gamma \in (0, 1)$ denotes the \textit{Alt-Phillips} parameter. 
    \item $\mathbf{w}_{\gamma} = (w_{1}, \cdots, w_{N-1}) \in (0, \infty)^{N-1}$ are weights for the \textit{Alt-Phillips} potential depending on $\gamma$.
    \item $G$ is a perturbation with Lipschitz bound $K > 0$ in the $U$-variable. Denote
   \begin{equation}\label{eq:g-lip}
\mathbf G_K:=
\left\{
G:\O\times\R^N\to\R
\;\middle|\;
\begin{aligned}
 &G(x,0)=0,\ 
 x \mapsto G(x,U)\text{ is measurable for every }U\in\R^N,\\
 &|G(x,U)-G(x,V)|\leq K|U-V| \ \text{for a.e. }x\in\O\text{ and all }U,V\in\R^N
\end{aligned}
\right\}.
\end{equation}
\end{enumerate}
We write $\J_{N, d}^\gamma$ if $G=0$, and suppress arguments \(\mathbf{w}_\gamma, \O\) when no ambiguity arises. Also set $W_1^\gamma \equiv0$.

We're interested in the minimization problem for $\J_{N, d}^{\gamma, G}(\cdot, \mathbf{w}_\gamma, \O)$ over the admissible class 
\begin{equation}\label{eq:admissible-class}
    \A =\{U\in H^1(\O; \R^N):u_1\ge\cdots\ge u_N\}
\end{equation}
with respect to some boundary data $\Psi \in \A$. The natural homogeneity for \(\J_{N, d}^{\gamma, G}\) is
\begin{equation}\label{eq:homogeneity}
 \beta=\frac{2}{2-\gamma} \in (1, 2).
\end{equation} Indeed, $U$ locally minimizes \(\J_{N, d}^{\gamma, G}(\cdot,\mathbf{w}_\gamma, B_r(x_0))\) iff the rescaling \(U_{x_0,r}(x)=r^{-\beta}U(x_0+rx)\) locally minimizes \(\J_{N, d}^{\gamma, G_{x_0, r}}(\cdot,\mathbf{w}_\gamma, B_1)\) with the rescaled perturbation 
\begin{equation}\label{eq:rescale-G}
G_{x_0,r}(x,V)=r^{2-2\beta}G(x_0+rx,r^\beta V) \in \mathbf{G}_{Kr^{2-\beta}} \subseteq \mathbf{G}_K, \qquad \forall \ r \in (0, 1).
\end{equation}
 In particular, the perturbation vanishes as $r \downarrow 0$ for $\gamma \in (0, 1)$.

A key motivation for our model is to interpolate between the two endpoint cases. For $\gamma =0$ consider the adhesive harmonic junctions model \cite{DSS24},
\[\J_{N, d}^0(U, \O) = \int_\O \left(|\nabla U|^2 + \sum_{i=1}^{N-1} \mathbf{1}_{\{u_i > u_{i+1}\}}\right)dx,\] and for $\gamma = 1$ consider the forced multiple-membrane problem \cite{SY19,SY20,SY23}
\[\J_{N, d}^1(U, \{f_j\}, \O) = \int_\O \left(|\nabla U|^2 + 2 \sum_{i = 1}^N f_i u_i\right)dx, \qquad f_1 > f_2 > \cdots > f_N. \]
In the following, we discuss our choice of $\J_{N, d}^{\gamma, G}$ and the assumptions that make such interpolation precise.

\subsubsection{Lipschitz perturbation $\mathbf{G}_K$}\label{sec:Lipschitz perturbation}

For `membrane models', one show optimal regularity by inducting on the number of membranes \cite[Lemma~2.9]{DSS24}\cite[Proposition~2.8]{SY19}.
The induction relies on that strictly separated blocks locally minimize
the same type of energy with fewer membranes, as in the endpoint
models $\gamma=0$ \cite[Lemma~2.5]{DSS24} and
$\gamma=1$ \cite[Remark~2.2]{SY19}.
For $0<\gamma<1$, nonlinear gap interactions persist under
strict separation, so this property can fail for
$\J_{N,d}^{\gamma}$; see \Cref{rmk:separated-blocks-interaction}. These interactions, however, are Lipschitz in gaps bounded away from zero. Therefore allowing perturbations in $\mathbf G_K$ absorbs them,
and separated blocks remain minimizers within the same energy
class;
see \Cref{prop:local-splitting}.

\subsubsection{Weight assumptions $\mathbf{w}_\gamma$}\label{sec:Weight assumptions}
We assume positive weights $\mathbf w_\gamma\in(0,\infty)^{N-1}$,
which penalize separation between adjacent membranes. Positivity alone suffices for the results in
\Cref{sec:preliminaries,sec:optimal-regularity,sec:monotonicity-blowups}. However for `membrane-collapsing' arguments in
\Cref{sec:nondegenracy,sec:single-gap-regularity},
we additionally impose strict discrete concavity
\begin{equation}\label{eq:weight-strict-concavity}
  2w_i>w_{i-1}+w_{i+1},
  \qquad i=1,\ldots,N-1,
  \qquad w_0=w_N=0.
\end{equation}
The motivation comes from $\gamma =1$, where $\mathbf{w}_1 $ satisfies 
\[\sum_{i = 1}^{N-1} w_{ i}(u_i - u_{i+1}) = 2 \sum_{i=1}^N f_i u_i, \qquad f_i = \frac{1}{2}(w_{i} - w_{i-1}).\]
\eqref{eq:weight-strict-concavity} corresponds exactly to the nondegenerate constant-force assumption $f_i > f_{i+1}$ from \cite{SY19}. In \Cref{sec:N3-d1-minimizing-cones} we study $1$-d cones with three membranes, where \eqref{eq:weight-strict-concavity} is fully-characterized by the ratio
\begin{equation*}
     \mathbf{w}_\gamma=(w_{ 1},w_{2}),\qquad \frac12< \rho = \frac{w_{ 2}}{w_{ 1}}<2.
\end{equation*}
Strict concavity allows certain cones to inherit minimality from $\J_{3, 1}^1$ for $\gamma$
close to $1$; see \Cref{thm:type-ii-minimality,thm:type-iv-minimality}.

\subsection{Background and related work}\label{sec:background}
The scalar Alt--Phillips problem has a rich theory on
minimizer and free boundaries \cite{Phi83,AP86,DSS21,RR25}, with recent works on generic regularity, stable and
minimizing cones \cite{FY23,KSP26,SY25,SY26Concentration}.
Variants include negative-power \cite{DSS22},
two-phase \cite{LP08}, quasilinear \cite{LQT15},
fully nonlinear \cite{AT13,WY22,Zahl26}, nonlocal \cite{Yang13}, and parabolic \cite{AS24,AS26}.

For vectorial problems, the two endpoint models with ordering constraint are most relevant to our work. At $\gamma=1$, classical results include existence, uniqueness
and $W^{2,p}$ regularity \cite{ChipotVergaraCaffarelli1985},
two-membrane $C^{1,1}$ and free-boundary regularity \cite{Sil05},
and quadratic gap growth \cite{LR09}.
Variants include quasilinear systems \cite{ARS05},
nonlocal constraints \cite{CCV05}, and different operators
\cite{CDS17}.
Savin and Yu proved optimal $C^{1,1}$ regularity and planar
blow-up classification \cite{SY19}, followed by $C^{1,\log}$
regularity near regular intersections for $N=3$ \cite{SY20}
and for general $N$ in the plane \cite{SY23}. At $\gamma=0$ with unit weights, De Silva and Savin
\cite{DSS24} established optimal Lipschitz and partial
free-boundary regularity, classified planar three-membrane
minimizing junction cones, and proved piecewise $C^{1,\alpha}$
regularity near planar junctions.

Other related vectorial problems include Bernoulli systems
\cite{CSY18,MTV20,DT20,TV26}, prescribed contact volume
\cite{SV26Rigidity}, transmission \cite{JS25Transmission},
local and long-range segregation
\cite{CTV05,CPQT19,CPT26,CT26Segregation},
and optimal eigenvalue partitions \cite{CL07}.
In the vectorial Alt--Phillips setting, results include
regularity of minimizers and almost minimizers
\cite{FSW21,DSJS23}, analyticity of regular free boundaries
\cite{FK24,DSJS26}, and partial regularity of minimizing
constraint maps with distance-power potentials \cite{Zig26}.

\subsection{Main results}\label{sec:main-results}Our main results are as following.
\begin{enumerate}[label=\textup{(\roman*)},leftmargin=*]
\item \emph{Optimal regularity and blowups.}
For positive weights and $G\in\mathbf G_K$, minimizers are
locally $C^{1,\beta-1}$, where $\beta=2/(2-\gamma)$.
After restricting to a maximal contact block and subtracting
its mean, free boundary blowups are nontrivial
$\beta$-homogeneous global minimizers with $G=0$
(\Cref{thm:optimal-regularity,thm:gap-blowup-cone}).

\item \emph{Flat single-gap cone.}
Under concavity
\eqref{eq:weight-strict-concavity}, flat single-gap cones are
global minimizers (\Cref{prop:single-gap-minimizers}).
Near such blowups, gluing reduces the problem to scalar
Alt--Phillips with a Lipschitz perturbation, yielding
$C^{1,\alpha}$ free boundary regularity
(\Cref{lem:membrane-block-collapse,cor:single-gap-regularity,sec:scalar-perturbed-flatness}).

\item \emph{One-dimensional cones.}
Even classifying $1$-d minimizing cones is complicated due to the nonlinearity in \eqref{eq:weighted-critical-coefficient-system}. For $G=0$, arbitrary $N$, and positive weights, we first classify
all $1$-d critical cones and exclude minimality
when a gap is active on both half-lines
(\Cref{thm:one-dimensional-critical-cones,lem:common-active-gap-instability}).
For $N=3$, the remaining candidates beyond single-gap cones are Type~II (one-sided, two-gap) and Type~IV (two-sided,
different active gaps). Type~II is non-minimizing for $\gamma<(14+4\sqrt6)/25$,
and Type~IV for $\gamma\leq2/3$; both minimize sufficiently
near $\gamma=1$ when $w_2/w_1\to\rho_*\in(1/2,2)$
(\Cref{thm:type-ii-minimality,thm:type-iv-minimality}).
\end{enumerate}

\subsection*{Acknowledgments}
The author thanks his advisors Ovidiu Savin and Daniela De Silva
for introducing the problem, fruitful discussions,
and continued support.

\subsection*{AI Disclosure}
The author used GPT-6.0 Astra to construct the Type~II
and Type~IV calibrations
\eqref{eq:type-ii-explicit-calibration} and
\eqref{eq:typeIVcalibration}, search for background and references,
check writing and computations, and generate figures.
The author has reviewed and understood all AI-assisted
material and takes full responsibility for the work.

\section{Preliminaries}\label{sec:preliminaries}

In this section, $N \ge 2,\, \O\subseteq \R^d \text{ is bounded Lipschitz}, \gamma \in (0, 1), G \in \mathbf{G}_K \text{ for } K >0, \mathbf{w}_\gamma \in (0, \infty)^{N-1},$ unless otherwise specified. We prove some first results for a minimizer to \(\J_{N, d}^{\gamma, G}\) among $\A$ as in \eqref{eq:admissible-class}.

\subsection{Existence, \texorpdfstring{$L^\infty$}{J} and \texorpdfstring{$C^{0, \alpha}$}{J} estimates.}

One begin with existence of minimizer.
\begin{proposition}\label{prop:existence}
For $\Psi\in\A$, a minimizer to \(\inf\{\J_{N, d}^{\gamma, G}(V,\mathbf{w}_\gamma,\O):
          V\in\A,\, V-\Psi\in H_0^1(\O; \R^N)\}\) exists.
\end{proposition}
\begin{proof}
Note \(\J_{N, d}^{\gamma, G}(\Psi)< \infty.\) Also, for any $V \in \A$ with boundary $\Psi$, using Poincaré and Young's Inequality,
\begin{align*}
    \J_{N, d}^{\gamma, G}(V) &\geq \norm{\nabla V}_{L^2}^2 - K |\O|^{1/2}\norm{V}_{L^2}\geq \norm{\nabla V}_{L^2}^2 -  K |\O|^{1/2}\left(\norm{V - \Psi}_{L^2} + \norm{\Psi}_{L^2} \right)\\
    &\geq \norm{\nabla V}_{L^2}^2 - C(\O) K\left(\norm{\nabla (V -\Psi)}_{L^2} + \norm{\Psi}_{L^2} \right)\geq \frac{1}{2}\norm{\nabla V}_{L^2}^2  - C(\O)\left(K^2 + \norm{\Psi}_{H^1}^2\right).
\end{align*}
Thus $-\infty < \inf\J_{N, d}^{\gamma, G}(V) < \infty.$ One take minimizing sequence $\{U_m\} \subseteq \A$ with trace $\Psi$ that has a weak $H^1$, strong $L^1$ and a.e. pointwise limit $U \in \A.$ $U$ attains the infimum due to weak lower-semicontinuity of Dirichlet energy, $\norm{W_N^\gamma(U_m)-W_N^\gamma(U)}_{L^1}\to 0$ from sub-additivity, and $\norm{G(x, U_m) - G(x, U)}_{L^1}\leq K\norm{U_m - U}\to 0.$
\end{proof}

One has one-sided bounds on the Laplacian for top $u_1$ and bottom $u_N$ membranes.
\begin{lemma}\label{lem:linfty-bound}
Let $U \in \A$ minimize $\J_{N, d}^{\gamma, G}$ with boundary $\Psi \in \A.$ Then 
\begin{equation*}
 \Delta u_1\ge-K,\qquad \Delta u_N\le K
 \qquad\hbox{weakly}.
\end{equation*}
Moreover, if $\Psi \in L^\infty(\O),$ then $\|U\|_{L^\infty(\O)}
\le C(d, N, \O)(\|\Psi\|_{L^\infty(\partial\O)}+K)$.
\end{lemma}
\begin{proof}
On any ball $B\Subset\O$, let $v-u_1\in H_0^1(B)$ solve
$\Delta v=-K$. Set $\widehat u_i=\min\{u_i,v\}$, and since ordering is preserved,
$\widehat{U} = (\hat{u}_1, \cdots, \hat{u}_N) \in \A$ with same boundary $\Psi$. Since every gap decreases,
$W_N^\gamma(\widehat{U}) \leq W_N^\gamma(U)$;
Also one may write $u_i = \widehat{u}_i + \phi_i$ for
$\phi_i=(u_i-v)_+ \in H_0^1(B)$, then $\nabla\widehat u_i=\nabla v$ where
$\phi_i>0$, and in particular
\[
\sum_i\int_B
\left(
|\nabla \widehat u_i|^2-|\nabla u_i|^2
\right)dx
=
\sum_i\int_B
\left(
-2\nabla \widehat u_i\cdot \nabla \phi_i
-|\nabla \phi_i|^2
\right)dx
\leq
-2K\sum_i\int_B \phi_i\,dx
-\sum_i\int_B |\nabla\phi_i|^2\,dx.
\]
Finally
$G(x, \widehat{U}) - G(x, U)
\leq K|\widehat{U} - U|
\leq K \sum_i \phi_i$,
so minimality gives
\[
 0\le \J_{N, d}^{\gamma,G}(\widehat U,B)
       -\J_{N, d}^{\gamma,G}(U,B)
 \le-\sum_i\int_B|\nabla\phi_i|^2
     -K\sum_i\int_B\phi_i \leq 0.
\]
Thus \(\phi_i\equiv0\) a.e. for all \(i\), and in particular
\(u_1\leq v\) in \(B\), which proves $\Delta u_1 \ge -K$ in $B$. \(\Delta u_N \leq K\) follows similarly, and maximum principle gives $L^\infty$ bound using boundary data.
\end{proof}

We show Hölder regularity for minimizer to small perturbations of Dirichlet energy with at most quadratic growth; See \cite{Ma26} for a detailed treatment. By rescaling, the result applies to our minimizer for $\J_{N, d}^{\gamma, G}$.

\begin{proposition}\label{prop:perturbative-holder}
Let $\alpha\in(0,1)$ and $\Lambda_F\ge1$. There are
$\delta_0,C_0>0$ depending only on $d,N,\alpha,\Lambda_F$ s.t. if $U$ minimize 
\begin{equation}\label{eq:f-growth}
\int_{B_1}(|\nabla V|^2+\delta F(x,V)), \qquad  |F(x,V)|\le \Lambda_F(1+|V|^2), \qquad 0<\delta\le\delta_0,
\end{equation}
among competitors taking values in a closed convex set $\mathcal K\subset\R^N$ with $\fint_{B_1}|U|^2\le1$, then \(
 \|U\|_{C^{0,\alpha}(\overline{B_{1/2}})}\le C_0.\)
\end{proposition}
\begin{proof}
\emph{Step 1}. First assume \eqref{eq:f-growth} with  $\Lambda \geq 1$ in place of $\Lambda_F$. Since $\K$ is convex, $H$ as harmonic replacement of $U$ in $B_1$ is admissible competitor. For $C = C(d, N) > 0$ changing from line to line, Poincar\'e gives
\[
\begin{aligned}
\int_{B_1}|\nabla (U - H)|^2
 &
 \le\delta\int_{B_1}\bigl(F(x,H)-F(x,U)\bigr)\le\delta\Lambda\int_{B_1}(2+|U|^2+|H|^2)
 \le C\delta\Lambda
 \left(1+\int_{B_1}|\nabla (U - H)|^2\right),
\end{aligned}
\]
so after absorbing the last term by choosing $\delta_0$ small depending on
$\Lambda$
\begin{equation*}
 \|\nabla (U - H)\|_{L^2(B_1)}^2+\|U - H\|_{L^2(B_1)}^2
 \le C\delta\Lambda.
\end{equation*}
In particular, harmonic estimates yields for $0<\rho<1/2$,
\[
    \fint_{B_\rho}|H - H(0)|^2   
    \leq C \rho^2, \qquad \fint_{B_\rho}|U-H(0)|^2
    \le C\rho^{-d}\delta\Lambda+C\rho^2.\]
Fix $C$ on the right-hand side. Choose $\rho=\rho(d,N,\alpha)$ s.t.
$C\rho^2\le\tfrac12\rho^{2\alpha}$, then
$\delta_0\Lambda\le (2C)^{-1}\rho^{d+2\alpha}$ s.t.
\begin{equation}
    \fint_{B_\rho}|U-M|^2\le\rho^{2\alpha},
    \qquad M = H(0),\qquad |M| \leq C_1 = C_1(d, N).
    \label{eq:harmonic-replacement-one-step}
\end{equation}

\emph{Step 2}. Now recover \eqref{eq:f-growth}. Fix
\[
 C_\alpha=\frac{C_1}{1-\rho^\alpha},\qquad
 \Lambda=\Lambda_F(3+2C_\alpha^2),\qquad
 \delta_0=\frac{\rho^{d+2\alpha}}{2C\Lambda}.
\]
Starting with $P_0=0$, define
\[
 U_k(x)=\rho^{-k\alpha}\bigl(U(\rho^kx)-P_k\bigr),
 \qquad
 \mathcal K_k=\rho^{-k\alpha}(\mathcal K-P_k),
 \qquad
 F_k(x,Z)=\rho^{2k(1-\alpha)}
 F(\rho^kx,P_k+\rho^{k\alpha}Z).
\]
Suppose inductively that
\[
 \fint_{B_1}|U_k|^2\le1,\qquad
 |P_k|\le C_1\sum_{j=0}^{k-1}\rho^{j\alpha}\le C_\alpha.
\]
Then since $\mathcal K_k$ is closed and convex, and
\( |F_k(x,Z)|\le \Lambda_F\rho^{2k(1-\alpha)}
       \bigl(1+2|P_k|^2+2\rho^{2k\alpha}|Z|^2\bigr)
 \le\Lambda(1+|Z|^2),
\)
applying \eqref{eq:harmonic-replacement-one-step} to $U_k$ gives $M_k$ with
$|M_k|\le C_1$ and
$\fint_{B_\rho}|U_k-M_k|^2\le\rho^{2\alpha}$. Set $P_{k+1}=P_k+\rho^{k\alpha}M_k$. Then
\[
 \fint_{B_1}|U_{k+1}|^2
 =\rho^{-2\alpha}\fint_{B_\rho}|U_k-M_k|^2\le1, \qquad |P_{k+1}|\leq C_1\left(\sum_{j=0}^{k-1}\rho^{j\alpha} + \rho^{k\alpha} \right) \leq C_\alpha,
\]
which closes the induction. Consequently, $P_k\to P_\infty \in \R^N$, with
\[
 \fint_{B_{\rho^k}}|U-P_\infty|^2
 \le 2(1+C_\alpha^2)\rho^{2k\alpha}, \qquad |P_\infty|\le C_\alpha.
\]
The argument applies uniformly for
$x_0\in B_{1/2}$. Conclude via Campanato's as in
\cite{Campanato1964Famiglia,HL11}.
\end{proof}

\begin{corollary}\label{cor:holder-minimizers}
For every $\alpha\in(0,1)$ and every minimizer $U$ to $\J_{N, d}^{\gamma,G}$ in $B_1$,
\begin{equation*}
 \|U\|_{C^{0,\alpha}(\overline{B_{1/2}})}
 \le C(d,N,\alpha,\gamma,K,\|\mathbf{w}_\gamma\|_{\ell^\infty})\left(1+\|U\|_{L^2(B_1)}\right).
\end{equation*}
\end{corollary}
\begin{proof}
Writing $a=A(1+\|U\|_{L^2(B_1)})$, notice $U/a$ minimizes
\eqref{eq:f-growth} with
\[
F_a(x,V)=a^{\gamma-1}W_N^\gamma(V,\mathbf{w}_\gamma)+a^{-1}G(x,aV), \qquad  0 < \delta=a^{-1} \leq A^{-1}.
\]
For $F_a$ with growth bound
$\Lambda_F = \Lambda_F(\gamma, N,\|\mathbf{w}_\gamma\|_{\ell^\infty}, K)$ that is independent of
$a$, choose $\delta_0 = \delta_0(d, N, \alpha, \Lambda_F)$ from
\Cref{prop:perturbative-holder}. Then set $A \geq \delta_0^{-1}$
up to constant depending on $d$, so \Cref{prop:perturbative-holder}
applies.
\end{proof}

\subsection{First Variation}\label{sec:euler-lagrange} 
\subsubsection{Vertical Perturbation}
We derive variational inequalities for minimizer to $\J_{N, d}^{\gamma, G}$. Given $U \in \A$ and $D \Subset \O,$ we say $\Phi\in H_0^1(D;\R^N)\cap C_c^0(D;\R^N)$ is a \textit{vertical perturbation in $D$} if 
$U+t\Phi \in \A$ for all small $|t| > 0$.

\begin{lemma}\label{lem:interior-equations}
Let $U\in\A$ locally minimize
$\J_{N,d}^{\gamma,G}$. Then for any $D\Subset\O$ and vertical perturbation $\Phi$ in $D$
\begin{equation}\label{eq:vertical-first-variation}
\begin{aligned}
 -K\int_D|\Phi|\,dx
 &\leq2\sum_{i=1}^N\int_D
       \nabla u_i\cdot\nabla\Phi_i\,dx+\gamma\sum_{i=1}^{N-1}w_{i}
       \int_{D\cap\{u_i>u_{i+1}\}}
       (u_i-u_{i+1})^{\gamma-1}
       (\Phi_i-\Phi_{i+1})\,dx\leq K\int_D|\Phi|\,dx.
\end{aligned}
\end{equation}

Moreover, for 
$I=\{k,\ldots,k+m-1\}\subseteq\{1,\ldots,N\}$ with \(
 \O_I=\{
 u_{k-1}>u_k\ge\cdots\ge
 u_{k+m-1}>u_{k+m}\},
\)
\begin{equation}\label{eq:vertical-first-variation-equation-perturbed}
\begin{aligned}
 -K\sqrt m\leq&
 2\Delta\left(\sum_{j\in I}u_j\right)
 -\gamma w_{k+m-1}
       (u_{k+m-1}-u_{k+m})^{\gamma-1}+\gamma w_{k-1}
       (u_{k-1}-u_k)^{\gamma-1}
 \leq K\sqrt m
 \qquad\text{weakly in }\O_I.
\end{aligned}
\end{equation}
We omit upper exterior term if $k = 1$ and the lower if $k +m-1 = N.$ 
\end{lemma}

\begin{proof}
Choose $\tau>0$ so $U\pm\tau\Phi$ are both ordered. Then
\(
 |\Phi_i-\Phi_{i+1}|
 \leq\tau^{-1}(u_i-u_{i+1})
\) for all $i.$
In particular, $\Phi_i=\Phi_{i+1}$ on
$\{u_i=u_{i+1}\}$.
Now on $\{u_i >u_{i+1}\},$ choose $0<|t|\leq\tau/2$ so mean value theorem gives
\[
\begin{aligned}
 \left|
 \frac{(u_i-u_{i+1}+t(\Phi_i-\Phi_{i+1}))^\gamma
       -(u_i-u_{i+1})^\gamma}{t}
 \right|
 &\leq
 \frac{\gamma\,2^{1-\gamma}}{\tau}
 (u_i-u_{i+1})^\gamma.
\end{aligned}
\] On $\{u_i =u_{i+1}\}$ the above is trivial, so DCT applies to the potential. Also, minimality and $G \in \mathbf{G}_K$  imply
\[
\begin{aligned}
 \J_{N,d}^{\gamma}(U+t\Phi,\mathbf w_\gamma,D)
 -\J_{N,d}^{\gamma}(U,\mathbf w_\gamma,D)\geq
 -\int_D\bigl(G(x,U+t\Phi)-G(x,U)\bigr)\,dx
 \geq-K|t|\int_D|\Phi|\,dx.
\end{aligned}
\]
Dividing and passing $|t| \to 0$ proves \eqref{eq:vertical-first-variation}. Now fix $\varphi\in C_c^\infty(\O_I)$ and take
$\Phi_j=\varphi$ for $j\in I$ and $\Phi_j=0$ otherwise. Since $|\Phi|=\sqrt m\,|\varphi|$ and only the two
exterior gaps vary, \eqref{eq:vertical-first-variation-equation-perturbed} follows. 
\end{proof}

In particular, \eqref{eq:vertical-first-variation-equation-perturbed} with $m=N$ gives
\begin{equation}\label{eq:laplace-average-bound}
     -\frac{K}{2\sqrt N}
 \leq\Delta\left(\frac1N\sum_{j=1}^N u_j\right)
 \leq\frac{K}{2\sqrt N}
 \qquad\text{weakly in }\O.
\end{equation}

When $G =0,$ we may consider the family of critical points w.r.t. vertical perturbations. Set
\[
 \A_{\mathrm{loc}}
 =\{U\in H^1_{\mathrm{loc}}(\O;\R^N)\cap C^0(\O;\R^N):
       u_1\ge\cdots\ge u_N\}.
\]
\begin{definition}\label{def:critical-points}
$U\in\A_{\mathrm{loc}}$ is vertically critical for
$\J_{N,d}^{\gamma}$ if for any $D\Subset\O$ and vertical perturbation $\Phi$ in $D$
\begin{equation*}
 \left.\frac{d}{dt}\right|_{t=0}
 \J_{N,d}^{\gamma}(U+t\Phi,\mathbf w_\gamma,D)=0.
\end{equation*}
\end{definition}
Using the first-variation \eqref{eq:vertical-first-variation}, critical points satisfy
\begin{equation}\label{eq:vertical-first-variation-equation}
\begin{aligned}
 2\Delta\left(\sum_{j\in I}u_j\right)
 &=\gamma w_{k+m-1}
       (u_{k+m-1}-u_{k+m})^{\gamma-1}-\gamma w_{k-1}
       (u_{k-1}-u_k)^{\gamma-1}
 \qquad\text{weakly in }\O_I.
\end{aligned}
\end{equation}

\subsubsection{Domain Perturbation}
We also discuss critical points w.r.t. domain perturbations for $G =0$.
\begin{definition}\label{def:gap-domain-criticality}
$U\in\A_{\mathrm{loc}}$ is domain stationary for
$\J_{N,d}^{\gamma}$ if for all $D\Subset\O$ and $\xi\in C_c^\infty(D;\R^d)$,
\[
 \left.\frac{d}{dt}\right|_{t=0}
 \J_{N,d}^{\gamma}
 \bigl(U\circ(\mathrm{Id}+t\xi)^{-1},\mathbf w_\gamma,D\bigr)=0.
\]
\end{definition}

We record the stationarity
identities used in the Weiss formula
\Cref{subsec:Weiss-monotonicity}; see \cite[Lemma~3.2]{DSS24}.

\begin{lemma}\label{lem:gap-stationarity-g0}
Let $U\in\A_{\mathrm{loc}}$ be vertically critical and
domain stationary for $\J_{N,d}^{\gamma}$ in $B_1$.
Then, 
\begin{equation}\label{eq:gap-domain-stationarity}
 0=\int_{B_1}\left[
 \bigl(|\nabla U|^2+W_N^\gamma(U)\bigr)\Div\xi
 -2\sum_{i,j=1}^d
 (\partial_iU\cdot\partial_jU)\partial_i\xi_j
 \right]dx, \qquad \forall \ \xi\in C_c^\infty(B_1;\R^d).
\end{equation}
For a.e. $r\in(0,1)$, with $\nu=x/r$ on
$\partial B_r$,
\begin{align}
 \label{eq:gap-rellich-g0}
 \int_{B_r}\left[
 (d-2)|\nabla U|^2+dW_N^\gamma(U)
 \right]dx
 &=r\int_{\partial B_r}\left(
 |\nabla U|^2+W_N^\gamma(U)-2|\partial_\nu U|^2
 \right)dS,\\
 \label{eq:gap-vertical-g0}
 \int_{B_r}\left(
 2|\nabla U|^2+\gamma W_N^\gamma(U)
 \right)dx
 &=2\int_{\partial B_r}U\cdot\partial_\nu U\,dS.
\end{align}
\end{lemma}

\begin{proof}
Fix $\xi\in C_c^\infty(B_1;\R^d)$ and
$D\Subset B_1$ containing $\supp\xi$.
Set $U_t=U\circ(\mathrm{Id}+t\xi)^{-1}$.
Changing variables and using domain stationarity gives
\eqref{eq:gap-domain-stationarity}. Taking $\xi(x)=\eta(|x|)x$ with $\eta \in C _c^\infty$, $0 \leq \eta \leq 1$ constant near $0$,
\[
\begin{aligned}
 0=&\int_{B_1}\left[
 (d-2)|\nabla U|^2+dW_N^\gamma(U)
 \right]\eta(|x|)\,dx+\int_{B_1}\left(
 |\nabla U|^2+W_N^\gamma(U)-2|\partial_\nu U|^2
 \right)|x|\eta'(|x|)\,dx.
\end{aligned}
\]
Approximating $\eta=\mathbf1_{[0,r)}$ gives \eqref{eq:gap-rellich-g0}. Vertical criticality applied with $(1+t\varphi)U$ for $\varphi\in C_c^\infty(B_1)$ gives
\[
 0=\int_{B_1}\left[
 \left(2|\nabla U|^2+\gamma W_N^\gamma(U)\right)\varphi
 +2\sum_{i=1}^d
 (U\cdot\partial_iU)\partial_i\varphi
 \right]dx.
\]
Approximating $\varphi=\mathbf1_{B_r}$ by radial cutoffs
gives \eqref{eq:gap-vertical-g0}.
\end{proof}

\subsection{Basic Operations}
Near a contact point, one would like to isolate a maximal block of membranes that collapses there and study it on $B_1$, with zero mean and value $0$ at the origin. We conduct such reduction through three operations: Rescaling, projection onto zero-mean, and local splitting;
each of which preserves minimality for an energy of
the same form, with a possibly different Lipschitz
perturbation. See \Cref{tab:basic-operations}.

\begin{table}[htbp]
\centering
\small
\renewcommand{\arraystretch}{1.5}
\setlength{\tabcolsep}{4pt}
\begin{tabular}{@{}l c c c@{}}
\toprule
Operation & Minimizer & Energy & Perturbation\\
\midrule
Rescaling
 & $U_{x_0,r}$
 & $\J_{N,d}^{\gamma,G_{x_0,r}}
    (\cdot,\mathbf w_\gamma,B_1)$
 & $G_{x_0,r}\in\mathbf G_{Kr^{2-\beta}}$
\\
Zero-mean projection
 & $\overline U$
 & $\J_{N,d}^{\gamma,\overline G}
    (\cdot,\mathbf w_\gamma,\O)$
 & $\overline G\in\mathbf G_K$
\\
Local splitting
 & $U^I$
 & $\J_{m,d}^{\gamma,G_I}
    (\cdot,\mathbf w_\gamma^I,B_r(x_0))$
 & $G_I\in\mathbf G_{K_I}$
\\
\bottomrule
\end{tabular}
\setlength{\abovecaptionskip}{2pt}
\setlength{\belowcaptionskip}{1pt}
\caption{Basic operations on minimizers.}
\label{tab:basic-operations}
\end{table}

Rescaling is described in \eqref{eq:rescale-G}.
We establish zero-mean projection and local splitting below.

\subsubsection{Zero mean}

For $U\in \A$, one may always project it onto zero-mean subspace by subtracting its average
\begin{equation*}
 \overline{U} \defeq  U-\frac1N\left(\sum_{i=1}^N u_i\right)\mathbf1 \in  \A_0:=\left\{V\in\A :
      \sum_{i=1}^N v_i=0\ \text{a.e. in }\O\right\}.
\end{equation*}

\begin{lemma}\label{lem:zero-mean}
 Let $U $ minimize $\J_{N, d}^{\gamma, G}$ with trace $\Psi \in \A$. Then $\overline{U}$ minimizes
$\J_{N, d}^{\gamma,\overline G}$ with trace \(\overline{\Psi},\) where
\begin{equation}\label{eq:centered-perturbation}
 \overline G(x,V)
 =G(x,\overline{V}+h(x)\mathbf1)-G(x,h(x)\mathbf1)
 \in\mathbf G_K, \qquad h = \frac{1}{N}\sum_{i=1}^N u_i.
\end{equation}
\end{lemma}
\begin{proof}
Take any $V\in\A$ with trace $\overline\Psi$.
Its projection $\overline V$ has the same trace,
so $\overline V+h\mathbf1$ is an admissible competitor
for $U$ with trace $\Psi$.
Since $\overline V$ and $\overline U$ have zero mean,
\[
 |\nabla(\overline V+h\mathbf1)|^2
 =|\nabla\overline V|^2+N|\nabla h|^2,
 \qquad
 |\nabla U|^2
 =|\nabla\overline U|^2+N|\nabla h|^2.
\]
Moreover, projection preserves the membrane differences and,
by \eqref{eq:centered-perturbation},
\(
 \overline G(x,\overline V)=\overline G(x,V),\
 |\nabla\overline V|^2\leq|\nabla V|^2,
\)
so canceling the common terms
$N|\nabla h|^2+G(x,h\mathbf1)$, minimality of $U$ give
\[
\begin{aligned}
 0
 &\leq
 \J_{N,d}^{\gamma,G}(\overline V+h\mathbf1)
 -\J_{N,d}^{\gamma,G}(U)=
 \J_{N,d}^{\gamma,\overline G}(\overline V)
 -\J_{N,d}^{\gamma,\overline G}(\overline U)\leq
 \J_{N,d}^{\gamma,\overline G}(V)
 -\J_{N,d}^{\gamma,\overline G}(\overline U).
\end{aligned}
\]
Therefore $\overline U$ minimizes
$\J_{N,d}^{\gamma,\overline G}$ with boundary
$\overline\Psi$ among all ordered competitors.
\end{proof}

\subsubsection{Local splitting}\label{sec:splitting-mechanism}

Now we justify that by allowing perturbations in $\mathbf{G}_K$, an interior block of membranes that strictly separates from the top and bottom blocks remains a minimizer to an energy of fewer membranes and a different perturbation. We first explain how $G \equiv 0$ fails to retain minimality under local splitting.

\begin{figure}[htbp]
\centering
\vspace{-3pt} 
\setlength{\abovecaptionskip}{1pt}
\definecolor{membraneone}{RGB}{31,119,180}
\definecolor{membranetwo}{RGB}{255,127,14}
\definecolor{membranethree}{RGB}{44,160,44}
\begin{tikzpicture}
\pgfmathsetmacro{\betamotiv}{4/3}
\pgfmathsetmacro{\amotiv}{pow(27/32,2/3)}
\pgfmathsetmacro{\gapmark}{\amotiv*pow(1.45,\betamotiv)}

\begin{axis}[
    width=0.52\linewidth, height=0.33\linewidth,
    axis lines=left,
    domain=-0.25:2.1, samples=151,
    xmin=-0.25, xmax=2.1, ymin=-1.00, ymax=1.80,
    xtick={0,1,2}, ytick={0}, xlabel={$x$},
    tick label style={font=\scriptsize}
]
\path[fill=gray!12]
    (axis cs:1,-1.00) rectangle (axis cs:2,1.80);
\draw[gray!60,densely dotted]
    (axis cs:1,-1.00) -- (axis cs:1,1.80)
    (axis cs:2,-1.00) -- (axis cs:2,1.80);

\addplot[membraneone,very thick]
    {(2/3)*\amotiv*pow(max(x,0),\betamotiv)};
\addplot[membranetwo,very thick,dashed]
    {-(1/3)*\amotiv*pow(max(x,0),\betamotiv)};
\addplot[membranethree,very thick,dotted]
    {-(1/3)*\amotiv*pow(max(x,0),\betamotiv)};

\draw[<->,semithick]
    (axis cs:1.45,{-\gapmark/3})
    -- (axis cs:1.45,{2*\gapmark/3});

\node[font=\small,text=membraneone]
    at (axis cs:1.68,1.51) {$u_1$};
\node[font=\small]
    at (axis cs:1.58,-0.84) {$u_2=u_3$};
\end{axis}
\end{tikzpicture}
\caption{The one-dimensional single-gap cone for $N=3$,
$\gamma=1/2$, and $w_1=w_2=1$.}
\label{fig:motivation-local-splitting}
\end{figure}
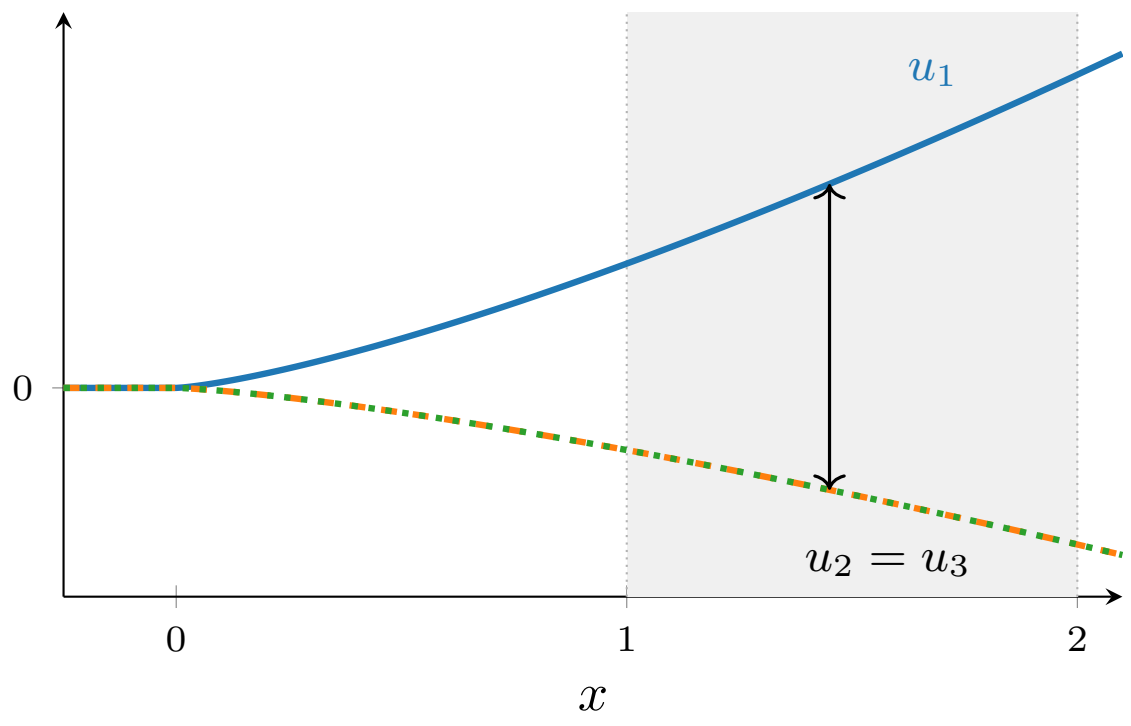
\begin{remark}
\label{rmk:separated-blocks-interaction}
By \Cref{prop:single-gap-minimizers}, the map
\[
  U(x)=a_1x_+^\beta
  \left(\frac23,-\frac13,-\frac13\right),
  \qquad
  a_1=\left(\frac{3}{2\beta^2}\right)^{1/(2-\gamma)},
\]
minimizes the three-membrane energy with $G=0$ and unit weights;
see \Cref{fig:motivation-local-splitting} for $\gamma=1/2$.

On $(1,2)$, the blocks $u_1$ and $(u_2,u_3)$ separate strictly, but neither minimizes its energy with $G = 0.$ Indeed, the model for $N =1$ is the Dirichlet energy, and for $N = 2$, $u_2 = u_3$ also enforces $(u_2, u_3)$ to minimize the Dirichlet energy; yet none of the two blocks are affine. This is because the nonlinear interaction $(u_1 - u_2)^\gamma$ contributes even under strict separation. At the endpoints, this interaction is
constant on an open gap when $\gamma=0$, while for $\gamma=1$
it splits into the linear forcing terms of the upper and lower blocks separately.
\end{remark}

Now we show how $\mathbf{G}_K$ preserves minimality under local splitting.
\begin{proposition}\label{prop:local-splitting}
    Let $U \in \A$ be a minimizer to $\J_{N, d}^{\gamma, G}(\cdot, \mathbf{w}_\gamma, \O).$ Denote  $I=\{k,\ldots,k+m-1\}\subsetneq\{1,\ldots,N\}$
as a consecutive block of length $m$, and suppose that at $x_0\in\O$, one has gap size
\[
 \kappa=\min\left\{
 u_{k-1}(x_0)-u_k(x_0),\,
 u_{k+m-1}(x_0)-u_{k+m}(x_0)
 \right\}>0.
\] 
with convention $u_0=+\infty, u_{N+1}=-\infty$. Set
\(
 U^I=(u_k,\ldots,u_{k+m-1})
\) as the interior block of membranes.

Then for any $B_R(x_0) \Subset \O$, there exists \(
 r_0=r_0\left(
d,K,\gamma,\|\mathbf{w}_\gamma\|_{\ell^\infty},\kappa,R,
 [U]_{C^{0,1/2}(B_R(x_0))}
 \right)\in(0,R/2]
\) s.t.  
\begin{equation*}
 \J_{m,d}^{\gamma,G_I}(U^I,\mathbf{w}_\gamma^I,B_r(x_0))
 =\inf\left\{
 \J_{m,d}^{\gamma,G_I}(V,\mathbf{w}_\gamma^I,B_r(x_0)):
v_1\geq\cdots\geq v_m, \ V - U^I \in H_0^1(B_r(x_0); \R^m)
 \right\}, \qquad \forall \ 0 < r \leq r_0.
\end{equation*}
where $\mathbf{w}_\gamma^I = (w_{k},\ldots,w_{k+m-2}) \in (0,\infty)^{m-1}$ and $G_I\in\mathbf G_{K_I}$ with
$K_I=K_I(K,\gamma,\|\mathbf{w}_\gamma\|_{\ell^\infty},\kappa)<\infty$.
\end{proposition}
\begin{figure}[htbp]
\centering
\vspace{-6pt} 
\setlength{\abovecaptionskip}{1pt}
\definecolor{membraneone}{RGB}{31,119,180}
\definecolor{membranetwo}{RGB}{255,127,14}
\definecolor{membranethree}{RGB}{44,160,44}
\resizebox{0.9\textwidth}{!}{%
\begin{tikzpicture}
\begin{groupplot}[
    group style={group size=2 by 1, horizontal sep=1.8cm},
    width=0.46\textwidth,
    height=0.36\textwidth,
    domain=-1.35:1.35,
    samples=201,
    xmin=-1.5, xmax=1.5,
    ymin=-1.55, ymax=1.55,
    axis x line=bottom,
    axis y line=none,
    xtick={-0.9,0,0.9},
    xticklabels={$x_0-r$,$x_0$,$x_0+r$},
    ytick=\empty,
    title style={font=\small},
    tick label style={font=\scriptsize},
    clip=false,
]

\nextgroupplot[title={Original block $U^I$}]
\fill[blue!3]
    (axis cs:-0.9,-1.4) rectangle (axis cs:0.9,1.4);
\draw[gray!50, densely dotted]
    (axis cs:-0.9,-1.4) -- (axis cs:-0.9,1.4)
    (axis cs:0.9,-1.4) -- (axis cs:0.9,1.4);
\addplot[gray, thick] {0.08*x+1+0.08*x^2};
\addplot[gray, thick] {0.08*x-1-0.08*x^2};
\addplot[membraneone, very thick] {0.08*x+0.28*x^2};
\addplot[membranetwo, very thick, dashed] {0.08*x};
\addplot[membranethree, very thick, dotted] {0.08*x-0.28*x^2};

\node[gray, font=\small] at (axis cs:0.95,1.35) {$u_{k-1}$};
\node[gray, font=\small] at (axis cs:0.95,-1.25) {$u_{k+m}$};
\node[membraneone, font=\small] at (axis cs:-1.15,0.48) {$u_k$};
\node[membranetwo, font=\small] at (axis cs:-1.15,0.05) {$u_{k+1}$};
\node[membranethree, font=\small] at (axis cs:-1.15,-0.64) {$u_{k+2}$};

\draw[<->, thick] (axis cs:0,0.03) -- (axis cs:0,0.97)
    node[midway,right=3pt,font=\small] {$\geq\kappa/2$};
\draw[<->, thick] (axis cs:0,-0.97) -- (axis cs:0,-0.03)
    node[midway,right=3pt,font=\small] {$\geq\kappa/2$};

\nextgroupplot[title={Replacement by $V^*$}]
\fill[blue!3]
    (axis cs:-0.9,-1.4) rectangle (axis cs:0.9,1.4);
\draw[gray!50, densely dotted]
    (axis cs:-0.9,-1.4) -- (axis cs:-0.9,1.4)
    (axis cs:0.9,-1.4) -- (axis cs:0.9,1.4);

\addplot[gray, thick] {0.08*x+1+0.08*x^2};
\addplot[gray, thick] {0.08*x-1-0.08*x^2};

\addplot[membraneone, very thick]
    {0.08*x+0.28*x^2+0.62*max(1-(x/0.9)^2,0)};
\addplot[membranetwo, very thick, dashed]
    {0.08*x+0.08*max(1-(x/0.9)^2,0)};
\addplot[membranethree, very thick, dotted]
    {0.08*x-0.28*x^2-0.52*max(1-(x/0.9)^2,0)};

\node[gray, font=\small] at (axis cs:0.95,1.35) {$u_{k-1}$};
\node[gray, font=\small] at (axis cs:0.95,-1.25) {$u_{k+m}$};
\node[membraneone, font=\small] at (axis cs:-0.45,0.63) {$v_1^*$};
\node[membranetwo, font=\small] at (axis cs:-0.45,0.16) {$v_2^*$};
\node[membranethree, font=\small] at (axis cs:-0.45,-0.65) {$v_3^*$};

\addplot[only marks, mark=*, mark size=1.5pt, membraneone]
    coordinates {(-0.9,0.1548) (0.9,0.2988)};
\addplot[only marks, mark=*, mark size=1.5pt, membranetwo]
    coordinates {(-0.9,-0.072) (0.9,0.072)};
\addplot[only marks, mark=*, mark size=1.5pt, membranethree]
    coordinates {(-0.9,-0.2988) (0.9,-0.1548)};

\draw[<->, thick] (axis cs:0,0.65) -- (axis cs:0,0.97)
    node[midway,right=3pt,font=\small] {$\geq\kappa/4$};
\draw[<->, thick] (axis cs:0,-0.97) -- (axis cs:0,-0.55)
    node[midway,right=3pt,font=\small] {$\geq\kappa/4$};

\end{groupplot}
\end{tikzpicture}%
}
\caption{Local splitting in \Cref{prop:local-splitting} along a line through $x_0$ for $m=3$.}
\label{fig:local-splitting}
\end{figure}
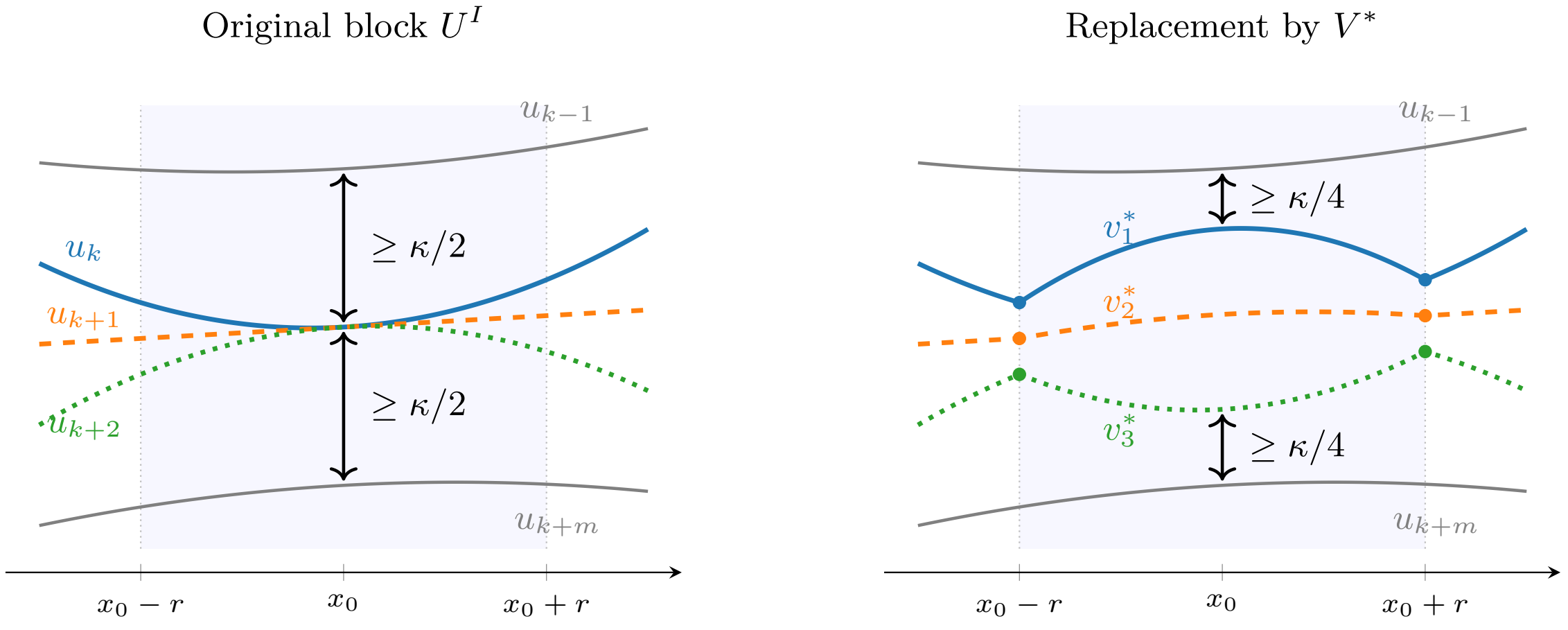
\begin{proof}
 \emph{Step 1.}  We begin with constructing the new perturbation $G_I.$  Define 
\begin{equation}\label{eq:g-kappa}
 g_\kappa(s)=
 \begin{cases}
 s^\gamma,&s\geq\kappa/4,\\
 \kappa/4-s+(\kappa/4)^\gamma,&s<\kappa/4,
 \end{cases}
 \qquad \seminorm{g_\kappa}_{C^{0, 1}}\leq 
 L_\kappa=\max\{1,\gamma(\kappa/4)^{\gamma-1}\}.
\end{equation}
For any $x \in B_R(x_0)$, set
\[
\begin{aligned}
 G_I(x,V)
 =&G(x,u_1, \cdots, u_{k-1}, V, u_{k+m}, \cdots, u_N)-G(x, u_1, \cdots, u_{k-1}, 0, \cdots, 0, u_{k+m}, \cdots, u_{N})\\
 &+w_{k-1}
 \bigl[g_\kappa(u_{k-1}-v_1)-g_\kappa(u_{k-1})\bigr]+w_{k+m-1}
 \bigl[g_\kappa(v_m-u_{k+m})-g_\kappa(-u_{k+m})\bigr].
\end{aligned}
\] We omit upper exterior term if $k = 1$ and the lower if $k +m-1 = N.$ Then $G_I\in\mathbf G_{K_I}$ due to
\[
 G_I(x,0)=0,\qquad
 |G_I(x,V)-G_I(x,W)|\le K_I|V-W|,
 \qquad
K_I=K+2\|\mathbf{w}_\gamma\|_{\ell^\infty}L_\kappa.
\]

\emph{Step 2.} By \Cref{cor:holder-minimizers}, set
\(
 M:=[U]_{C^{0,1/2}(B_R(x_0))}<\infty.
\) Choose \(r_0=\min\left\{
 \frac R2,\,
 \frac{\kappa^2}{16M^2},\,
 \sqrt{\frac{d\kappa}{2K_I}}
 \right\},\) so that
\begin{equation}\label{eq:g-kappa-radius-requirement}
Mr_0^{1/2}\leq\frac{\kappa}{4},
 \qquad
 \frac{K_Ir_0^2}{2d}\leq\frac{\kappa}{4}.
\end{equation}
The first Hölder bound in \eqref{eq:g-kappa-radius-requirement} and definition of $\kappa$ gives, throughout $B_{r_0}(x_0)$
\begin{equation}\label{eq:block-separation}
 \begin{aligned}
 u_{k+m-1}(x_0)-\frac{\kappa}{4}
 \leq u_{k+m-1}\leq u_k
 &\leq u_k(x_0)+\frac{\kappa}{4},\\
 u_{k-1}
 \geq u_{k-1}(x_0) - M r_0^{1/2}&\geq u_k(x_0)+\frac{3\kappa}{4}
 &&k>1,\\
 u_{k+m}
 \leq u_{k+m}(x_0) + M r_0^{1/2}&\leq u_{k+m-1}(x_0)-\frac{3\kappa}{4}
 &&k+m\leq N.
 \end{aligned}
\end{equation}
In particular, the exterior gaps of $U$,  $u_{k-1} - u_k$ and $u_{k+m-1}- u_{k+m}$, are at least
$\kappa/2$ in $B_{r_0}(x_0)$. 

Fix $0<r\leq r_0$. By \Cref{prop:existence}, there
exists a minimizer $V^*=(v_1^*,\ldots,v_m^*)$ of
$\J_{m,d}^{\gamma,G_I}(\cdot,\mathbf w_\gamma^I,B_r(x_0))$ among ordered membranes with boundary $U^I$.
By \Cref{lem:linfty-bound},
\[
 \Delta v_1^*\geq-K_I,
 \qquad
 \Delta v_m^*\leq K_I.
\]Applying the weak maximum principle, using
\eqref{eq:g-kappa-radius-requirement} and the first line in \eqref{eq:block-separation},
we obtain
\[
 \begin{aligned}
 v_1^*(x)
 &\leq\sup_{\partial B_r(x_0)}u_k
      +\frac{K_I}{2d}(r^2-|x-x_0|^2)\leq u_k(x_0)+\frac{\kappa}{4}
      +\frac{\kappa}{4}
 =u_k(x_0)+\frac{\kappa}{2},\\
 v_m^*(x)
 &\geq\inf_{\partial B_r(x_0)}u_{k+m-1}
      -\frac{K_I}{2d}(r^2-|x-x_0|^2)\geq u_{k+m-1}(x_0)-\frac{\kappa}{4}
      -\frac{\kappa}{4}
 =u_{k+m-1}(x_0)-\frac{\kappa}{2}.
 \end{aligned}
\]
Consequently, using the second and third lines in \eqref{eq:block-separation} (see \Cref{fig:local-splitting}), 
\[
 u_{k-1}-v_1^*\geq\frac{\kappa}{4}
 \qquad k>1,
 \qquad
 v_m^*-u_{k+m}\geq\frac{\kappa}{4}
 \qquad k+m\leq N.
\]
Now $(u_1, \ldots, u_{k-1},  V^*, u_{k+m}, \ldots, u_N)$ extended by $U$
outside $B_r(x_0)$ is an admissible competitor for $U$ among $\A$ in $\O$. Moreover, the exterior gaps of them both are at
least $\kappa/4$ in $B_r(x_0)$, so $g_\kappa(s)=s^\gamma$ there. Hence, 
\[
\begin{aligned}
 0
 &\leq
 \J_{N,d}^{\gamma,G}(u_1, \ldots, u_{k-1},  V^*, u_{k+m}, \ldots, u_N,\mathbf w_\gamma,\O)
 -\J_{N,d}^{\gamma,G}(U,\mathbf w_\gamma,\O)\\
 &=\int_{B_r(x_0)}\Biggl[
 \sum_{j=1}^{m}
 \left(|\nabla v_j^*|^2-|\nabla u_{k+j-1}|^2\right)
 +\sum_{j=1}^{m-1}w_{k+j-1}
 \left((v_j^*-v_{j+1}^*)^\gamma
       -(u_{k+j-1}-u_{k+j})^\gamma\right)\\
 &\qquad
 +w_{k-1}
 \left((u_{k-1}-v_1^*)^\gamma
       -(u_{k-1}-u_k)^\gamma\right)
 +w_{k+m-1}
 \left((v_m^*-u_{k+m})^\gamma
       -(u_{k+m-1}-u_{k+m})^\gamma\right)\\
 &\qquad
 +G(x,u_1, \ldots, u_{k-1},  V^*, u_{k+m}, \ldots, u_N)-G(x,U)
 \Biggr]\,dx=\J_{m,d}^{\gamma,G_I}(V^*,\mathbf w_\gamma^I,B_r(x_0))
 -\J_{m,d}^{\gamma,G_I}(U^I,\mathbf w_\gamma^I,B_r(x_0))
 \leq0.
\end{aligned}
\]
\end{proof}

\begin{remark}\label{rmk:reduction-from-maximal-contact-block}
We now reduce the problem near a contact point to a setting with
zero-mean minimizer in $B_1$ with $0$ at the origin. To be precise, let $U$ be a minimizer and $x_0\in\O$. We call
$I=\{k,\ldots,k+m-1\}$ a
\textit{maximal contact block at $x_0$ of size $m$} if
\[
 u_{k-1}(x_0)>
 u_k(x_0)=\cdots=u_{k+m-1}(x_0)>
 u_{k+m}(x_0).
\]
Then using \eqref{eq:rescale-G}, \Cref{lem:zero-mean,prop:local-splitting}, the composition of all three operations
\[
 \overline U_{x_0, r}^I(x)
 =r^{-\beta}\left(U^I(x_0+rx)- (\frac1m\sum_{j=k}^{k+m-1}u_j)(x_0+rx)\mathbf1\right)
\]
minimizes \(
 \J_{m,d}^{\gamma,\overline G_{I,x_0, r}}
 (\cdot,\mathbf w_\gamma^I,B_1),\) with \( 
 \overline G_{I, x_0, r}\in\mathbf G_{K_Ir^{2-\beta}},
\)
for $r > 0$ small, and
\(
 \sum_{j=1}^m(\overline U_{x_0, r}^I)_j=0, \
 \overline U_{x_0, r}^I(0)=0.
\)
\end{remark}

\section{Optimal regularity}\label{sec:optimal-regularity}

In this section, $N \geq2$, $\gamma \in (0, 1),$ $\beta=2/(2-\gamma)\in(1,2)$ as in \eqref{eq:homogeneity}, $G \in \mathbf{G}_K$ for $K > 0$, $\mathbf{w}_\gamma \in (0, \infty)^{N-1}$, unless otherwise specified. We prove optimal regularity for a minimizer to $\J_{N, d}^{\gamma, G}.$

\begin{theorem}\label{thm:optimal-regularity}
Let $U \in \A$ minimize $\J_{N,d}^{\gamma,G}$ in $B_1$, then
\begin{equation}\label{eq:opt-reg-singular}
 \|U\|_{C^{1,\beta-1}(\overline{B_{1/2}})}
 \le C(d,N,\gamma,K,\|\mathbf{w}_\gamma\|_{\ell^\infty})\left(1+\|U\|_{L^2(B_1)}\right).
\end{equation}
\end{theorem}

We first claim it suffices to work with $U \in \A_0$ only. Indeed, by \Cref{lem:zero-mean}, removing the mean remains minimizer among $\A$ with Lipschitz perturbation of the same $K$, and \eqref{eq:laplace-average-bound} gives for every $0 < \alpha < 1$
\begin{equation}\label{eq:average-local-estimate}
 \|\frac{1}{N}\sum_{i = 1}^N u_i\|_{C^{1,\alpha}(\overline{B_{3/4}})}
 \le C(d, N, K, \alpha)(1+\|U\|_{L^2(B_1)}).
\end{equation}

\subsection{Dichotomy}\label{sec:dichotomy}
The following dichotomy argument adopts from \cite[Theorem 1.1]{DS20}. The idea is to say, for a minimizer $U$, if
\[
    a=\left(\fint_{B_1}|U|^2\,dx\right)^{1/2}
\]
is large, then \(a^{-1}U\) minimizes small perturbation of the Dirichlet
energy, and \(a^{-1}U\) is close to its harmonic replacement $H$. Either the center $|H(0)|$ is small, which gives decay in $L^2$ average of $|U|$ at the
\(\beta\)-scale, or $|H(0)|$ is bounded away from zero, which
forces \(|U|\) to stay positive at possibly another scale. See also \cite{Ma26} for a comprehensive treatment for the scalar case Alt-Phillips.

\begin{proposition}[Dichotomy]\label{prop:dichotomy-singular}
For any $\xi\in(0,1/2]$, there exists \( a_*\ge1\) depending only on \(d,N,\gamma,K,\|\mathbf{w}_\gamma\|_{\ell^\infty},\xi\) s.t. if $U\in\A_0$ minimize
$\J_{N,d}^{\gamma,G}$ in $B_1$ with
\(a=\left(\fint_{B_1}|U|^2\,dx\right)^{1/2}\ge a_*,
\)
then either
\begin{equation}\label{eq:dichotomy-decay}
 \xi^{-\beta}
 \left(\fint_{B_\xi}|U|^2\,dx\right)^{1/2}
 \le\frac a2,
\end{equation}
or for \(\eta=\eta(d,N)\in(0,1/4), \ c_0=c(d,N)\xi^\beta>0,  \ C_0=C_0(d,N)\ge1,\)
\begin{equation}\label{eq:dichotomy-positive-ball}
 c_0a\le |U(x)|\le C_0a
 \qquad\text{for every }x\in B_\eta.
\end{equation}
\end{proposition}
\begin{proof}
Throughout the proof, write $C = C(d, N, \gamma, K, \|\mathbf{w}_\gamma\|_{\ell^\infty}) > 0$ and $C_1 = C_1(d, N)> 0$ as constants that change from line to line, independent of $a$ and $\xi$. 

\emph{Step 1.} Set $\widetilde U=U/a$. Dividing the energy by $a^2$, $\widetilde U$ minimizes
\[
 \int_{B_1}|\nabla V|^2+\delta F_a(x,V)\,dx,
 \qquad
 \delta=a^{-1},\qquad
 F_a(x,V)=a^{\gamma-1}W_N^\gamma(V,\mathbf{w}_\gamma)
          +a^{-1}G(x,aV),
\]
among competitors in the closed convex cone
\(
 \mathcal C_0
 =\{V\in\R^N:v_1\ge\cdots\ge v_N,\ \textstyle\sum_i v_i=0\}
\)
with same boundary, and that \(
 |F_a(x,V)|\le C(1+|V|^2).\) Let $H$ be the harmonic replacement of $\widetilde U$
in $B_1$. Convexity gives $H\in\mathcal C_0$. From minimality, we have
\[
\begin{aligned}
 \int_{B_1}|\nabla(\widetilde U-H)|^2\,dx
 &\le\delta\int_{B_1}
       \bigl(F_a(x,H)-F_a(x,\widetilde U)\bigr)\,dx\le C\delta\int_{B_1}(1+|\widetilde U|^2+|H|^2)\,dx
 \le C\delta\left(
       1+\int_{B_1}|\nabla(\widetilde U-H)|^2\,dx
       \right).
\end{aligned}
\]
Taking $a_*$ large and absorbing gives
\begin{equation}\label{eq:dichotomy-harmonic-close}
 \int_{B_1}|\nabla(\widetilde U-H)|^2\,dx
 +\int_{B_1}|\widetilde U-H|^2\,dx
 \le C\delta.
\end{equation}
Increasing $a_*$ if necessary and using
$\fint_{B_1}|\widetilde U|^2\,dx=1$, we obtain
\[
 |H(0)|^2
 \le\fint_{B_1}|H|^2\,dx
 \le2+2\fint_{B_1}|H-\widetilde U|^2\,dx
 \le C_1 .
\]
Since $H_1,-H_N\ge0$ and
$H_N\le H_j\le H_1$, Harnack's inequality and
harmonic estimates give
\[
 \|H\|_{L^\infty(B_{3/4})}
 +\|\nabla H\|_{L^\infty(B_{1/2})}
 \le C_1|H(0)|.
\]

\emph{Step 2.} First suppose $|H(0)|\le2c_0$ for $c_0$ small to choose. Then
\[
 \fint_{B_\xi}|H-H(0)|^2\,dx
 \le\xi^2\|\nabla H\|_{L^\infty(B_\xi)}^2
 \le C_1|H(0)|^2.
\]
Together with \eqref{eq:dichotomy-harmonic-close}, this gives
\[
 \begin{aligned}
 \xi^{-2\beta}\fint_{B_\xi}|\widetilde U|^2\,dx
 &\le 3\xi^{-2\beta}\left(
       \fint_{B_\xi}|\widetilde U-H|^2
       +\fint_{B_\xi}|H-H(0)|^2
       +|H(0)|^2\right)\le C\xi^{-d-2\beta}\delta
       +C_1c_0^2\xi^{-2\beta}.
 \end{aligned}
\]
Choose $c=c(d,N)>0$ small and take $a_*$ large, so that by setting $c_0=c\xi^\beta$,
\[
 C_1c_0^2\xi^{-2\beta}\le\frac18,\qquad
 C\xi^{-d-2\beta}a_*^{-1}\le\frac18.
\]
This proves \eqref{eq:dichotomy-decay}.

\emph{Step 3.} Now suppose $|H(0)|>2c_0$. Choose $\eta=\eta(d,N)\in(0,1/4)$ s.t.
\[
 \sup_{B_\eta}|H-H(0)|
 \le\eta\|\nabla H\|_{L^\infty(B_{1/2})}
 \le C_1\eta|H(0)|
 \le\frac14|H(0)|.
\]
By \Cref{prop:perturbative-holder} with $\alpha=1/2$
and harmonic estimates,
\[
 \|\widetilde U-H\|_{C^{0,1/2}(B_{1/2})}\le C.
\]
For $x\in B_{1/4}$ and $0<r\le1/4$, the H\"older bound gives \(|(\tilde{U} - H)(x)| \leq |(\tilde{U} - H)(y)| + Cr^{1/2}\) for $y \in B_r(x)$. Averaging in $y$ and using \eqref{eq:dichotomy-harmonic-close} then yield
\[
\begin{aligned}
 |(\widetilde U-H)(x)|
 &\le \left(\fint_{B_r(x)}
       |\widetilde U-H|^2\,dy\right)^{1/2}
       +Cr^{1/2}\le C\left(r^{-d/2}\delta^{1/2}+r^{1/2}\right).
\end{aligned}
\]
Choosing $r=\delta^{1/(d+1)}\le1/4$ balances the two terms, and gives a uniform bound 
\[
 \|\widetilde U-H\|_{L^\infty(B_{1/4})}
 \le C\delta^{1/(2d+2)}.
\]
Increase $a_*$ so that $Ca_*^{-1/(2d+2)}\le c_0/2$.
For every $x\in B_\eta$, we have
\[
 c_0
 \le\frac34|H(0)|-\frac{c_0}{2}
 \le|\widetilde U(x)|
 \le\frac54|H(0)|+\frac{c_0}{2}
 \le\frac32|H(0)|
 \le C_0,
\]
where $C_0=C_0(d,N)$. Multiplying by $a$ proves
\eqref{eq:dichotomy-positive-ball}.
\end{proof}

Using the Dichotomy, one obtain the key Harnack Inequality and $\beta$-growth rate from the free boundary.
\begin{proposition}
\label{prop:harnack-growth}
There are \(
 c_H=c_H(d,N,\gamma)>0,\) and \( C_H\ge1, r_H\in(0,\frac{1}{8})
\) depending only on 
\(d,N,\gamma,K,\|\mathbf{w}_\gamma\|_{\ell^\infty}\)
s.t. if $U\in\A_0$ minimize
$\J_{N,d}^{\gamma,G}$ in $B_1$ with
\(
 \left(\fint_{B_1}|U|^2\,dx\right)^{1/2}\le1,
\) then 
\begin{equation}\label{eq:harnack-optimal-reg}
 c_H|U(0)|\le |U(x)|\le C_H|U(0)|
 \qquad \forall \ |x|\le r_H|U(0)|^{1/\beta}.
\end{equation}
If moreover, $U(0)=0$, one has
\begin{equation}\label{eq:l2-beta-growth}
 \left(\fint_{B_r}|U|^2\,dx\right)^{1/2}
 \le C_Hr^\beta,\qquad 0<r<1/2.
\end{equation}
\end{proposition}
\begin{proof}
Take $\xi=1/2$, and fix the constants $a_*,c_0,C_0,\eta$ from \Cref{prop:dichotomy-singular}. For $k\ge0$, set
\[
 r_k=2^{-k},\qquad
 a_k=r_k^{-\beta}
     \left(\fint_{B_{r_k}}|U|^2\,dx\right)^{1/2}.
\]
The dichotomy  \Cref{prop:dichotomy-singular} applies whenever $a_k\ge a_*$. We fix \(C_{d,\gamma}=2^{\beta+d/2}\) so that
\[
 a_{k+1}
 =2^\beta r_k^{-\beta}
   \left(\fint_{B_{r_k/2}}|U|^2\,dx\right)^{1/2}
 \le C_{d,\gamma} a_k, \qquad \forall \ k.
\]

\emph{Step 1.} Consider the set of indices for which the following bound holds
\begin{equation}\label{eq:dyadic-stopping-bound}
 a_k\le C_{d,\gamma} a_*.
\end{equation}
It holds at $k=0$, since $a_0\le1\le a_*$.
If it holds for every $k$, then for $r_{j+1}<r\le r_j$,
\[
 r^{-\beta}
 \left(\fint_{B_r}|U|^2\,dx\right)^{1/2}
 \le C_{d,\gamma} a_j
 \le C_{d,\gamma}^2a_*.
\]
In this case, continuity implies $U(0)=0$.

\emph{Step 2.} 
Otherwise, let $k+1$ be the first index where
\eqref{eq:dyadic-stopping-bound} fails.
Then $a_k>a_*$, otherwise
$a_{k+1}\le C_{d,\gamma} a_*$. Now \Cref{prop:dichotomy-singular} applies.
Note the decay alternative \eqref{eq:dichotomy-decay}
is also impossible, otherwise
\(
 a_{k+1}\le\frac12a_k\le\frac12C_{d,\gamma} a_*.
\)
Hence \eqref{eq:dichotomy-positive-ball}, applied to
$U_{0,r_k}$, gives
\begin{equation}\label{eq:stopping-scale-nonzero}
 c_0r_k^\beta a_k\le |U(x)|\le C_0r_k^\beta a_k
 \qquad\text{for }x\in B_{\eta r_k}.
\end{equation}
In particular, $U(0)\ne0$, and comparison with the bounds at
$x=0$ yields
\[
 \frac{c_0}{C_0}|U(0)|
 \le |U(x)|\le\frac{C_0}{c_0}|U(0)|
 \qquad\text{for }x\in B_{\eta r_k}.
\]
Since $a_k\le C_{d,\gamma} a_*$, combining with \eqref{eq:stopping-scale-nonzero} we also have an estimate of $r_k$ from below
\begin{equation}\label{eq:stopping-estimate-rk-below}
    |U(0)|\le C_0C_{d,\gamma} a_*r_k^\beta,
 \qquad
 r_k\ge
 \left(\frac{|U(0)|}{C_0C_{d,\gamma} a_*}\right)^{1/\beta}.
\end{equation}

Set
\[
 c_H=\frac{c_0}{C_0},\qquad
 C_H=\max\left\{\frac{C_0}{c_0},\,C_{d,\gamma}^2a_*\right\},
 \qquad
 r_H=\eta
      (C_0C_{d,\gamma} a_*)^{-1/\beta}.
\]
If $U(0)\ne0$, the stopping index exists, then the above choice and estimate \eqref{eq:stopping-estimate-rk-below} give
\[
 r_H|U(0)|^{1/\beta}\le \eta r_k,
\]
which proves \eqref{eq:harnack-optimal-reg}.
If $U(0)=0$, no stopping index exists, so
\eqref{eq:l2-beta-growth} follows from the bound in the
first case.
\end{proof}

Under the assumption in \Cref{prop:harnack-growth}, increasing $a_*$ and adjusting $r_H,C_H$ depending on $d, \gamma$, gives
\eqref{eq:harnack-optimal-reg} uniformly in
$B_{3/4}$. Now let $x_0\in B_{1/2}$ and $x\in B_r(x_0)$ for $0<r<1/4$.
If $|x-x_0|\le r_H|U(x)|^{1/\beta}$, \eqref{eq:harnack-optimal-reg} centered at $x$ gives $c_H|U(x)|\le|U(x_0)|$.
Otherwise, $|U(x)|<r_H^{-\beta}|x-x_0|^\beta$.
Increasing $C_H$ therefore gives
\begin{equation}\label{eq:general-beta-growth}
 \sup_{B_r(x_0)}|U|
 \le C_H\bigl(|U(x_0)|+r^\beta\bigr),
 \qquad x_0\in B_{1/2},\quad 0<r<1/4.
\end{equation}
In particular, since $U\in\A_0$, 
\begin{equation}\label{eq:pointwise-beta-growth}
 |U(x)|\le C_H|x-x_0|^\beta,
 \qquad x_0\in\bigcap_{i=1}^{N-1}\{u_i=u_{i+1}\}\cap B_{1/2},\quad
 |x-x_0|<1/4.
\end{equation}

\subsection{Induction on the number of membranes}
\label{sec:induction-singular}

The base case is $N=1$, for which $W_1^\gamma=0$.
Using two-sided variations and the Lipschitz bound on $G$,
a minimizer $u$ satisfies
\[
 \left|2\int_{B_1}\nabla u\cdot\nabla\varphi\,dx\right|
 \le K\int_{B_1}|\varphi|\,dx
 \qquad\text{for every }\varphi\in C_c^\infty(B_1).
\]
Thus $\|\Delta u\|_{L^\infty(B_1)}\le K/2$, and 
\[
 \|u\|_{C^{1,\beta-1}(\overline{B_{1/2}})}
 \le C(d,\gamma)\left(\|u\|_{L^2(B_1)}+K\right),
\]
which proves \eqref{eq:opt-reg-singular} for $N=1$.

\begin{proof}[Proof of \Cref{thm:optimal-regularity}]
Throughout the proof,
$C=C(d,N,\gamma,K,\|\mathbf{w}_\gamma\|_{\ell^\infty}) > 0$
denotes constant that may change from line to line. Assume \Cref{thm:optimal-regularity} for every membrane
number between $1$ and $N-1$.
By \Cref{lem:zero-mean} and \eqref{eq:average-local-estimate}, assume for zero-mean minimizer $U \in \A_0$. Also assume $\fint_{B_1}|U|^2 \leq 1$ by rescaling.

\emph{Step 1.}
By Campanato's characterization \cite{Campanato1964Famiglia}, it suffices to show
\begin{equation}\label{eq:campanato-goal-singular}
 \inf_{P\ \mathrm{affine}}
 \left(\fint_{B_r(x_0)}|U-P|^2\,dx\right)^{1/2}
 \le Cr^\beta,
 \qquad x_0\in B_{1/2},\quad 0<r<1/4.
\end{equation}
If $U(x_0)=0$, take $P=0$ and use
\eqref{eq:pointwise-beta-growth}.

\emph{Step 2.} Now suppose $U(x_0)\ne0$, and set
\[
 R=c|U(x_0)|^{1/\beta},
 \qquad
 U_{x_0,R}(x)=R^{-\beta}U(x_0+Rx).
\]
By \Cref{cor:holder-minimizers}, $|U(x_0)|\le C$.
We may therefore choose a fixed
$c=c(d,N,\gamma,K,\|w_\gamma\|_{\ell^\infty})>0$
sufficiently small that $4c\le r_H$ and $4R<1/4$, so one may apply Harnack inequality \eqref{eq:harnack-optimal-reg} and obtain
\[
 C^{-1}\le |U_{x_0,R}|\le C
 \qquad\text{in }B_4.
\]
Since $U_{x_0,R}$ is ordered and zero-mean, there exists
$\ell\in\{1,\ldots,N-1\}$ such that
\[
\begin{aligned}
 \kappa
 &\defeq
 (U_{x_0,R})_\ell(0)-(U_{x_0,R})_{\ell+1}(0)\geq
 \frac{(U_{x_0,R})_1(0)-(U_{x_0,R})_N(0)}{N-1}
 \geq\frac{|U_{x_0,R}(0)|}{\sqrt N(N-1)}
 \geq C^{-1}.
\end{aligned}
\]
Then
\(
 L_\kappa
 =\max\{1,\gamma(\kappa/4)^{\gamma-1}\}\leq C
\) in \eqref{eq:g-kappa}.
Also, \Cref{cor:holder-minimizers} again gives \(
    [U_{x_0,R}]_{C^{0,1/2}(B_1)}\le C.\)
One may then choose
\(
 \sigma=\sigma\left(
 d,N,\gamma,K,\|\mathbf w_\gamma\|_{\ell^\infty}
 \right)\in(0,1/4)
\) independent of $x_0$ and $R$, such that
\[
 [U_{x_0,R}]_{C^{0,1/2}(B_1)}\sigma^{1/2}
 \leq\frac{\kappa}{4},
 \qquad
 \frac{
 \left(K+2\|\mathbf w_\gamma\|_{\ell^\infty}L_\kappa\right)
 \sigma^2}{2d}
 \leq\frac{\kappa}{4}.
\]
These are precisely requirements \eqref{eq:g-kappa-radius-requirement} for \(U_{x_0,R},\) so one may apply \Cref{prop:local-splitting} to each of the blocks
\[
 I=\{1,\ldots,\ell\},
 \qquad
 I=\{\ell+1,\ldots,N\}.
\]
For either block, its only exterior gap at the origin
has size $\kappa$. Hence $(U_{x_0,R})^I$ minimizes
\(
 \J_{m,d}^{\gamma,G_I}
 (\cdot,\mathbf w_\gamma^I,B_\sigma)
\)
among ordered maps with \(m=|I|<N,\)
\(
 \|\mathbf w_\gamma^I\|_{\ell^\infty}
 \leq\|\mathbf w_\gamma\|_{\ell^\infty},\) and \(
 K_I\leq
 K+2\|\mathbf w_\gamma\|_{\ell^\infty}L_\kappa
 \leq C.
\)
Rescale both to $B_1$ and apply the induction hypothesis for both blocks, then scaling back gives
\[
 \|U_{x_0,R}\|_{C^{1,\beta-1}
                    (\overline{B_{\sigma/2}})}\leq C\left(1 + \norm{U_{x_0, R}}_{L^2(B_\sigma)}\right)
 \leq C.
\]

Further scaling back, for $0<r<\sigma R/2$,
\[
\begin{aligned}
 \sup_{x\in B_r(x_0)}
 \left|U(x)-U(x_0)
       -R^{\beta-1}\nabla U_{x_0,R}(0)(x-x_0)\right|\le CR^\beta(r/R)^\beta
 =Cr^\beta.
\end{aligned}
\]
For $\sigma R/2\le r<1/4$, we have
\(
 |U(x_0)|=(R/c)^\beta\le Cr^\beta,
\)
and \eqref{eq:general-beta-growth} gives
\[
 \sup_{B_r(x_0)}|U|
 \le C\bigl(|U(x_0)|+r^\beta\bigr)
 \le Cr^\beta.
\]
These complete
\eqref{eq:campanato-goal-singular}, which gives
\(
 \|U\|_{C^{1,\beta-1}(\overline{B_{1/2}})}\le C.
\)
Adding back the mean using
\eqref{eq:average-local-estimate} proves
\eqref{eq:opt-reg-singular}.
\end{proof}

\section{Non-degeneracy and gluing phenomenon}\label{sec:nondegenracy}

In this section, $N \geq 2$, $\gamma \in (0, 1)$, $\beta =2/(2 - \gamma)\in (1, 2)$ as in \eqref{eq:homogeneity}, $G \in \mathbf{G}_K$ for $K > 0$, \(\mathbf{w}_\gamma\in (0, \infty)^{N-1}\) unless otherwise stated. We first discuss non-degeneracy for $U\in \A_0$ that only require positive weights. 

\begin{proposition}\label{prop:gap-nondegeneracy}
There are
$c_0=c_0(d,N,\gamma,\min_i w_{i})>0$
and $r_0=r_0(K,\gamma)\in(0,\frac{1}{2})$ such that if $U\in\A_0$ minimize
$\J_{N,d}^{\gamma,G}$ in $B_1$, and
$0\in\overline{\{|U|>0\}}$, then
\begin{equation}\label{eq:gap-nondegeneracy-y}
 \sup_{B_r}u_1\ge c_0r^\beta,
 \qquad 0<r<r_0.
\end{equation}
\end{proposition}

\begin{proof}
First assume $K\le1$ and $\sup_{B_1}u_1\le\ep_0$ for 
$\ep_0>0$ small to choose. Ordering and zero mean give
\[
 0\le u_1,\qquad
 u_1-u_N\le Nu_1\le N\ep_0,\qquad
 \{u_1=0\}=\{ U=0\}.
\]
At any point of $\{u_1>0\}$, choose $1\le i<N$ s.t.
$u_1=\cdots=u_i>u_{i+1}$ there. By continuity, $u_i>u_{i+1}$ near that point, so applying \eqref{eq:vertical-first-variation-equation-perturbed}
to the first $i$ membranes gives
\[
 \Delta\left(\frac{u_1+\cdots+u_i}{i}\right)
 \ge
 \frac{\gamma w_{i}}{2i}
 (u_i-u_{i+1})^{\gamma-1}
 -\frac{K}{2\sqrt i} \qquad \text{weakly in $B_1 \cap \{u_i > u_{i+1}\}$}.
\]
Every test function touching $u_1$ from above also touches this average from above. Repeat at each point so
\[
 \Delta u_1
 \ge
 \frac{\gamma N^{\gamma-1}\min_i w_{i}}{2(N-1)}
 \ep_0^{\gamma-1}-\frac12
 >8d\ep_0
 \qquad\text{in }B_1\cap\{u_1>0\} \text{ in the viscosity sense},
\]
after choosing $\ep_0> 0$ small. Set $\psi(x)=\ep_0(2|x|-1)_+^2,$ so
$\Delta\psi\le8d\ep_0$ in $B_1$,
$\psi=\ep_0$ on $\partial B_1$, and $\psi=0$ in $B_{1/2}$. Now maximum principle on
$\{u_1>\psi\}\subset\{u_1>0\}$ gives $u_1\le\psi$, thus $u_1=0$ in $B_{1/2}$. In general, set
\(
 c_0=\ep_0,\
 r_0=\frac14(1+K)^{-1/(2-\beta)}.
\)
If \eqref{eq:gap-nondegeneracy-y} fails for some $0<r<r_0$, $U_{0,r}(x)=r^{-\beta}U(rx)$ minimizes a rescaled energy
with 
$\mathbf G_{Kr^{2-\beta}}\subset\mathbf G_1$
and 
\(
 \sup_{B_1}(U_{0,r})_1<c_0.
\)
The previous applies, and contradicts $0 \in \overline{\{|U| > 0\}}.$
\end{proof}

The same collapsing idea yields the following gluing result
under the strict concavity assumption
\eqref{eq:weight-strict-concavity}, i.e.,  a block with sufficiently
small total separation in $B_1$ must coincide in $B_{1/2}$.
\begin{proposition}
\label{lem:membrane-block-collapse}
Assume \eqref{eq:weight-strict-concavity}.
There are $\ep_0,\delta_0>0$, depending only on
$d,N,\gamma,\mathbf w_\gamma$ s.t. if $U\in\A$ minimizes $\J_{N,d}^{\gamma,G}$ in $B_1$
with $G\in\mathbf G_{\delta_0}$ and
\(
  \|u_k-u_{k+m-1}\|_{L^\infty(B_1)}\leq\ep_0 \ \text{for } 1 \leq k < k+m-1 \leq N,
\)
then
\begin{equation}\label{eq:membrane-block-collapse}
  u_k=\cdots=u_{k+m-1}
  \qquad\text{in }B_{1/2}.
\end{equation}
\end{proposition}

\begin{proof}
 Denote 
\begin{equation}\label{eq:membrane-block-collapse-mu}
     \mu:=
  \min_{1\leq i<N}
  \left(w_i-\frac{w_{i-1}+w_{i+1}}2\right)>0.
\end{equation} 
All constants below depend only on $d, N, \gamma , \mathbf{w}_\gamma$ and change from line to line.

\emph{Step 1.}
Fix $1\leq i<N$ and assume
$\|u_i-u_{i+1}\|_{L^\infty(B_1)}\leq\ep_0$.
We first prove that $u_i=u_{i+1}$ in $B_{1/2}$. For $0<s<t<1$, choose $\zeta\in C_c^\infty(B_t)$ with
$0\leq\zeta\leq1$, $\zeta=1$ in $B_s$, and
$|\nabla\zeta|\leq C/(t-s)$. Define 
\[
  v_j^\pm=
  \begin{cases}
    (1-\zeta^2)u_j
    +\dfrac{\zeta^2}{2}
      \bigl(u_i+u_{i+1}\pm\theta(u_i-u_{i+1})\bigr),
      &j=i,i+1,\\
    u_j,&j\notin\{i,i+1\},
  \end{cases}
\]
where $\theta\in(0,1)$ will be chosen close to $1$. With $V^\pm = (v_1^\pm, \cdots, v_N^\pm),$ minimality gives
\begin{equation}
    \label{eq:averaged minimality inequality-nondegeneracy}
  0\leq
  \frac12\sum_\pm
  \J_{N,d}^{\gamma,G}(V^\pm,\mathbf w_\gamma,B_t)
  -\J_{N,d}^{\gamma,G}(U,\mathbf w_\gamma,B_t).
\end{equation}
\begin{figure}[htbp]
\centering
\vspace{-15pt} 
\setlength{\abovecaptionskip}{1pt}
\begin{tikzpicture}[x=1.55cm,y=1.3cm,font=\footnotesize,
  declare function={
    ramp(\x)=min(1,max(0,(abs(\x)-.55)/.50));
    collapse(\x)=1-10*ramp(\x)^3+15*ramp(\x)^4-6*ramp(\x)^5;
  }]
\foreach \panelx/\direction/\signtext/\linecolor in
  {0/1/+/blue!70!black,6/-1/-/orange!85!black}{
  \begin{scope}[xshift=\panelx cm]

    \fill[\linecolor,opacity=.06]
      (-.55,-.85) rectangle (.55,.85);
    \foreach \a in {-1.05,-.55,.55,1.05}
      \draw[black!25,densely dotted]
        (\a,-.90) -- (\a,.85);

    \draw[black!45,dashed] (-1.4,.65) -- (1.4,.65);
    \draw[black!45,dashed] (-1.4,-.65) -- (1.4,-.65);
    \node[anchor=west] at (1.43,.65) {$u_i$};
    \node[anchor=west] at (1.43,-.65) {$u_{i+1}$};

    \foreach \branch in {-1,1}
      \draw[\linecolor,very thick]
        plot[domain=-1.4:1.4,samples=181,variable=\x]
        (\x,{.65*(\branch*(1-collapse(\x))
                         +\direction*.7*collapse(\x))});

    \node[text=\linecolor] at (0,{\direction*.20})
      {$v_i^{\signtext}=v_{i+1}^{\signtext}$};
    \node at (0,1.10) {$V^{\signtext}$};

    \draw[->,>=stealth,black!65]
      (-1.5,-.95) -- (1.65,-.95) node[right] {$x$};
    \foreach \a/\ticktext in {-1.05/-t,-.55/-s,.55/s,1.05/t}{
      \draw[black!65] (\a,-.92) -- (\a,-.98);
      \node[anchor=north] at (\a,-1.00) {$\ticktext$};
    }
  \end{scope}
}
\end{tikzpicture}
\caption{Competitors $V^\pm$ with $\theta = 0.7$.}
\label{fig:pair-collapse-competitors}
\end{figure}
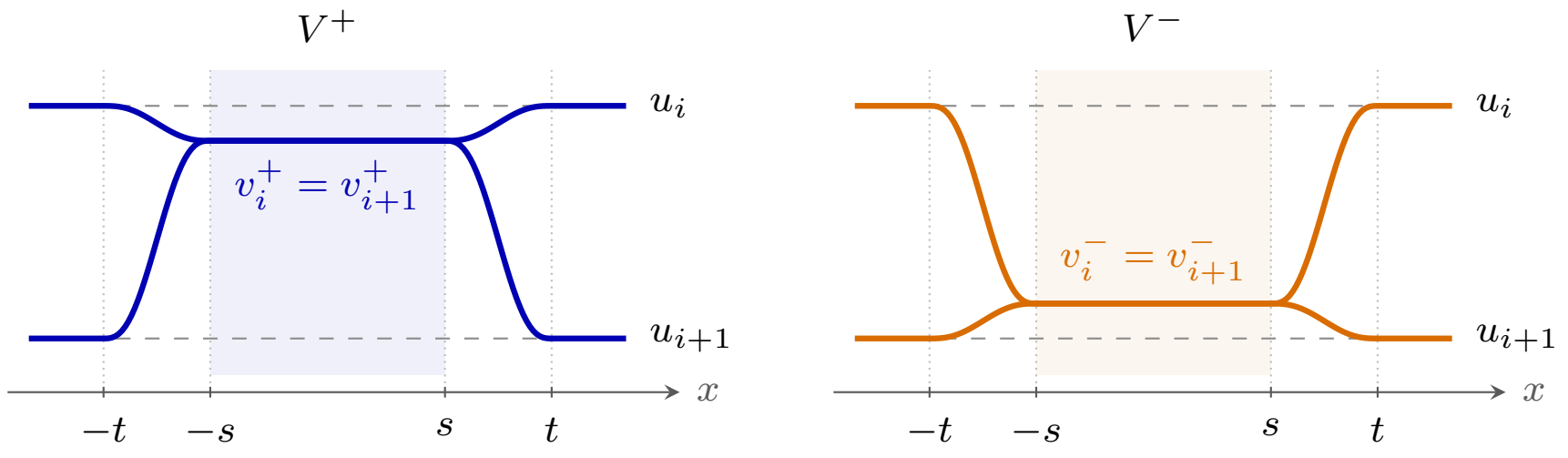

We estimate the three contributions to this difference. For the Dirichlet terms, 
\begin{align*}
  \frac12\sum_\pm|\nabla V^\pm|^2-|\nabla U|^2
  &=
  -\nabla(u_i-u_{i+1})\cdot
    \nabla\bigl(\zeta^2(u_i-u_{i+1})\bigr)
  +\frac{1+\theta^2}{2}
    \left|\nabla\bigl(\zeta^2(u_i-u_{i+1})\bigr)\right|^2\\
  &\leq
  -\frac{1-\theta^2}{4}\zeta^2
     |\nabla(u_i-u_{i+1})|^2
  +C(u_i-u_{i+1})^2|\nabla\zeta|^2.
\end{align*}
For the gap potential, the only modified gaps are
\begin{align*}
  v_i^\pm-v_{i+1}^\pm
  &=(1-\zeta^2)(u_i-u_{i+1}),\\
  v_{i-1}^\pm-v_i^\pm
  &=u_{i-1}-u_i
    +\frac{1\mp\theta}{2}\zeta^2(u_i-u_{i+1}),
    \qquad i>1,\\
  v_{i+1}^\pm-v_{i+2}^\pm
  &=u_{i+1}-u_{i+2}
    +\frac{1\pm\theta}{2}\zeta^2(u_i-u_{i+1}),
    \qquad i<N-1.
\end{align*}
Using $(a+b)^\gamma-a^\gamma\leq b^\gamma$ and averaging
over the two signs, we obtain
\begin{align*}
  \frac12\sum_\pm W_N^\gamma(V^\pm,\mathbf w_\gamma)
    -W_N^\gamma(U,\mathbf w_\gamma)\leq
  \left[
    -w_i\bigl(1-(1-\zeta^2)^\gamma\bigr)
    +A_i(\theta)\zeta^{2\gamma}
  \right](u_i-u_{i+1})^\gamma,
\end{align*}
where
\[
  A_i(\theta):=
  \frac{w_{i-1}+w_{i+1}}2
  \left[
    \left(\frac{1-\theta}{2}\right)^\gamma
    +\left(\frac{1+\theta}{2}\right)^\gamma
  \right].
\]
Using \eqref{eq:membrane-block-collapse-mu}, one may choose
$\theta=\theta(N,\gamma,\mathbf w_\gamma)\in(0,1)$
sufficiently close to $1$ that
\[
  A_i(\theta)\leq w_i-\frac{\mu}{2},
  \qquad i=1,\ldots,N-1.
\]
For the perturbation, note
\[
  |V^\pm-U|
  =\sqrt{\frac{1+\theta^2}{2}}\,
     \zeta^2(u_i-u_{i+1})
  \leq\zeta^2(u_i-u_{i+1}).
\]
Thus taking $\delta_0\leq\mu/4$ and $\ep_0\leq1$, we obtain
\[
  \frac12\sum_\pm G(x,V^\pm)-G(x,U)
  \leq\delta_0\zeta^2(u_i-u_{i+1})
  \leq\frac{\mu}{4}\zeta^2(u_i-u_{i+1})^\gamma.
\]

In $B_s$, where $\zeta=1$, the combined potential and
perturbation terms are at most
\[
  \left(-w_i+A_i(\theta)+\frac{\mu}{4}\right)
       (u_i-u_{i+1})^\gamma
  \leq-\frac{\mu}{4}(u_i-u_{i+1})^\gamma.
\]
On $B_t\setminus B_s$, they are bounded above by
$C(u_i-u_{i+1})^\gamma$. Thus, setting
\[
  E(r):=\int_{B_r}
  \left(
    |\nabla(u_i-u_{i+1})|^2
    +(u_i-u_{i+1})^\gamma
  \right)\,dx,
\]
\eqref{eq:averaged minimality inequality-nondegeneracy} gives
\[
  E(s)\leq C_0\bigl(E(t)-E(s)\bigr)
  +\frac{C_0}{(t-s)^2}
     \int_{B_t}(u_i-u_{i+1})^2\,dx.
\]
Rearranging and applying \cite[Lemma~4.3]{HL11}, we obtain
\begin{equation}
\label{eq:membrane-block-collapse-energy-comparison-bound}
  E(s)\leq\frac{C}{(t-s)^2}
     \int_{B_t}(u_i-u_{i+1})^2\,dx,
  \qquad 0<s<t<1.
\end{equation}
In particular taking $s=3/4$ and $t=7/8$ gives $E(3/4)\leq C\ep_0^2$. 

\emph{Step 2.}
We show that small $\ep_0$ forces $E(1/2)=0$. Since $0\leq u_i-u_{i+1}\leq1$,
Sobolev embedding gives
\[
  \left(\int_{B_t}(u_i-u_{i+1})^q\,dx\right)^{2/q}
  \leq C\int_{B_t}
  \left(
    |\nabla(u_i-u_{i+1})|^2+(u_i-u_{i+1})^2
  \right)\,dx
  \leq CE(t), \qquad q=2+2/d,
\]
uniformly for $1/2\leq t\leq3/4$.
H\"older's then yields
\begin{align*}
  \int_{B_t}(u_i-u_{i+1})^2\,dx
  &\leq
  \left(\int_{B_t}(u_i-u_{i+1})^\gamma\,dx\right)^{
    \frac{q-2}{q-\gamma}}
  \left(\int_{B_t}(u_i-u_{i+1})^q\,dx\right)^{
    \frac{2-\gamma}{q-\gamma}}\\
  &\leq
  CE(t)^{\frac{q-2}{q-\gamma}
          +\frac{q(2-\gamma)}{2(q-\gamma)}}
  =CE(t)^{1+\sigma}, \qquad  \sigma=\frac{(q-2)(2-\gamma)}{2(q-\gamma)}
  =\frac{2-\gamma}{d(2-\gamma)+2}>0.
\end{align*}
Combining this with
\eqref{eq:membrane-block-collapse-energy-comparison-bound},
we obtain
\[
  E(s)\leq\frac{C}{(t-s)^2}E(t)^{1+\sigma},
  \qquad \frac12\leq s<t\leq\frac34.
\]
For $R_n=\frac12+2^{-n-2}$, this gives
\[
  E(R_{n+1})\leq C4^nE(R_n)^{1+\sigma}.
\]
Since $E(3/4)\leq C\ep_0^2$, we may choose $\ep_0$
so small that $CE(3/4)^\sigma\leq4^{-1/\sigma}$. Induction then gives
\[
  E(R_n)\leq4^{-n/\sigma}E(3/4)\longrightarrow0.
\]
Therefore $E(1/2)=0$, and continuity implies
$u_i=u_{i+1}$ in $B_{1/2}$. Finally, all choices above are uniform in $i$. For $0\leq u_i-u_{i+1}\leq u_k-u_{k+m-1}\leq\ep_0$, applying the above to each internal gap gives \eqref{eq:membrane-block-collapse}.
\end{proof}

As a consequence, we obtain non-degeneracy of an individual gap.
\begin{corollary}
\label{prop:individual-gap-nondegeneracy}
Assume \eqref{eq:weight-strict-concavity}.
There exist
$c_1=c_1(d,N,\gamma,\mathbf w_\gamma)>0$
and
$r_1=r_1(d,N,\gamma,K,\mathbf w_\gamma)\in(0,1/2)$
such that if $U\in\A$ minimize
$\J_{N,d}^{\gamma,G}$ in $B_1$ with
\(
 0\in\partial\{u_i>u_{i+1}\},
\) then
\begin{equation*}
 \sup_{B_r}(u_i-u_{i+1})\geq c_1r^\beta,
 \qquad 0<r<r_1.
\end{equation*}
\end{corollary}

\begin{proof}
Apply \Cref{lem:membrane-block-collapse} to
$I=\{i,i+1\}$ and argue by contradiction as in \Cref{prop:gap-nondegeneracy}.
\end{proof}

\section{Weiss monotonicity and blowups}\label{sec:monotonicity-blowups}
In this section, $N \geq 2$, $\gamma \in (0, 1)$, $\beta = 2/(2 - \gamma) \in (1, 2)$ as in \eqref{eq:homogeneity}, $G \in \mathbf{G}_K$ for $K > 0$, $\mathbf{w}_\gamma \in (0, \infty)^{N-1},$ unless otherwise stated. We develop Weiss formula, and show that blowups at full contact free boundary points are $\beta$-homogeneous nontrivial global minimizing cones. We first prove compactness of `minimality'.
\begin{proposition}\label{prop:gap-compactness}
Let $U_j$ locally minimize
$\J_{N,d}^{\gamma,G_j}(\cdot,\mathbf w_\gamma,B_1)$ with $G_j\in\mathbf G_K$ for fixed $K\geq0$.
Suppose
\begin{equation}\label{eq:gm-to-zero}
 U_j\longrightarrow U
 \quad\text{locally uniformly in }B_1,
 \qquad
 G_j(x,V)\longrightarrow0
 \quad\text{for a.e. }x\in B_1
 \text{ and every }V\in\R^N.
\end{equation}
Then $U$ locally minimizes $\J_{N,d}^{\gamma}$, $U_j\to U$ in $C^1_{\mathrm{loc}}(B_1;\R^N)$, and
\begin{equation}\label{eq:gap-energy-convergence}
 \J_{N,d}^{\gamma,G_j}(U_j,\mathbf w_\gamma,B_r)
 \longrightarrow
 \J_{N,d}^{\gamma}(U,\mathbf w_\gamma,B_r), \qquad \forall \ B_r\Subset B_1.
\end{equation}
\end{proposition}

\begin{proof}
\Cref{thm:optimal-regularity} and Ascoli-Arzela upgrade local uniform convergence to $C_{\text{loc}}^1$. Combining with \eqref{eq:gm-to-zero} and dominated convergence, \eqref{eq:gap-energy-convergence} follows. To prove minimality, fix $B_r\Subset B_R\Subset B_1$, then take competitor
$V$ for $U$ in $B_r$ and extend by $U$ outside $B_r$. Choose
$\zeta\in C_c^\infty(B_R)$ s.t. $0\leq\zeta\leq1$ and $\zeta=1$ near $\overline{B_r}$.
Set
\(
 V_j=\zeta V+(1-\zeta)U_j,
\)
which is an admissible competitor for $U_j$ in $B_R$.
Since $V=U$ wherever $1-\zeta\ne0$,
\[
 V_j-V=(1-\zeta)(U_j-U)
 \longrightarrow0
 \qquad\text{strongly in }H^1(B_R;\R^N).
\]
Thus the Dirichlet energies converge. Also \(
 \int_{B_R}
 \left|W_N^\gamma(V_j,\mathbf w_\gamma)
       -W_N^\gamma(V,\mathbf w_\gamma)\right|dx
 \leq C\|V_j-V\|_{L^2(B_R)}^\gamma
 \to 0,
\)
and from
\(
 |G_j(x,V_j)|
 \leq K|V_j-V|+|G_j(x,V)|,\ |G_j(x,V)|\leq K|V|\in L^1(B_R),
\)
dominated convergence gives
$G_j(\cdot,V_j)\to0$ in $L^1(B_R)$. Passing to the limit using \eqref{eq:gap-energy-convergence} and minimality then yields
\[
\begin{aligned}
 \J_{N,d}^{\gamma}(U,\mathbf w_\gamma,B_R)
 &=\lim_{j\to\infty}
   \J_{N,d}^{\gamma,G_j}(U_j,\mathbf w_\gamma,B_R)\leq\lim_{j\to\infty}
   \J_{N,d}^{\gamma,G_j}(V_j,\mathbf w_\gamma,B_R)
 =\J_{N,d}^{\gamma}(V,\mathbf w_\gamma,B_R).
\end{aligned}
\]
Since $V=U$ on $B_R\setminus B_r$, the terms outside
$B_r$ cancel. This proves local minimality of $U$.
\end{proof}

\subsection{Weiss formula}
\label{subsec:Weiss-monotonicity}
Define the Weiss
energy associated with $G =0$ by
\begin{equation}\label{eq:membrane-weiss-g0}
 \W_U(r)
 =r^{-d-2\beta+2}
   \int_{B_r}\left(|\nabla U|^2+W_N^\gamma(U)\right)dx
 -\beta r^{-d-2\beta+1}
   \int_{\partial B_r}|U|^2\,dS,
 \qquad 0<r<1.
\end{equation}

\begin{proposition}
\label{prop:gap-weiss-g0}
Let $U\in\A_{\mathrm{loc}}$ be vertically critical
and domain stationary for $\J_{N,d}^{\gamma}$ as in \Cref{def:critical-points,def:gap-domain-criticality}.
Then $\W_U$ is locally absolutely continuous on $(0,1)$, and
\begin{equation}\label{eq:membrane-weiss-g0-derivative}
 \W_U'(r)
 =\frac{2}{r^{d+2\beta}}
   \int_{\partial B_r}
   |r\partial_\nu U-\beta U|^2\,dS, \qquad \text{for a.e. $r \in (0, 1)$.}
\end{equation}
In particular, $\W_U$ is nondecreasing.
It is constant if and only if
$U$ is $\beta$-homogeneous.
\end{proposition}
\begin{proof}
Set
\(
 a=d+2\beta-2.\) We write $\W_U(r)=r^{-a}E(r)-\beta r^{-a-1}H(r)$ for  \(
 E(r)=\int_{B_r}\left(|\nabla U|^2+W_N^\gamma(U)\right)dx\) and
 \(
 H(r)=\int_{\partial B_r}|U|^2\,dS.
\)
For a.e. $r$, using \(H'(r)
=\frac{d-1}{r}H(r)
   +2\int_{\partial B_r}U\cdot\partial_\nu U\,dS,\)
\begin{equation}\label{eq:gap-weiss-differentiation}
\begin{aligned}
 \W_U'(r)
 &=r^{-a}\left(E'(r)-\frac arE(r)\right)
   -2\beta r^{-a-1}
     \int_{\partial B_r}U\cdot\partial_\nu U\,dS
   +2\beta^2r^{-a-2}H(r).
\end{aligned}
\end{equation}
 With $\beta\gamma=2\beta-2$ and 
\(  E'(r)
 =\int_{\partial B_r}
   \left(|\nabla U|^2+W_N^\gamma(U)\right)dS,
\) one obtain
\[
\begin{aligned}
 rE'(r)-aE(r)
 &=2r\int_{\partial B_r}|\partial_\nu U|^2\,dS
   -2\beta\int_{B_r}|\nabla U|^2\,dx
   -(2\beta-2)\int_{B_r}W_N^\gamma(U)\,dx\\
 &=2r\int_{\partial B_r}|\partial_\nu U|^2\,dS
   -\beta\int_{B_r}
       \left(2|\nabla U|^2+\gamma W_N^\gamma(U)\right)dx=2r\int_{\partial B_r}|\partial_\nu U|^2\,dS
   -2\beta\int_{\partial B_r}U\cdot\partial_\nu U\,dS,
\end{aligned}
\]
 using \eqref{eq:gap-rellich-g0} for the first equality, and \eqref{eq:gap-vertical-g0} for the last.
Substituting into \eqref{eq:gap-weiss-differentiation} gives
\[
\begin{aligned}
 \W_U'(r)
 &=2r^{-a}\int_{\partial B_r}|\partial_\nu U|^2\,dS
   -4\beta r^{-a-1}
     \int_{\partial B_r}U\cdot\partial_\nu U\,dS
   +2\beta^2r^{-a-2}H(r)=2r^{-a-2}\int_{\partial B_r}
       |r\partial_\nu U-\beta U|^2\,dS.
\end{aligned}
\]
This proves
\eqref{eq:membrane-weiss-g0-derivative}. Finally, $\W_U'(r)= 0$ a.e. on
$(r_1,r_2)$ iff
\(
 \partial_r\left(r^{-\beta}U(r\theta)\right)=0
\) for a.e. \((r,\theta)
 \in(r_1,r_2)\times\S^{d-1}.\)
\end{proof}

For $G\in\mathbf G_K$, we avoid differentiating the energy and use minimality to compare $U$ with the $\beta$-homogeneous extension of its boundary values, which gives an almost monotonicity formula for \eqref{eq:membrane-weiss-g0}.

\begin{proposition}\label{prop:gap-weiss}
Let $U$ locally minimize
$\J_{N,d}^{\gamma,G}(\cdot,\mathbf w_\gamma,B_1)$. Suppose
\begin{equation}\label{eq:beta-growth-gk-weiss}
 |U(x)|\leq M|x|^\beta
 \qquad\text{in }B_1.
\end{equation}
Then $\W_U$ is locally absolutely continuous on $(0,1)$,
and, for a.e. $r\in(0,1)$,
\begin{equation}\label{eq:membrane-weiss-gk-derivative}
 \W_U'(r)
 \geq\frac{1}{r^{d+2\beta}}
       \int_{\partial B_r}|r\partial_\nu U-\beta U|^2\,dS
       -CKMr^{1-\beta},
 \qquad C=C(d,\gamma)>0.
\end{equation}
Consequently, after increasing $C$, the following is nondecreasing on $(0,1)$,
\begin{equation}\label{eq:gap-weiss-gk-adjusted}
 r\longmapsto\W_U(r)+CKMr^{2-\beta}.
\end{equation}
\end{proposition}

\begin{proof}
    Throughout this proof, $C=C(d,\gamma)>0$ denotes constant changing from line to line.
Use $a,E,H$ as in proof of \Cref{prop:gap-weiss-g0}.
The absolute continuity and \eqref{eq:gap-weiss-differentiation} remain valid. For a.e. $r\in(0,1)$, define $\beta$-homogeneous extension of boundary data on $\p B_r$ as
\[
 U^{(r)}(\rho\theta)
 =\left(\frac\rho r\right)^\beta U(r\theta),
 \qquad 0<\rho<r,\quad\theta\in\S^{d-1},
 \qquad U^{(r)}(0)=0.
\] 
Writing $\nabla_\tau$ for the tangential gradient, since $\beta\gamma=2\beta-2$ and $a=d+2\beta-2>0$, we compute
\begin{equation}\label{eq:hom-replacement-energy1}
\begin{aligned}
 \J_{N,d}^{\gamma}(U^{(r)},\mathbf w_\gamma,B_r)
 &=\int_0^r\int_{\S^{d-1}}
   \left(\frac\rho r\right)^{2\beta-2}
   \left(\frac{\beta^2}{r^2}|U(r\theta)|^2
         +|\nabla_\tau U(r\theta)|^2
         +W_N^\gamma(U(r\theta))\right)
   \rho^{d-1}\,dS(\theta)\,d\rho\\
 &=\frac r a\int_{\partial B_r}
   \left(|\nabla_\tau U|^2
         +\frac{\beta^2}{r^2}|U|^2
         +W_N^\gamma(U)\right)dS.
\end{aligned}
\end{equation}
Using \eqref{eq:beta-growth-gk-weiss},
$|U^{(r)}(x)|\leq M|x|^\beta$ in $B_r$,
thus
\[
\begin{aligned}
 \int_{B_r}\left(G(x,U^{(r)})-G(x,U)\right)dx
 &\leq K\int_{B_r}|U^{(r)}-U|\,dx\leq2KM\int_{B_r}|x|^\beta\,dx
 \leq CKMr^{d+\beta}.
\end{aligned}
\]
Minimality and \eqref{eq:hom-replacement-energy1}
therefore imply
\[
\begin{aligned}
 E(r)
 &\leq\J_{N,d}^{\gamma}(U^{(r)},\mathbf w_\gamma,B_r)
       +\int_{B_r}\left(G(x,U^{(r)})-G(x,U)\right)dx\\
 &\leq\frac r a\int_{\partial B_r}
   \left(|\nabla_\tau U|^2
         +\frac{\beta^2}{r^2}|U|^2
         +W_N^\gamma(U)\right)dS
       +CKMr^{d+\beta}\\
 &=\frac r a\left(E'(r)
       -\int_{\partial B_r}|\partial_\nu U|^2\,dS
       +\frac{\beta^2}{r^2}H(r)\right)
       +CKMr^{d+\beta}.
\end{aligned}
\]
Rearranging gives
\[
 E'(r)-\frac arE(r)
 \geq\int_{\partial B_r}|\partial_\nu U|^2\,dS
       -\frac{\beta^2}{r^2}H(r)
       -CKMr^{d+\beta-1}.
\]
Substituting into \eqref{eq:gap-weiss-differentiation},
we obtain
\[
\begin{aligned}
 \W_U'(r)
 &\geq r^{-a}\int_{\partial B_r}|\partial_\nu U|^2\,dS
       -2\beta r^{-a-1}
        \int_{\partial B_r}U\cdot\partial_\nu U\,dS
       +\beta^2r^{-a-2}H(r)
       -CKMr^{1-\beta}\\
 &=r^{-a-2}\int_{\partial B_r}
        |r\partial_\nu U-\beta U|^2\,dS
       -CKMr^{1-\beta}.
\end{aligned}
\]
This proves \eqref{eq:membrane-weiss-gk-derivative}.
Since $\beta<2$, the error $r^{1-\beta}$ is integrable
at $0$, so increasing $C$ proves \eqref{eq:gap-weiss-gk-adjusted}.
\end{proof}

\subsection{Homogeneous blowups}

To study free boundaries for any minimizer, one zoom in at $x_0 \in \p\{u_i > u_{i+1}\}$ and would like to understand its blow-up profile. Using \Cref{rmk:reduction-from-maximal-contact-block}, one may take a maximal contact block around $x_0$, rescale to $B_1$ and subtract its average, so it suffices to work with $U \in \A_0$ with $0 \in \p\{|U| > 0\}.$ We show that every such blow-up limit is a nontrivial $\beta$-homogeneous global minimizing cone.

\begin{theorem}
\label{thm:gap-blowup-cone}
Let $U\in\A_{0}$ locally minimize
$\J_{N,d}^{\gamma,G}(\cdot,\mathbf w_\gamma,B_1)$, and $0\in\partial\{|U|>0\}$.
Then there is $r_j\downarrow 0$ s.t.
\begin{equation}\label{eq:blowup-convergence}
 U_{0,r_j}(x)=r_j^{-\beta}U(r_jx)
 \longrightarrow U_0(x)
 \qquad\text{in }C^1_{\mathrm{loc}}(\R^d;\R^N),
\end{equation}
where $U_0 \in C_{\text{loc}}^{1, \beta - 1}(\R^d; \R^N)$ is a zero-mean local minimizer to $\J_{N,d}^{\gamma}$ in $\R^d$ with $0 \in \p\{|U_0| > 0\}$. Moreover, 
\begin{equation}\label{eq:blowup-homogeneity}
 U_0(tx)=t^\beta U_0(x),
 \qquad x\in\R^d,\quad t>0.
\end{equation}
\end{theorem}

\begin{proof}
By \Cref{thm:optimal-regularity},
$U\in C^{1,\beta-1}_{\mathrm{loc}}$.
Then using $U(0)=0$ and $U \in \A_0$, $\nabla U(0)=0$, so for some $M > 0$,
\[
 |U(x)|+|x||\nabla U(x)|
 \leq M|x|^\beta
 \qquad\text{in }B_{1/2}.
\]
Arzelà--Ascoli then yields \eqref{eq:blowup-convergence},
with $U_0\in C^{1, \beta - 1}$, preserving the order and zero mean. The rescaled perturbations vanish in the limit, so \Cref{prop:gap-compactness} shows that
$U_0$ minimizes $\J_{N,d}^{\gamma}$ among all ordered maps in every ball of $\R^d$.
Moreover, \Cref{prop:gap-nondegeneracy} gives
\(\lim_{j\to\infty}
   r_j^{-\beta}\sup_{B_{r_j}}|U|
 \geq c_0,
\)
so $U_0\not\equiv0$ and $0 \in \partial \{|U_0| > 0\}.$

Since $\W_U(r)\geq-CM^2$ for
$0<r<1/2$, \Cref{prop:gap-weiss} implies
$\W_U(r)+CKMr^{2-\beta}$ is nondecreasing, and there exists a finite limit
$\W_U(0+)=\lim_{r\downarrow0}\W_U(r)$.
For every $\rho>0$, scaling and
\eqref{eq:blowup-convergence} therefore give
\[
 \W_{U_0}(\rho)
 =\lim_{j\to\infty}\W_{U_{0,r_j}}(\rho)
 =\lim_{j\to\infty}\W_U(r_j\rho)
 =\W_U(0+).
\]
Hence $\W_{U_0}$ is constant. Since $U_0$ minimizes
the unperturbed energy, \Cref{prop:gap-weiss-g0}
implies \eqref{eq:blowup-homogeneity}.
\end{proof}

\section{One-dimensional critical cones and flat single-gap minimizing cones}
\label{sec:single-gap-regularity}

In this section, $N\geq2$, $\gamma\in(0,1)$,
$\beta=2/(2-\gamma)\in(1,2)$ as in \eqref{eq:homogeneity},
$G\equiv0$, and $\mathbf w_\gamma\in(0,\infty)^{N-1}$,
unless otherwise stated. We write
\[
 \Gamma_i(U)=\partial\{u_i>u_{i+1}\}\cap\O,
 \qquad i=1,\ldots,N-1.
\]
At $x_0\in\Gamma_i(U)$,
\Cref{rmk:reduction-from-maximal-contact-block,thm:gap-blowup-cone}
give nontrivial $\beta$-homogeneous minimizing blowup cones for the maximal contact block containing
$u_i,u_{i+1}$.

We first describe one-dimensional critical cones using
\Cref{def:critical-points,lem:interior-equations}.
Then we show minimality of the flat single-gap cones under
concave weight assumptions, and use the gluing phenomenon
to reduce minimizers near such cones to scalar minimizers,
whose regularity is developed in
\Cref{sec:scalar-perturbed-flatness}.

\subsection{One-dimensional critical cones}
\label{sec:one-dimensional-critical-beta-cones}
Throughout this subsection, $d=1$. Let $U$ be ordered $\beta$-homogeneous map that is vertically critical for $\J_{N,1}^\gamma$. On $(0, \infty)$, write its general form as
\[
 U(x)=\mathbf C x^\beta,\qquad
 \mathbf C=
 \left(
 \underbrace{c_1,\ldots,c_1}_{n_1},
 \ldots,
 \underbrace{c_m,\ldots,c_m}_{n_m}
 \right),
 \qquad c_1>\cdots>c_m,
\]
where $n_1,\ldots,n_m\geq1$ and
$n_1+\cdots+n_m=N$.
Set the indices separating consecutive blocks as
\[
 S=\{k_1,\ldots,k_{m-1}\} \subseteq \{1, \cdots, N-1\}, \qquad
 k_\ell=\sum_{j=1}^{\ell}n_j
 \qquad \text{for }\ell=1,\ldots,m.
\] Also set \(k_0=0.\)
Thus the $\ell$th contact block
$I_\ell = \{k_{\ell-1}+1,\ldots,k_\ell\}$ consists of
$n_\ell$ membranes of size $c_\ell x^\beta$. Each block is separated
from its neighboring blocks throughout $(0,\infty)$, thus \eqref{eq:vertical-first-variation-equation}
with $G=0$ gives
\[
\begin{aligned}
 2\frac{d^2}{dx^2}
 \left(\sum_{j=k_{\ell-1}+1}^{k_\ell}u_j\right)
 &=
 \gamma w_{k_\ell}
 (u_{k_\ell}-u_{k_\ell+1})^{\gamma-1}-
 \gamma w_{k_{\ell-1}}
 (u_{k_{\ell-1}}-u_{k_{\ell-1}+1})^{\gamma-1}
 \qquad\text{in }(0,\infty).
\end{aligned}
\]
Substituting $u_j(x)=c_\ell x^\beta$ for $j\in I_\ell$
and using $\beta(\gamma-1)=\beta-2$, we obtain \begin{equation*}
\begin{aligned}
 2\beta(\beta-1)n_\ell c_\ell x^{\beta-2}
 =\gamma\left[
 w_{k_\ell}(c_\ell-c_{\ell+1})^{\gamma-1}
 -w_{k_{\ell-1}}(c_{\ell-1}-c_\ell)^{\gamma-1}
 \right]x^{\beta-2}.
\end{aligned}
\end{equation*}
We omit the first term on right-hand side if $\ell = m$ and omit the second if $\ell = 1.$ In view of \(\frac{\gamma}{2\beta(\beta-1)}=\beta^{-2},\) the block heights $(c_1, \ldots, c_m)$ therefore solves the constrained system
\begin{equation}\label{eq:weighted-critical-coefficient-system}
 \left\{
 \begin{aligned}
 n_\ell c_\ell
 &=\beta^{-2}\left[
 w_{k_\ell}(c_\ell-c_{\ell+1})^{\gamma-1}
 -w_{k_{\ell-1}}(c_{\ell-1}-c_\ell)^{\gamma-1}
 \right],
 \qquad \ell=1,\ldots,m,\\
 c_1&>\cdots>c_m.
 \end{aligned}
 \right.
\end{equation}
In particular, summing these equations gives
$\sum_{\ell=1}^m n_\ell c_\ell=0$.

\begin{lemma}
\label{lem:eighted-critical-coefficient-system}
For every nonempty $S\subseteq\{1,\ldots,N-1\}$,
\eqref{eq:weighted-critical-coefficient-system}
has a unique solution $c^S$.
\end{lemma}
\begin{proof}
Consider
\[
 F_S(c)
 =\frac12\sum_{\ell=1}^m n_\ell c_\ell^2
 -\frac{\beta^{-2}}{\gamma}
  \sum_{\ell=1}^{m-1}
  w_{k_\ell}(c_\ell-c_{\ell+1})^\gamma
\]
on the closed convex set
$\{c\in\R^m:c_1\geq\cdots\geq c_m\}$.
Since $s\mapsto-s^\gamma$ is convex on $[0,\infty)$,
$F_S$ is strictly convex. It is also coercive,
so it has a unique minimizer $c$. We claim this minimizer has $c_\ell>c_{\ell+1}$ for every $\ell<m$.
Indeed, if $c_j=c_{j+1}$ for some $j<m$, define
\(
 c(t)=(c_1+t,\ldots,c_j+t,c_{j+1},\ldots,c_m)
\) for $t > 0.$ We obtain
\[
\begin{aligned}
 F_S(c(t))-F_S(c)
 &=\frac12\sum_{\ell=1}^j
   n_\ell\bigl((c_\ell+t)^2-c_\ell^2\bigr)
   -\frac{\beta^{-2}}{\gamma}w_{k_j}t^\gamma=t\sum_{\ell=1}^j n_\ell c_\ell
   +\frac{k_j}{2}t^2
   -\frac{\beta^{-2}}{\gamma}w_{k_j}t^\gamma<0
\end{aligned}
\]
for small $t>0$, which contradicts minimality. Thus it satisfies the Euler--Lagrange equations 
\eqref{eq:weighted-critical-coefficient-system}. Conversely, every solution to \eqref{eq:weighted-critical-coefficient-system} is a global minimizer of $F_S$ by convexity. Strict convexity then gives uniqueness.
\end{proof}
Now for any set of indices $S \subseteq \{1, \cdots, N-1\}$,  define the membrane coefficient vector
\[
 \mathbf{C}^S=
 \left(
 \underbrace{c_1^S,\ldots,c_1^S}_{n_1},
 \ldots,
 \underbrace{c_m^S,\ldots,c_m^S}_{n_m}
 \right)\in\R^N,
 \qquad \mathbf{C}^\varnothing=0.
\]

\begin{theorem}
\label{thm:one-dimensional-critical-cones}
A nonzero $\beta$-homogeneous map $U$ over $\R$ is vertically critical
for $\J_{N,1}^\gamma$ if and only if
\begin{equation}\label{eq:weighted-critical-coefficients-active}
 U(x)=\mathbf{C}^{S^+}x_+^\beta+\mathbf{C}^{S^-}x_-^\beta,
 \qquad x_-=(-x)_+,\qquad S^+\cup S^-\ne\varnothing,
\end{equation}
for some $S^+,S^-\subseteq\{1,\ldots,N-1\}$.
\end{theorem}
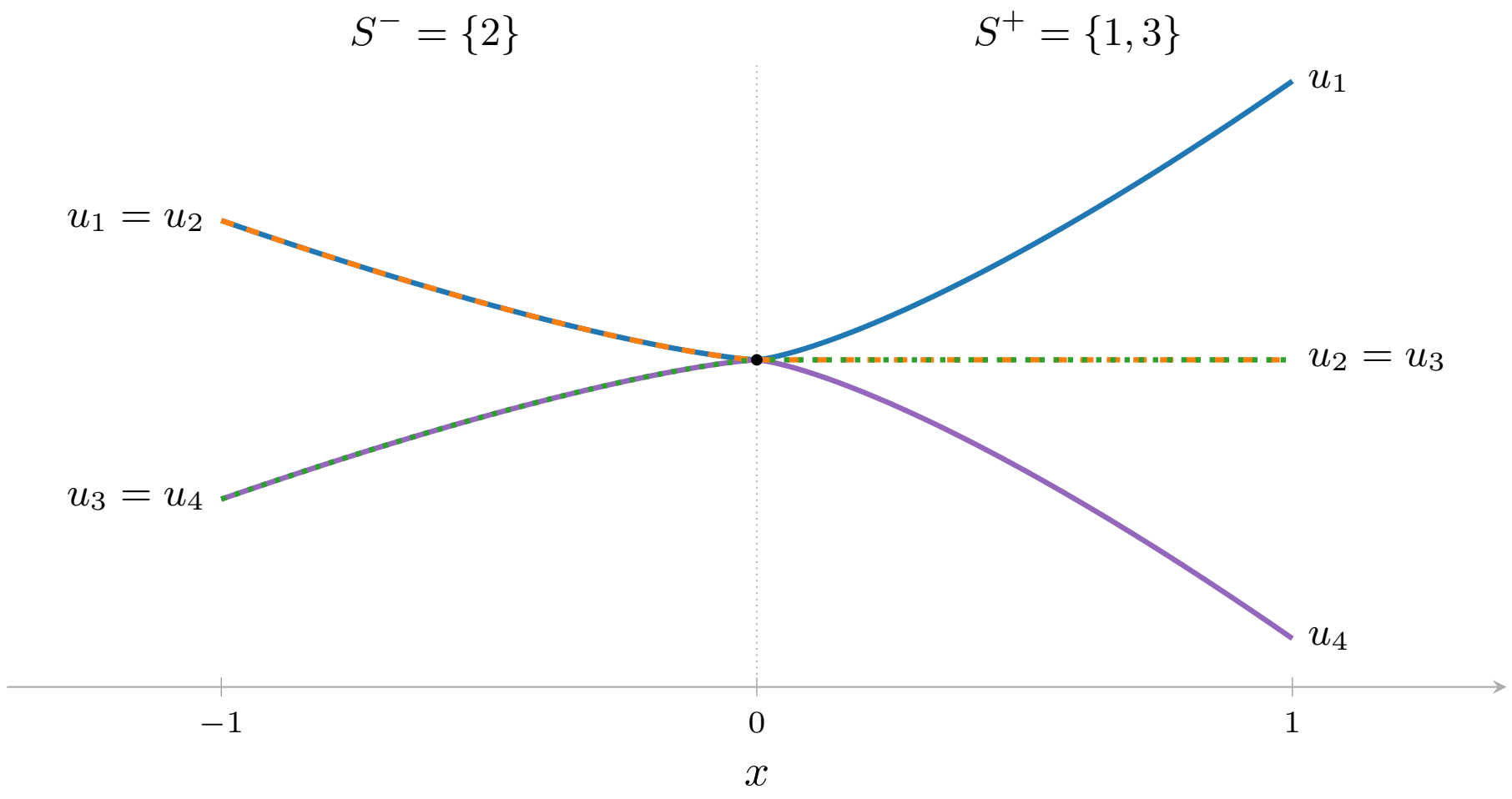
\begin{figure}[htbp]
\centering
\vspace{-6pt} 
\setlength{\abovecaptionskip}{1pt}
\definecolor{membraneone}{RGB}{31,119,180}
\definecolor{membranetwo}{RGB}{255,127,14}
\definecolor{membranethree}{RGB}{44,160,44}
\definecolor{membranefour}{RGB}{148,103,189}
\begin{tikzpicture}
\pgfmathsetmacro{\betacone}{4/3}
\pgfmathsetmacro{\bcone}{pow(3/4,4/3)}
\def\leftprofile{(\bcone/2)*pow(max(-x,0),\betacone)}
\def\rightprofile{\bcone*pow(max(x,0),\betacone)}

\begin{axis}[
    width=0.72\linewidth, height=0.33\linewidth,
    scale only axis,
    domain=-1:1, samples=201,
    xmin=-1.40, xmax=1.40, ymin=-0.80, ymax=0.88,
    axis x line=bottom, axis y line=none,
    axis line style={gray!65},
    xlabel={$x$}, xtick={-1,0,1}, ytick=\empty,
    tick label style={font=\scriptsize}, clip=false
]
\draw[gray!55,densely dotted]
    (axis cs:0,-0.77) -- (axis cs:0,0.72);

\addplot[membraneone,very thick]
    {(\leftprofile)+(\rightprofile)};
\addplot[membranefour,very thick]
    {-(\leftprofile)-(\rightprofile)};
\addplot[membranetwo,very thick,dashed]
    {\leftprofile};
\addplot[membranethree,very thick,dotted]
    {-(\leftprofile)};

\fill (axis cs:0,0) circle (1.3pt);

\node[font=\small] at (axis cs:-0.60,0.80) {$S^-=\{2\}$};
\node[font=\small] at (axis cs: 0.60,0.80) {$S^+=\{1,3\}$};

\node[left,font=\small]  at (axis cs:-1,{\bcone/2})  {$u_1=u_2$};
\node[left,font=\small]  at (axis cs:-1,{-\bcone/2}) {$u_3=u_4$};
\node[right,font=\small] at (axis cs:1,\bcone)       {$u_1$};
\node[right,font=\small] at (axis cs:1,0)            {$u_2=u_3$};
\node[right,font=\small] at (axis cs:1,{-\bcone})    {$u_4$};
\end{axis}
\end{tikzpicture}
\caption{A four-membrane critical cone from
\Cref{thm:one-dimensional-critical-cones},
with $\gamma=1/2$ and $\mathbf w_\gamma=\mathbf 1$.}
\label{fig:one-dimensional-critical-cone}
\end{figure}
\begin{proof}
Homogeneity gives
$U(x)=\mathbf{C}^+x_+^\beta+\mathbf{C}^-x_-^\beta$.
On each half-line, they solve
\eqref{eq:weighted-critical-coefficient-system}, and uniqueness gives $\mathbf{C}^\pm=\mathbf{C}^{S^\pm}$. Conversely, for $U$ as in \eqref{eq:weighted-critical-coefficients-active}, take any admissible two-sided variation $\Phi$
as in \Cref{def:critical-points}. The components
$\Phi_j$ coincide on $\{u_j = u_{j+1}\}$, then
group \eqref{eq:vertical-first-variation}
by these blocks and integrate by parts on each
half-line gives $0$ by
\eqref{eq:weighted-critical-coefficient-system}.  Since $\beta > 1$ gives
\(
 U'(0+)=U'(0-)=0,\ U''\in L^1_{\mathrm{loc}}(\R;\R^N),
\) integration by parts is valid and there is no boundary contribution at the origin.
\end{proof}

We may also compute the Weiss energy for $U$ as in \eqref{eq:weighted-critical-coefficients-active}.
Multiplying \eqref{eq:weighted-critical-coefficient-system}
by $c_\ell^S$ and summing gives
\[
 \beta^2|\mathbf{C}^S|^2
 =\beta^2\sum_{\ell=1}^m n_\ell(c_\ell^S)^2
 =\sum_{\ell=1}^{m-1}
   w_{k_\ell}(c_\ell^S-c_{\ell+1}^S)^\gamma
 =W_N^\gamma(\mathbf{C}^S,\mathbf w_\gamma).
\]
The same identity holds for $S=\varnothing$.
Thus the Dirichlet and potential densities agree on
each half-line. Using $\beta\gamma=2\beta-2$ in
\eqref{eq:membrane-weiss-g0}, we obtain the Weiss energy independent of $r > 0$
\begin{equation}\label{eq:weighted-critical-cone-weiss-energy}
\begin{aligned}
 \W_U
 &=\left(\frac2{2\beta-1}-\frac1\beta\right)
 \left[
 W_N^\gamma(\mathbf{C}^{S^+},\mathbf w_\gamma)
 +W_N^\gamma(\mathbf{C}^{S^-},\mathbf w_\gamma)
 \right]=\frac{
 W_N^\gamma(\mathbf{C}^{S^+},\mathbf w_\gamma)
 +W_N^\gamma(\mathbf{C}^{S^-},\mathbf w_\gamma)}
 {\beta(2\beta-1)}.
\end{aligned}
\end{equation}

\subsection{Flat single-gap minimizing cones}
\label{sec:one-sided-single-gap-flat-type}

\begin{figure}[htbp]
\centering
\vspace{-6pt} 
\setlength{\abovecaptionskip}{1pt}
\definecolor{flatupper}{RGB}{31,119,180}
\definecolor{flatlower}{RGB}{44,160,44}
\begin{tikzpicture}[
    x={(2.05cm,0cm)}, y={(0.62cm,0.90cm)}, z={(0cm,2.65cm)},
    line join=round, line cap=round, font=\small
]
\pgfmathsetmacro{\flatgamma}{1/2}
\pgfmathsetmacro{\flatbeta}{2/(2-\flatgamma)}
\pgfmathsetmacro{\flata}{pow(5/(6*\flatbeta^2),1/(2-\flatgamma))}
\pgfmathsetmacro{\flattop}{3*\flata/5}
\pgfmathsetmacro{\flatbottom}{-2*\flata/5}

\filldraw[fill=gray!12,draw=gray!60]
    (-1.45,0,0) -- (0,0,0) -- (0,1,0) -- (-1.45,1,0) -- cycle;
\draw[gray!75,thick] (-1.45,0,0) -- (0,0,0);
\node at (-0.80,0.52,0) {$u_1=\cdots=u_5=0$};

\foreach \coef/\col in {\flatbottom/flatlower,\flattop/flatupper}{
    \filldraw[fill=\col!18,draw=\col!80!black]
        plot[domain=0:1.65,samples=51,variable=\t]
            (\t,0,{\coef*pow(\t,\flatbeta)})
        -- plot[domain=1.65:0,samples=51,variable=\t]
            (\t,1,{\coef*pow(\t,\flatbeta)})
        -- cycle;

    \foreach \s in {0.45,0.85,1.25}{
        \draw[\col!55,thin]
            (\s,0,{\coef*pow(\s,\flatbeta)})
            -- (\s,1,{\coef*pow(\s,\flatbeta)});
    }
    \draw[\col!85!black,very thick]
        plot[domain=0:1.65,samples=51,variable=\t]
            (\t,0,{\coef*pow(\t,\flatbeta)});
}

\draw[black,dashed,thick] (0,0,0) -- (0,1,0)
    node[above left] {$x_1=0$};

\node[right,text=flatupper!80!black]
    at (1.65,1,{\flattop*pow(1.65,\flatbeta)}) {$u_1=u_2$};
\node[right,text=flatlower!75!black]
    at (1.65,1,{\flatbottom*pow(1.65,\flatbeta)}) {$u_3=u_4=u_5$};
\end{tikzpicture}
\caption{The flat single-gap minimizing cone $U_2^{e_1}$
for $N=5$, $\gamma=1/2$, and $\mathbf w_\gamma=\mathbf 1$.}
\label{fig:flat-single-gap-cone}
\end{figure}
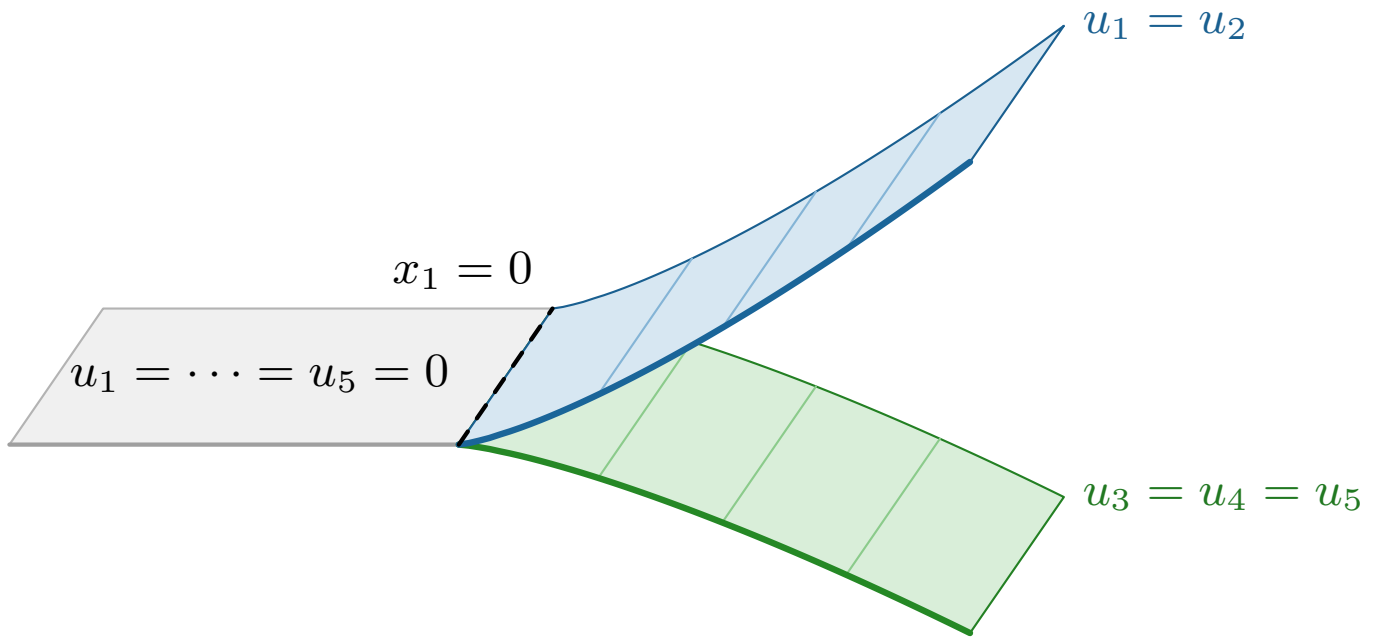

We next consider cones with only two blocks of membranes. For the single gap indexed at $S=\{k\}$, the block sizes are
$n_1=k$ and $n_2=N-k$. Set \(
 a_k=c_1-c_2>0.
\)
The zero-mean condition $kc_1+(N-k)c_2=0$ then gives
\(
 c_1=\frac{N-k}{N}a_k,\
 c_2=-\frac{k}{N}a_k.
\)
Thus
\[\mathbf{C}^{\{k\}} = a_k  \mathbf{e}_{N, k}, \qquad  \mathbf{e}_{N, k}=\frac1N
 \left(
 \underbrace{N-k,\ldots,N-k}_{k},
 \underbrace{-k,\ldots,-k}_{N-k}
 \right).\]
Solving for $a_k$ using \(|\mathbf{e}_{N, k}|^2 = \frac{k(N-k)}{N}\) and the first equation in
\eqref{eq:weighted-critical-coefficient-system} gives
\begin{equation}
\label{eq:single-gap-ak}
 a_k=
 \left(\frac{w_{k}}{\beta^2|\mathbf{e}_{N, k}|^2}\right)^{1/(2-\gamma)}, \qquad  \text{satisfying} \qquad |\mathbf{e}_{N, k}|^2 a_k=\beta^{-2} w_{k}a_k^{\gamma-1}.
\end{equation}
One therefore define the flat single-gap critical cone in all dimensions $d \geq 1$ as
\begin{equation}\label{eq:flat-gap-cone}
 U_k^\nu(x)=a_k(x\cdot\nu)_+^\beta\mathbf{e}_{N, k}, \qquad \forall \ \nu \in \S^{d-1}.
\end{equation}
These are trivial extensions of $1$-d cones whose first $k$ and last $N-k$ membranes coincide; see \Cref{fig:flat-single-gap-cone}.

\subsubsection{Minimality of $U_k^\nu$}

For this result, we assume only non-strict  concavity from \eqref{eq:weight-strict-concavity}. The weights therefore lie above the chords joining
$(0,0)$ to $(k,w_k)$ and $(k,w_k)$ to $(N,0)$ 
\begin{equation}
    \label{eq:flat-weight-compatibility-chord}
  w_i\geq\frac{i}{k}w_k
  \qquad 1\leq i<k,
  \qquad
  w_i\geq\frac{N-i}{N-k}w_k
  \qquad k<i<N.
\end{equation}

\begin{proposition}\label{prop:single-gap-minimizers}
Under \eqref{eq:flat-weight-compatibility-chord}, the cone $U_k^\nu$ from \eqref{eq:flat-gap-cone} is a global minimizing cone for  $\J_{N,d}^\gamma$.
\end{proposition}

We first record Federer's dimension reduction.

\begin{lemma}
\label{lem:translation-invariant-minimizers}
Let $d\geq2$ and let
$U\in H^1_{\mathrm{loc}}(\R^{d-1};\R^N)$ be ordered.
Define
\(
  \widetilde U(x',x_d)=U(x'),\ (x',x_d)\in\R^{d-1}\times\R.
\)
Then $U$ globally minimizes $\J_{N,d-1}^\gamma$
if and only if $\widetilde U$ globally minimizes
$\J_{N,d}^\gamma$.
\end{lemma}

\begin{proof}
The argument is the same as in \cite[Lemma~10.10]{Vel23}.
\end{proof}

Let's explain our strategy. 
After reducing to $U=a_kx_+^\beta\mathbf e_{N,k}$, scalar minimality of $a_kx_+^\beta$ gives minimality of $U$
among non-negative competitors on the ray \(\mathbf e_{N,k}\). We need to project any ordered competitor
$V=(v_1,\ldots,v_N)$ onto $\mathbf e_{N,k}$ without increasing the energy. Consider
$(v_j-v_\ell)\mathbf e_{N,k}$, $j\leq k<\ell$.
Using \eqref{eq:flat-weight-compatibility-chord}, these are admissible competitors whose averaged energy
does not exceed that of $V$, thus at least one gives the desired projection.

\begin{proof}[Proof of \Cref{prop:single-gap-minimizers}]
By
\Cref{lem:translation-invariant-minimizers}, it suffices
to prove minimality of
\[
  U(x)=v_0(x)\mathbf e_{N,k},
  \qquad v_0(x)=a_kx_+^\beta,
  \qquad x\in\R.
\]
Fix a bounded interval $I = (r, s)$ and an admissible
competitor $V=(v_1,\ldots,v_N)$ for $U$ in $I$, one has
\[
  v_j-v_\ell\geq0,
  \qquad (v_j-v_\ell)-v_0\in H_0^1(I),
  \qquad
  v_j-v_\ell=\sum_{i=j}^{\ell-1}(v_i-v_{i+1}), \qquad \forall \ 1 \leq j \leq k < \ell \leq N.
\]

Counting the appearances of $(v_i-v_{i+1})^\gamma$
according to whether $i<k$, $i=k$, or $i>k$, gives
\begin{equation}\label{eq:single-gap-minimizers-identity}
\begin{aligned}
 \sum_{j=1}^k\sum_{\ell=k+1}^N
 \sum_{i=j}^{\ell-1}(v_i-v_{i+1})^\gamma
 &=(N-k)\sum_{i=1}^{k-1}i(v_i-v_{i+1})^\gamma+k(N-k)(v_k-v_{k+1})^\gamma+k\sum_{i=k+1}^{N-1}
 (N-i)(v_i-v_{i+1})^\gamma.
\end{aligned}
\end{equation}
Combining \eqref{eq:single-gap-minimizers-identity}
with subadditivity and
\eqref{eq:flat-weight-compatibility-chord} , we obtain
\[
\begin{aligned}
 \frac{w_{k}}{k(N-k)}
 \sum_{j=1}^k\sum_{\ell=k+1}^N(v_j-v_\ell)^\gamma&\leq
 w_{k}\left[
 \sum_{i=1}^{k-1}\frac{i}{k}(v_i-v_{i+1})^\gamma
 +(v_k-v_{k+1})^\gamma
 +\sum_{i=k+1}^{N-1}
 \frac{N-i}{N-k}(v_i-v_{i+1})^\gamma
 \right]\leq W_N^\gamma(V,\mathbf w_\gamma).
\end{aligned}
\]
For the Dirichlet energy, using
$|\mathbf e_{N,k}|^2=k(N-k)/N$, we obtain
\begin{align*}
  \frac{|\mathbf e_{N,k}|^2}{k(N-k)}
  \sum_{j=1}^k\sum_{\ell=k+1}^N|(v_j-v_\ell)'|^2
  &\leq\frac1N
    \sum_{1\leq j<\ell\leq N}|v_j'-v_\ell'|^2=|V'|^2-\frac1N\left|\sum_{j=1}^Nv_j'\right|^2\leq|V'|^2.
\end{align*}
Consequently,
\[
  \J_{N,1}^\gamma(V,\mathbf w_\gamma,I)
  \geq
  \frac1{k(N-k)}
  \sum_{j=1}^k\sum_{\ell=k+1}^N
  \int_I\left[
    |\mathbf e_{N,k}|^2|(v_j-v_\ell)'|^2
    +w_k(v_j-v_\ell)^\gamma
  \right]dx.
\]

It remains to use the scalar minimality of $v_0$.
By \eqref{eq:single-gap-ak},
\[
  |\mathbf e_{N,k}|v_0'
  =\sqrt{w_k}\,v_0^{\gamma/2}
  \qquad\text{a.e. in }I.
\]
Thus, for any nonnegative $v\in H^1(I)$ with
$v-v_0\in H_0^1(I)$, completing the square gives
\begin{equation}
\label{eq:single-gap-minimizers-comparison-energy-0}
\begin{aligned}
  \int_I\left(
    |\mathbf e_{N,k}|^2|v'|^2+w_kv^\gamma
  \right)dx
  &\geq
  2|\mathbf e_{N,k}|\sqrt{w_k}
  \int_r^s v^{\gamma/2}v'\,dx=
  \frac{2|\mathbf e_{N,k}|\sqrt{w_k}}{1+\gamma/2}
  \left[v_0^{1+\gamma/2}\right]_r^s=
  \int_I\left(
    |\mathbf e_{N,k}|^2|v_0'|^2+w_kv_0^\gamma
  \right)dx.
\end{aligned}
\end{equation}
Applying \eqref{eq:single-gap-minimizers-comparison-energy-0} to every $v_j-v_\ell$ and averaging,
we conclude
\[
  \J_{N,1}^\gamma(V,\mathbf w_\gamma,I)
  \geq
  \int_I\left(
    |\mathbf e_{N,k}|^2|v_0'|^2+w_kv_0^\gamma
  \right)dx
  =\J_{N,1}^\gamma(U,\mathbf w_\gamma,I).
\]
\end{proof}

\begin{remark}\label{rmk:single-gap-weight-assumption}
Positive weights alone do not guarantee the result of
\Cref{prop:single-gap-minimizers}. For example 
take $N=3$, $d=1$, and $\mathbf w_\gamma=(1,\varepsilon)$
with $\varepsilon>0$. Consider the cone
\[
 U(x)=a_1x_+^\beta\mathbf e_{3,1}
 =\left(\frac23v_0,-\frac13v_0,-\frac13v_0\right),
 \qquad v_0=a_1x_+^\beta.
\]
By \eqref{eq:single-gap-ak}, $a_1$ is independent of
$\varepsilon$. Now fix $D=(1,2)$ and a nonzero
$\varphi\in C_c^\infty(D)$ with $\varphi\geq0$. Set
\[
 V_t=(\frac23v_0,-\frac13v_0+t\varphi,-\frac13v_0-t\varphi).
\]
For sufficiently small $t>0$, \( V_t\) is admissible competitor for $U$. By concavity of $s\mapsto s^\gamma$, we obtain
\[
\begin{aligned}
 \J_{3,1}^\gamma(V_t,\mathbf w_\gamma,D)
 -\J_{3,1}^\gamma(U,\mathbf w_\gamma,D)
 &=2t^2\int_D|\varphi'|^2\,dx
   +\int_D\bigl[(v_0-t\varphi)^\gamma-v_0^\gamma\bigr]\,dx
+\varepsilon(2t)^\gamma\int_D\varphi^\gamma\,dx\\ &\leq
 2t^2\int_D|\varphi'|^2\,dx
 -\gamma t\int_Dv_0^{\gamma-1}\varphi\,dx
 +\varepsilon(2t)^\gamma\int_D\varphi^\gamma\,dx.
\end{aligned}
\]
Choose $t>0$ small, then $\varepsilon>0$ small so the whole expression is negative. Therefore opening the gap with weight $\ep$ while reducing
the first gap lowers the energy, yet $U$ is not minimizing.
\end{remark}

\subsubsection{Regularity of $\Gamma_i$ near $U_{k}^\nu$}
\label{sec:flat-type-partial-regularity}

We now study regularity of the free boundary $\Gamma_i$ for minimizers near $U_k^\nu$. Under concave weights, gluing reduces the single gap around $\Gamma_i$ to a scalar Alt-Phillips problem with Lipschitz perturbation. Improvement of flatness developed from \Cref{sec:scalar-perturbed-flatness} then gives $C^{1, \alpha}$ of $\Gamma_i$.

We record the direct consequence of \Cref{lem:membrane-block-collapse} applied to $\{1, \ldots, k\}$ and $\{k+1, \ldots N\}.$
\begin{corollary}
\label{lem:flat-gluing}
Assume \eqref{eq:weight-strict-concavity}.
There are $\ep_0,\delta_0>0$ depending only on
$d,N,\gamma,\mathbf w_\gamma$ s.t. if $U$
locally minimizes
$\J_{N,d}^{\gamma,G}$ in $B_1$ with 
\(
 G\in\mathbf G_{\delta_0},\
 \|U-U_k^\nu\|_{L^\infty(B_1)}\leq\ep_0
 \ \text{for $1\leq k<N$ and }\nu\in\S^{d-1},
\)
then
\[
 u_1=\cdots=u_k,\qquad
 u_{k+1}=\cdots=u_N \qquad \text{in }B_{1/2}.
\]
\end{corollary}

Let $U$ be a minimizer and $x_0\in\Gamma_i(U)$.
Restricting to the maximal contact block, rescaling,
and subtracting its mean as in
\Cref{rmk:reduction-from-maximal-contact-block},
we may work with $U\in\A_0$ in $B_1$, with
$U(0)=0$ and $0\in\Gamma_k(U)$, after relabeling
the membranes and weights. In this setting we obtain the following.
\begin{corollary}\label{cor:single-gap-regularity}
Assume \eqref{eq:weight-strict-concavity}. If $U\in\A_0$ locally minimize
$\J_{N,d}^{\gamma,G}$ in $B_1$ with $0\in\Gamma_k(U)$, and for $r_j\downarrow0$,
\[
 r_j^{-\beta}U(r_jx)\longrightarrow U_k^\nu(x)
 \qquad\text{locally uniformly in }\R^d,
\]
then $\Gamma_k(U)$ is locally $C^{1,\alpha}$ for $\alpha=\alpha(d,\gamma)\in(0,2-\beta)$, with unit normal $\nu$ at $0$ pointing into $\{u_k > u_{k+1}\}$.
\end{corollary}

\begin{proof}
By \eqref{eq:rescale-G}, the rescaled perturbation has
Lipschitz bound $Kr_j^{2-\beta}\to0$. The convergence
and \Cref{lem:flat-gluing} therefore give
$U=(u_k-u_{k+1})\mathbf e_{N,k}$ near $0$.
Comparing with non-negative scalar competitors shows that
\[
 g=\left(
 \frac{|\mathbf e_{N,k}|^2}{w_{k}}
 \right)^{1/(2-\gamma)}(u_k-u_{k+1})
\]
locally minimizes \eqref{eq:scalar-perturbed-energy}
with a Lipschitz perturbation. By
\eqref{eq:single-gap-ak}, its blowup along $r_j$ is
$\beta^{-\beta}(x\cdot\nu)_+^\beta$.
The conclusion follows from
\Cref{thm:scalar-flat-blowup-regularity}, since
$\Gamma_k(U)=\partial\{g>0\}$ near $0$.
\end{proof}

Undoing the reduction gives $C^{1,\alpha}$ regularity of $\Gamma_i(U)$ near $x_0$ provided assumptions for \Cref{cor:single-gap-regularity} holds.

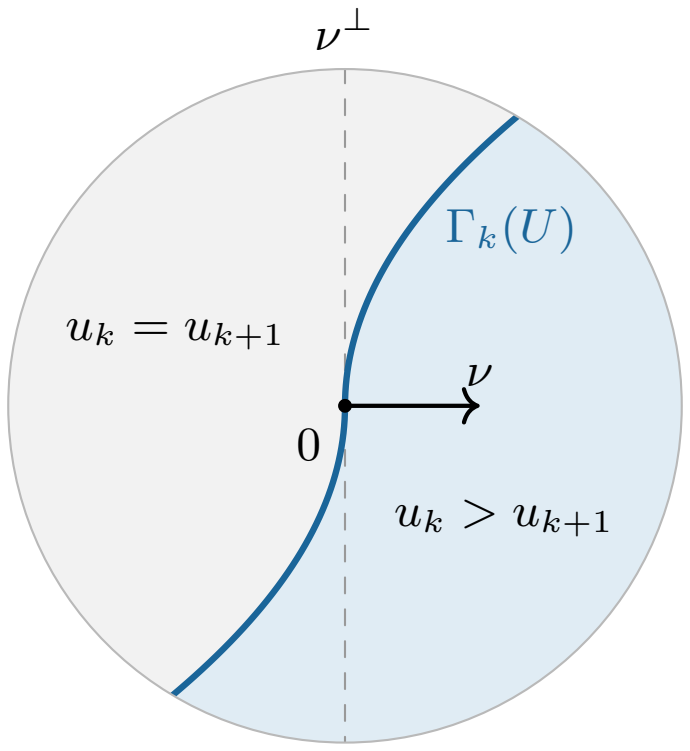
\begin{figure}[htbp]
\centering
\vspace{-6pt} 
\setlength{\abovecaptionskip}{1pt}
\definecolor{gapblue}{RGB}{31,119,180}
\begin{tikzpicture}[
    x=1.35cm, y=1.35cm, font=\small,
    line cap=round, line join=round
]
\begin{scope}
    \clip (0,0) circle (1.65);
    \fill[gray!10] (-1.7,-1.7) rectangle (1.7,1.7);

    % The gap is open to the right of the free boundary.
    \fill[gapblue!14]
        plot[domain=-1.7:1.7,samples=101,variable=\t]
            ({0.42*\t*abs(\t)},\t)
        -- (1.7,1.7) -- (1.7,-1.7) -- cycle;

    % Tangent line and free boundary.
    \draw[gray!80,dashed] (0,-1.7) -- (0,1.7);
    \draw[gapblue!85!black,very thick]
        plot[domain=-1.7:1.7,samples=101,variable=\t]
            ({0.42*\t*abs(\t)},\t);
\end{scope}
\draw[gray!55] (0,0) circle (1.65);

\fill (0,0) circle (1.3pt) node[below left=1pt] {$0$};
\draw[->,thick] (0,0) -- (0.66,0) node[above] {$\nu$};

\node at (-0.83,0.35) {$u_k=u_{k+1}$};
\node at (0.78,-0.55) {$u_k>u_{k+1}$};
\node[above] at (0,1.65) {$\nu^\perp$};
\node[right,text=gapblue!85!black] at (0.40,0.86) {$\Gamma_k(U)$};
\end{tikzpicture}
\caption{The free boundary in \Cref{cor:single-gap-regularity}
for $d=2$. The dashed line is its tangent at $0$.}
\label{fig:single-gap-regularity}
\end{figure}

\section{One-dimensional minimizing cones for three membranes}
\label{sec:N3-d1-minimizing-cones}

Throughout this section, $N=3$, $d=1$, $\gamma \in (0, 1)$, $\beta = 2/(2 - \gamma) \in (1, 2)$ as in \eqref{eq:homogeneity}, and $G=0$ unless otherwise stated. For weights, we assume either $\mathbf{w}_\gamma = (w_{1}, w_{2}) \in (0, \infty)^2$ or further impose
\begin{equation}
     \label{eq:weights_notcanonical}
 \qquad \frac12<\rho=\frac{ w_{ 2}}{ w_{1}}<2, \qquad \text{and} \qquad \rho(\gamma) \to \rho_* \in (\frac{1}{2}, 2), \qquad \gamma \uparrow 1.
\end{equation}

We first list the six types of critical cones in
\Cref{tab:n3-gap-cones}. Type~I is minimizing for every
$0<\gamma<1$ under \eqref{eq:weights_notcanonical}
by \Cref{prop:single-gap-minimizers}, while
\Cref{lem:common-active-gap-instability} rules out
Types~III, V, and VI for arbitrary positive weights.
For the remaining Types~II and IV, we prove
non-minimality for small $\gamma$ with arbitrary positive
weights, and minimality for $\gamma$ sufficiently near $1$ under
\eqref{eq:weights_notcanonical}; see
\Cref{thm:type-ii-minimality,thm:type-iv-minimality}. In particular, for $\rho\equiv1$, Type~II is not minimizing for
$0<\gamma<(14+4\sqrt6)/25\approx 0.952$ and is minimizing for
$\gamma_{\mathrm{II}}\approx 0.968\leq\gamma<1$, while Type~IV is not minimizing
for $0<\gamma\leq2/3$ and is minimizing for
$\gamma_{\mathrm{IV}}\approx 0.911\leq\gamma<1$; see \Cref{rmk:typeIIrho=1,rmk:typeIVrho=1}.

\subsection{One-dimensional critical cones for \texorpdfstring{$N=3$}{J}}
\label{sec:n3-critical-cones}

By \Cref{thm:one-dimensional-critical-cones}, nonzero critical cones take the form
\begin{equation*}
 U(x)=\mathbf C^+x_+^\beta+\mathbf C^-x_-^\beta,
 \qquad
 \mathbf C^\pm\in
 \{0,\mathbf C^{\{1\}}, \mathbf C^{\{2\}},\mathbf C^{\{1,2\}}\},
 \qquad (\mathbf C^+,\mathbf C^-)\ne(0,0),
\end{equation*}
where $\mathbf{C}^{\{k\}} = a_k \mathbf{e}_{3, k}$ are given by \eqref{eq:single-gap-ak}, and $\mathbf C^{\{1,2\}}=(c_1,c_2,c_3)$ is the unique
solution from \Cref{lem:eighted-critical-coefficient-system}
\begin{equation}\label{eq:n3-two-gap-coefficients}
 \left\{
 \begin{aligned}
 c_1&=\beta^{-2} w_{1}(c_1-c_2)^{\gamma-1},\\
 c_2&=\beta^{-2}\left[
 w_{2}(c_2-c_3)^{\gamma-1}
 -w_{1}(c_1-c_2)^{\gamma-1}\right],\\
 c_3&=-\beta^{-2} w_{2}(c_2-c_3)^{\gamma-1},\\
 c_1&>c_2>c_3.
 \end{aligned}
 \right.
\end{equation}
In particular, using $c_1+c_2+c_3=0$ and \(\mathbf e_{3,1}=\frac13(2,-1,-1),  \mathbf e_{3,2}=\frac13(1,1,-2),\) one may rewrite 
\begin{equation}
    \label{eq:n3-two-gap-coefficients-rewrite}\mathbf{C}^{\{1, 2\}} = (c_1 -c_2)\mathbf{e}_{3,1} + (c_2 - c_3)\mathbf{e}_{3, 2}.
\end{equation}

Up to reflection in $x$, the six types of critical cones are listed in \Cref{tab:n3-gap-cones} below.
\begin{table}[htbp]
\centering
\renewcommand{\arraystretch}{1.5}
\setlength{\tabcolsep}{4pt}
\begin{tabular}{@{}c l c@{}}
\toprule
\multicolumn{2}{c}{Type} & $U(x)$\\
\midrule
$\mathrm I$
 & one-sided single gap
 & $\mathbf{C}^{\{k\}}x_+^\beta$\\
$\mathrm{II}$
 & one-sided two gaps
 & $\mathbf{C}^{\{1,2\}}x_+^\beta$\\
$\mathrm{III}$
 & two-sided, same single gap
 & $\mathbf{C}^{\{k\}}|x|^\beta$\\
$\mathrm{IV}$
 & two-sided, different single gaps
 & $\mathbf{C}^{\{1\}}x_+^\beta
    +\mathbf{C}^{\{2\}}x_-^\beta$\\
$\mathrm V$
 & two-sided, two gaps / single gap
 & $\mathbf{C}^{\{1,2\}}x_+^\beta
    +\mathbf{C}^{\{k\}}x_-^\beta$\\
$\mathrm{VI}$
 & two-sided two gaps
 & $\mathbf{C}^{\{1,2\}}|x|^\beta$\\
\bottomrule
\end{tabular}
\setlength{\abovecaptionskip}{2pt}
\setlength{\belowcaptionskip}{1pt}
\caption{One-dimensional critical cones for three membranes, $k=1,2$.}
\label{tab:n3-gap-cones}
\end{table}

\begin{figure}[htbp]
\centering
\definecolor{membraneone}{RGB}{31,119,180}
\definecolor{membranetwo}{RGB}{255,127,14}
\definecolor{membranethree}{RGB}{44,160,44}
\resizebox{\textwidth}{!}{%
\begin{tikzpicture}
\pgfmathsetmacro{\acon}{pow(27/32,2/3)}
\pgfmathsetmacro{\bcon}{pow(9/16,2/3)}
\pgfmathsetmacro{\betacon}{4/3}

\begin{groupplot}[
    group style={
        group size=3 by 2,
        horizontal sep=1.6cm,
        vertical sep=1.5cm
    },
    width=0.29\textwidth,
    height=0.22\textwidth,
    domain=-1:1,
    samples=201,
    axis lines=middle,
    xmin=-1.05, xmax=1.05,
    ymin=-0.95, ymax=0.95,
    xtick={-1,1},
    ytick=\empty,
    xlabel={$x$},
    ylabel={$u$},
    every axis x label/.style={
        at={(ticklabel* cs:1.03)}, anchor=west
    },
    every axis y label/.style={
        at={(ticklabel* cs:1.03)}, anchor=south
    },
    title style={font=\small},
    tick label style={font=\scriptsize},
    clip=false,
]

\nextgroupplot[title={Type $\mathrm I$}]
\addplot[membraneone, very thick]
    {(2/3)*\acon*pow(max(x,0),\betacon)};
\addplot[membranetwo, very thick, dashed]
    {-(1/3)*\acon*pow(max(x,0),\betacon)};
\addplot[membranethree, very thick, dotted]
    {-(1/3)*\acon*pow(max(x,0),\betacon)};

\nextgroupplot[title={Type $\mathrm{II}$}]
\addplot[membraneone, very thick]
    {\bcon*pow(max(x,0),\betacon)};
\addplot[membranetwo, very thick, dashed] {0};
\addplot[membranethree, very thick, dotted]
    {-\bcon*pow(max(x,0),\betacon)};

\nextgroupplot[title={Type $\mathrm{III}$}]
\addplot[membraneone, very thick]
    {(2/3)*\acon*pow(abs(x),\betacon)};
\addplot[membranetwo, very thick, dashed]
    {-(1/3)*\acon*pow(abs(x),\betacon)};
\addplot[membranethree, very thick, dotted]
    {-(1/3)*\acon*pow(abs(x),\betacon)};

\nextgroupplot[title={Type $\mathrm{IV}$}]
\addplot[membraneone, very thick]
    {\acon*((2/3)*pow(max(x,0),\betacon)
           +(1/3)*pow(max(-x,0),\betacon))};
\addplot[membranetwo, very thick, dashed]
    {\acon*(-(1/3)*pow(max(x,0),\betacon)
           +(1/3)*pow(max(-x,0),\betacon))};
\addplot[membranethree, very thick, dotted]
    {\acon*(-(1/3)*pow(max(x,0),\betacon)
           -(2/3)*pow(max(-x,0),\betacon))};

\nextgroupplot[title={Type $\mathrm V$}]
\addplot[membraneone, very thick]
    {\bcon*pow(max(x,0),\betacon)
     +(2/3)*\acon*pow(max(-x,0),\betacon)};
\addplot[membranetwo, very thick, dashed]
    {-(1/3)*\acon*pow(max(-x,0),\betacon)};
\addplot[membranethree, very thick, dotted]
    {-\bcon*pow(max(x,0),\betacon)
     -(1/3)*\acon*pow(max(-x,0),\betacon)};

\nextgroupplot[title={Type $\mathrm{VI}$}]
\addplot[membraneone, very thick]
    {\bcon*pow(abs(x),\betacon)};
\addplot[membranetwo, very thick, dashed] {0};
\addplot[membranethree, very thick, dotted]
    {-\bcon*pow(abs(x),\betacon)};

\end{groupplot}

\begin{scope}[shift={($(group c3r2.east)+(1.2cm,1.5cm)$)}]
    \draw[membraneone, very thick] (0,0) -- (0.9,0);
    \node[anchor=west] at (1.05,0) {$u_1$};

    \draw[membranetwo, very thick, dashed]
        (0,-0.45) -- (0.9,-0.45);
    \node[anchor=west] at (1.05,-0.45) {$u_2$};

    \draw[membranethree, very thick, dotted]
        (0,-0.9) -- (0.9,-0.9);
    \node[anchor=west] at (1.05,-0.9) {$u_3$};
\end{scope}
\end{tikzpicture}%
}
\caption{The six critical cones for $w_ 1=w_2=1$ and $\gamma=1/2$. The single active gap is $k=1$.}
\label{fig:N3-six-cones-gamma-half}
\end{figure}
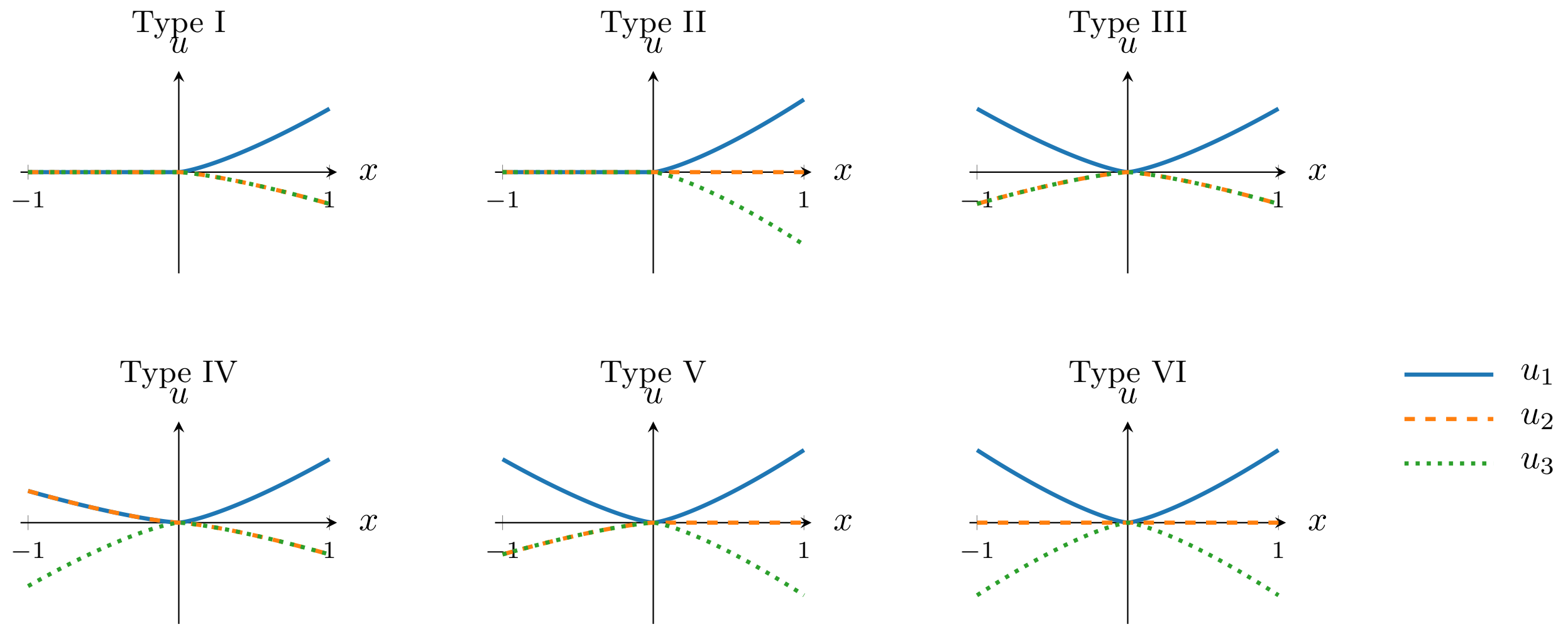

Note the Weiss energy for $1$-d critical cones follows from
\eqref{eq:weighted-critical-cone-weiss-energy}
\(
 \W_U=
 \frac{W_3^\gamma(\mathbf C^+,\mathbf w_\gamma)
       +W_3^\gamma(\mathbf C^-,\mathbf w_\gamma)}
 {\beta(2\beta-1)},
\)
with
\[  W_3^\gamma(\mathbf{C}^{\{k\}},\mathbf w_\gamma)
 =w_{k}a_k^\gamma
 =\left(\frac32\beta^{-2}\right)^{\gamma/(2-\gamma)}
   w_{k}^{2/(2-\gamma)}, \qquad  W_3^\gamma(\mathbf C^{\{1,2\}},\mathbf w_\gamma)
 =\beta^2|\mathbf C^{\{1,2\}} |^2.
\]
Under \eqref{eq:weights_notcanonical} with $\rho \equiv 1$,
symmetry gives
\(
 c_2=0,\
 c_1=-c_3=(\beta^{-2}w_{1})^{1/(2-\gamma)}.
\)
Consequently,
\[
 1<
 \frac{W_3^\gamma(\mathbf C^{\{1,2\}},\mathbf w_\gamma)}
      {W_3^\gamma(\mathbf C^{\{k\}},\mathbf w_\gamma)}
 =2\left(\frac23\right)^{\gamma/(2-\gamma)}
 <2,
 \qquad k=1,2.
\]
In this case, the six types in \Cref{tab:n3-gap-cones} are ordered
by their Weiss energies as
\[
 \W_{\mathrm I}<\W_{\mathrm{II}}<\W_{\mathrm{III}}
 =\W_{\mathrm{IV}}<\W_{\mathrm V}<\W_{\mathrm{VI}}.
\]
In particular, for the canonical choice $w_{1}=w_{2}=1$, coefficients write
\begin{equation*}
 a_1=a_2
 =\left(\frac32\beta^{-2}\right)^{1/(2-\gamma)},
 \qquad
 \mathbf C^{\{1, 2\}}=\beta^{-\beta}(1,0,-1).
\end{equation*}
We illustrate these profiles for
$\gamma=1/2$ in \Cref{fig:N3-six-cones-gamma-half}.

Now we discuss minimality for \Cref{tab:n3-gap-cones}. Under \eqref{eq:weights_notcanonical}, Type~I is minimizing for
every $0<\gamma<1$ by \Cref{prop:single-gap-minimizers}.
At $\gamma=0$, the weights are $(1,1)$, and the only nonzero
minimizing cones are of Type~I, with $a_1=a_2=\sqrt{3/2}$;
see \cite[Proposition~1.3]{DSS24}.
At $\gamma=1$, the coefficients become
\[
 a_k=\frac38w_{k},\qquad
 \mathbf C^{\{1, 2\}}=\frac14(w_{ 1},w_{2}-w_{ 1},-w_{2}).
\]
All six types satisfy the variational inequality for the
forcing $(f_1,f_2,f_3)=\frac12(w_{1},w_{ 2}-w_{1},-w_{2})$,
and hence minimize by convexity;
see \cite[Remarks~2.2 and~2.4]{SY19}. To understand how `minimality' transits as $\gamma$ ranges from $0$ to $1$, it thus remains to study minimality for Type~II to Type~VI cones under \eqref{eq:weights_notcanonical}.

\subsection{Cones with a common active gap}
\label{sec:type-iii-v-vi-nonminimality}

We first rule out Types~III, V, and VI for $0<\gamma<1$.
The following applies to cones of any number of membranes that contains a common active gap on both half lines.

\begin{lemma}
\label{lem:common-active-gap-instability}
For $N\geq2$, $\mathbf w_\gamma\in(0,\infty)^{N-1}$, a critical cone
\(
 U(x)=\mathbf C^{S^+}x_+^\beta +\mathbf C^{S^-}x_-^\beta
\) from \Cref{thm:one-dimensional-critical-cones}
is not minimizing if $S^+\cap S^-\ne\varnothing$.
\end{lemma}

\begin{figure}[htbp]
\centering
\vspace{-6pt} 
\setlength{\abovecaptionskip}{1pt}
\definecolor{membraneone}{RGB}{31,119,180}
\definecolor{membranetwo}{RGB}{255,127,14}
\definecolor{membranethree}{RGB}{44,160,44}
\resizebox{0.75\textwidth}{!}{%
\begin{tikzpicture}
\pgfmathsetmacro{\betacon}{4/3}
\pgfmathsetmacro{\acon}{pow(27/32,2/3)}
\pgfmathsetmacro{\bcon}{pow(9/16,2/3)}
\pgfmathsetmacro{\hcon}{\bcon/3}

\begin{groupplot}[
    group style={
        group size=2 by 1,
        horizontal sep=1.5cm
    },
    width=0.45\textwidth,
    height=0.32\textwidth,
    domain=-1:1,
    samples=401,
    axis lines=middle,
    xmin=-1.05, xmax=1.05,
    ymin=-0.8, ymax=0.8,
    xtick={-1,1},
    ytick=\empty,
    xlabel={$x$},
    ylabel={$u$},
    every axis x label/.style={
        at={(ticklabel* cs:1.03)}, anchor=west
    },
    every axis y label/.style={
        at={(ticklabel* cs:1.03)}, anchor=south
    },
    title style={font=\small,yshift=10pt},
    tick label style={font=\scriptsize},
    clip=false,
]

\nextgroupplot[title={Original cone $U_{\mathrm V}$}]
\addplot[membraneone, very thick]
    {\bcon*pow(max(x,0),\betacon)
     +(2/3)*\acon*pow(max(-x,0),\betacon)};
\addplot[membranetwo, very thick, dashed]
    {-(1/3)*\acon*pow(max(-x,0),\betacon)};
\addplot[membranethree, very thick, dotted]
    {-\bcon*pow(max(x,0),\betacon)
     -(1/3)*\acon*pow(max(-x,0),\betacon)};

\nextgroupplot[title={Competitor $U_h$}]
\addplot[membraneone, very thick]
    {\bcon*pow(max(x,0),\betacon)
     +(2/3)*\acon*pow(max(-x,0),\betacon)
     +(2/3)*\hcon*(1-pow(abs(x),\betacon))};
\addplot[membranetwo, very thick, dashed]
    {-(1/3)*\acon*pow(max(-x,0),\betacon)
     -(1/3)*\hcon*(1-pow(abs(x),\betacon))};
\addplot[membranethree, very thick, dotted]
    {-\bcon*pow(max(x,0),\betacon)
     -(1/3)*\acon*pow(max(-x,0),\betacon)
     -(1/3)*\hcon*(1-pow(abs(x),\betacon))};

% The first separation at x=0 is h.
\draw[gray!65, densely dotted]
    (axis cs:0,{2*\hcon/3})
    -- (axis cs:-0.18,{2*\hcon/3});
\draw[gray!65, densely dotted]
    (axis cs:0,{-\hcon/3})
    -- (axis cs:-0.18,{-\hcon/3});
\draw[<->, thick]
    (axis cs:-0.18,{-\hcon/3})
    -- (axis cs:-0.18,{2*\hcon/3})
    node[pos=1,left=3pt,font=\small] {$h$};

\end{groupplot}

\begin{scope}[
    shift={($(group c1r1.south)!0.5!(group c2r1.south)
             +(0,-0.7cm)$)}
]
    \draw[membraneone, very thick]
        (-2.4,0) -- (-1.7,0);
    \node[anchor=west] at (-1.55,0) {$u_1$};

    \draw[membranetwo, very thick, dashed]
        (-0.7,0) -- (0,0);
    \node[anchor=west] at (0.15,0) {$u_2$};

    \draw[membranethree, very thick, dotted]
        (1,0) -- (1.7,0);
    \node[anchor=west] at (1.85,0) {$u_3$};
\end{scope}
\end{tikzpicture}%
}
\caption{Type~V competitor in \Cref{lem:common-active-gap-instability} with $\gamma=1/2$, $w_{1}=w_{2}=1$,
$j=1$, and $h=\beta^{-\beta}/3$.}
\label{fig:common-active-gap-competitor}
\end{figure}
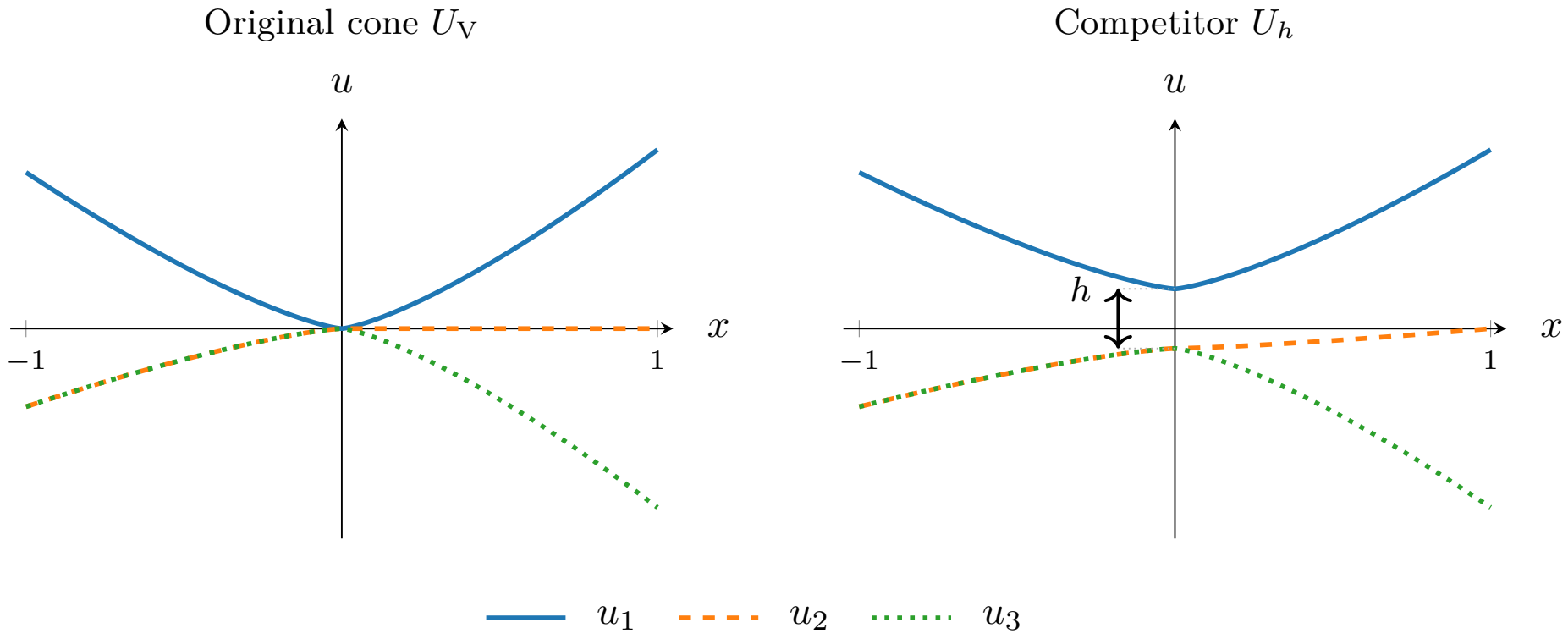
\begin{proof}
    Fix $j \in S^+ \cap S^-,$ so the first $j$ membranes form a union of contact blocks on each half-line. Fix either half-line and write
\(
 0=k_0<k_1<\cdots<k_r=j
\)
for the endpoints of these blocks. Summing
\eqref{eq:vertical-first-variation-equation} over them,
\[
\begin{aligned}
 2\left(\sum_{i=1}^j u_i\right)''
 &=\gamma\sum_{\ell=1}^r
 w_{k_\ell}
 (u_{k_\ell}-u_{k_\ell+1})^{\gamma-1}-\gamma\sum_{\ell=1}^{r-1}
 w_{k_\ell}
 (u_{k_\ell}-u_{k_\ell+1})^{\gamma-1}=\gamma w_{j}(u_j-u_{j+1})^{\gamma-1}.
\end{aligned}
\]
On the other hand, the definition of $\mathbf e_{N,j}$
and $\sum_{i=1}^N u_i=0$ give
\[
\begin{aligned}
 U''\cdot\mathbf e_{N,j}
 &=\frac{N-j}{N}\sum_{i=1}^j u_i''
   -\frac{j}{N}\sum_{i=j+1}^N u_i''=\sum_{i=1}^j u_i''
   -\frac{j}{N}\sum_{i=1}^N u_i''
 =\left(\sum_{i=1}^j u_i\right)''.
\end{aligned}
\]
Combining these identities on both half-lines yields, for any interval $I \subseteq \R$ containing $0,$
\begin{equation}\label{eq:common-active-gap-instability-equation}
     2U''\cdot\mathbf e_{N,j}
 =\gamma w_{j}(u_j-u_{j+1})^{\gamma-1}
 \qquad\text{in }I\setminus\{0\}.
\end{equation}

Take $I = (-1, 1).$ For $\phi(x) = 1 - |x|^\beta,$ consider $U_h = U + h \phi \mathbf{e}_{N, j}$ for $h > 0$; see \Cref{fig:common-active-gap-competitor}. Compute
\[
\begin{aligned}
 \J_{N,1}^\gamma(U_h,I)
 -\J_{N,1}^\gamma(U,I)
 &=
 2h\int_I\phi'\,U'\cdot\mathbf e_{N,j}\,dx
 +h^2|\mathbf e_{N,j}|^2\int_I|\phi'|^2\,dx+
 w_{j}\int_I
 \left[
 (u_j-u_{j+1}+h\phi)^\gamma
 -(u_j-u_{j+1})^\gamma
 \right]dx.
\end{aligned}
\]
We integrate the first term by parts on each half-interval using $\phi(\pm1)=0$, $U'(0^\pm)=0$, and \eqref{eq:common-active-gap-instability-equation},
\[\begin{aligned}
 2h\int_I\phi'\,U'\cdot\mathbf e_{N,j}\,dx
 &=-2h\int_I\phi\,U''\cdot\mathbf e_{N,j}\,dx=-\gamma h w_{j}
   \int_I(u_j-u_{j+1})^{\gamma-1}\phi\,dx.
\end{aligned}
\]
Moreover, 
\(
 \int_I|\phi'|^2\,dx
 =2\beta^2\int_0^1x^{2\beta-2}\,dx
 =\frac{2\beta^2}{2\beta-1}.
\)
Substituting these into the energy difference, we obtain
\[
\begin{aligned}
 \J_{N,1}^\gamma(U_h,I)
 -\J_{N,1}^\gamma(U,I)=
 \frac{2\beta^2|\mathbf e_{N,j}|^2}{2\beta-1}h^2+
 w_{j}\int_I  R_hdx,
\end{aligned}
\]
where
\[
 R_h=
 (u_j-u_{j+1}+h\phi)^\gamma
 -(u_j-u_{j+1})^\gamma
 -\gamma h(u_j-u_{j+1})^{\gamma-1}\phi.
\]
Concavity gives $R_h \leq 0$ a.e. Moreover for $a = u_j(1) - u_{j+1}(1)> 0$, one take $0 < h < a/2$ so that 
\(\phi \ge \frac{1}{2},\) \(u_j - u_{j+1} + h \phi \leq 2 h\) for all $0 < x < (h/a)^{1/\beta}$. Then on such region, Taylor gives
\[R_h \leq - \frac{\gamma (1 - \gamma)}{2}(u_j - u_{j+1} + h \phi)^{\gamma - 2} h^2 \phi^2 \leq  - \gamma (1 - \gamma) 2^{\gamma - 5} h^\gamma.\]
Since $\int_I R_h \lesssim -h^{\gamma + 1/\beta}$ and $\gamma + 1/\beta < 2$, the energy difference is negative for small $h>0.$
\end{proof}

\subsection{The one-sided cone with two active gaps}
\label{sec:type-ii-minimality}

We consider
$U_{\mathrm{II}}=\mathbf C^{\{1,2\}}x_+^\beta$,
where $\mathbf C^{\{1,2\}}$ solves
\eqref{eq:n3-two-gap-coefficients}.

\begin{theorem}[Type II]
\label{thm:type-ii-minimality}
For every $\mathbf w_\gamma\in(0,\infty)^2$,
$U_{\mathrm{II}}$ is not minimizing when
\[
 0<\gamma<\gamma_+:=\frac{14+4\sqrt6}{25}.
\]
Under \eqref{eq:weights_notcanonical}, 
$U_{\mathrm{II}}$ is minimizing for
$\gamma_{\mathrm{II}}<\gamma<1$ with \(\gamma_{\mathrm{II}}\) depending on $\rho(\cdot),\rho_* \in (1/2, 2)$.
\end{theorem}

\subsubsection{Non-minimality for $\gamma$ near $0$}

We do vertical perturbation away from $0$ to rule out middle range of $\gamma$.
\begin{lemma}
\label{lem:type-ii-stability-condition}
If $U_{\mathrm{II}}$ is minimizing, then
\begin{equation}\label{eq:type-ii-necessary-stability}
 2\gamma(1-\gamma)\sum_{i=1}^2
 w_{i}(c_i-c_{i+1})^{\gamma-2}(z_i-z_{i+1})^2
 \leq|\mathbf z|^2
 \qquad\forall\,\mathbf z\in\R^3.
\end{equation}
For every $\mathbf{w}_\gamma \in (0, \infty)^2$, this condition
fails when $\gamma_-<\gamma<\gamma_+$, where
$\gamma_\pm=(14\pm4\sqrt6)/25$.
\end{lemma}

\begin{proof}
Fix $\mathbf z\in\R^3$ and
$0\ne\phi\in C_c^\infty(0,1)$.
The membranes are strictly ordered on the support of
$\phi$, so $U_\ep= U_{\mathrm{II}}+\ep\phi\mathbf z$ is
admissible for small $|\ep|$. In view of
\[
\begin{aligned}
 \J_{3,1}^\gamma(U_\ep,(0,1))
 =\int_0^1\Bigg[
 &\sum_{i=1}^3
 \bigl(\beta c_i x^{\beta-1}+\ep\phi'z_i\bigr)^2+\sum_{i=1}^2w_{i}
 \bigl((c_i-c_{i+1})x^\beta
       +\ep\phi(z_i-z_{i+1})\bigr)^\gamma
 \Bigg]dx,
\end{aligned}
\]
one compute the second variation
\[
\begin{aligned}
 0
 &\leq
 \left.\frac{d^2}{d\ep^2}\right|_{\ep=0}
 \J_{3,1}^\gamma(U_\ep,(0,1))=2|\mathbf z|^2\int_0^1|\phi'|^2\,dx
 +\gamma(\gamma-1)
 \sum_{i=1}^2w_{i}(z_i-z_{i+1})^2
 \int_0^1
 \bigl((c_i-c_{i+1})x^\beta\bigr)^{\gamma-2}
 \phi^2\,dx.
\end{aligned}
\]
Since $\beta(\gamma-2)=-2$, rearranging gives
\[
 |\mathbf z|^2\int_0^1|\phi'|^2\,dx
 \geq
 \frac{\gamma(1-\gamma)}2
 \sum_{i=1}^2w_{i}
 (c_i-c_{i+1})^{\gamma-2}(z_i-z_{i+1})^2
 \int_0^1\frac{\phi^2}{x^2}\,dx.
\]
The sharp Hardy constant \cite{HLP52},
\(
 \inf_{0\ne\phi\in C_c^\infty(0,1)}
 \frac{\int_0^1|\phi'|^2\,dx}
      {\int_0^1\phi^2x^{-2}\,dx}
 =\frac14,
\)
therefore gives \eqref{eq:type-ii-necessary-stability}.

Moreover, \eqref{eq:n3-two-gap-coefficients} and
$c_1+c_2+c_3=0$ give
\[
\begin{aligned}
 w_{1}(c_1-c_2)^{\gamma-2}
 +w_{2}(c_2-c_3)^{\gamma-2}=
 \frac{\beta^2}{3}\left(
 4+\frac{c_2-c_3}{c_1-c_2}
  +\frac{c_1-c_2}{c_2-c_3}\right)
 \geq2\beta^2.
\end{aligned}
\]
Taking $\mathbf z=(1,-2,1)$ in
\eqref{eq:type-ii-necessary-stability} thus requires
\[
 6\gamma(1-\gamma)\beta^2\leq1,
 \qquad\text{equivalently}\qquad
 25\gamma^2-28\gamma+4\geq0.
\]
This fails for $\gamma_-<\gamma<\gamma_+$.
\end{proof}

We then change topology by delaying the separation of one gap until $\tau >0$ to rule out small $\gamma.$
\begin{lemma}
\label{lem:type-ii-delayed-branching}
For every $\mathbf w_\gamma\in(0,\infty)^2$,
$U_{\mathrm{II}}$ is not minimizing when $0<\gamma<3/4$.
\end{lemma}

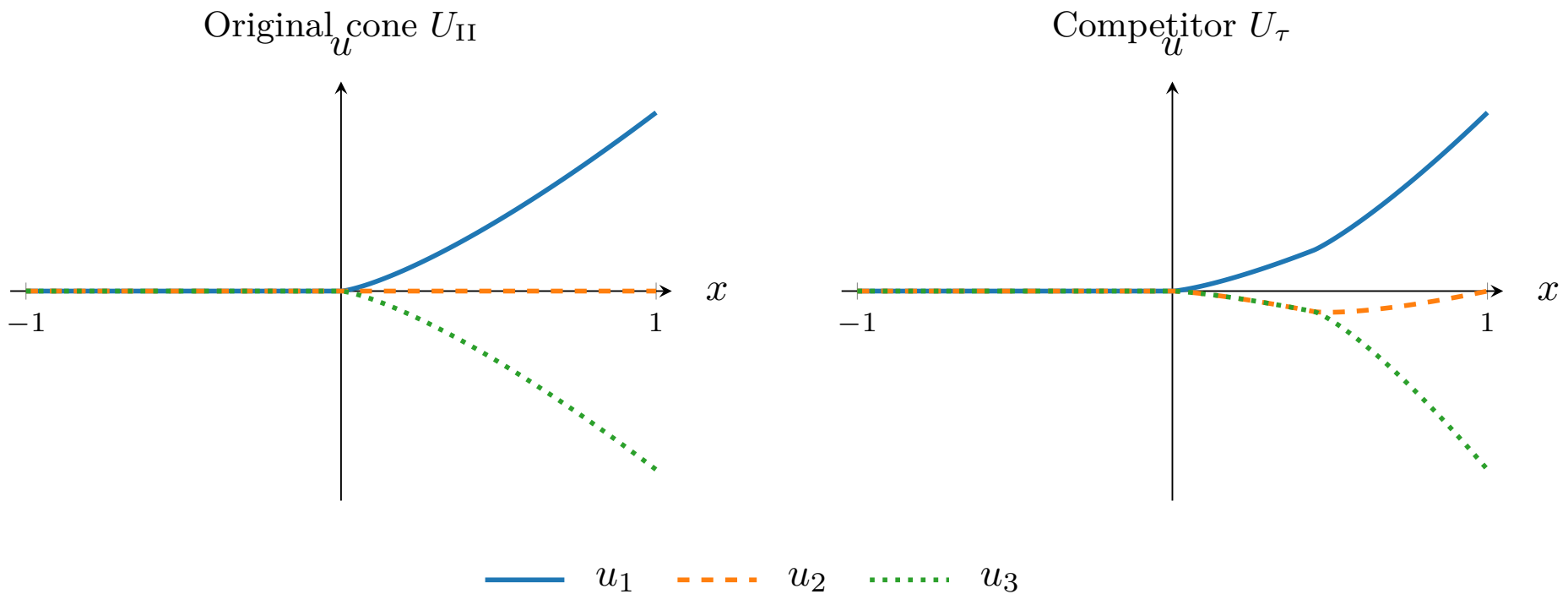
\begin{figure}[htbp]
\centering
\vspace{-6pt} 
\setlength{\abovecaptionskip}{1pt}
\definecolor{membraneone}{RGB}{31,119,180}
\definecolor{membranetwo}{RGB}{255,127,14}
\definecolor{membranethree}{RGB}{44,160,44}
\resizebox{0.75\textwidth}{!}{%
\begin{tikzpicture}
\pgfmathsetmacro{\betacon}{4/3}
\pgfmathsetmacro{\bcon}{pow(9/16,2/3)}
\pgfmathsetmacro{\taucon}{0.45}

\begin{groupplot}[
    group style={
        group size=2 by 1,
        horizontal sep=1.5cm
    },
    width=0.45\textwidth,
    height=0.32\textwidth,
    domain=-1:1,
    samples=401,
    axis lines=middle,
    xmin=-1.05, xmax=1.05,
    ymin=-0.8, ymax=0.8,
    xtick={-1,1},
    ytick=\empty,
    xlabel={$x$},
    ylabel={$u$},
    every axis x label/.style={
        at={(ticklabel* cs:1.03)}, anchor=west
    },
    every axis y label/.style={
        at={(ticklabel* cs:1.03)}, anchor=south
    },
    title style={font=\small},
    tick label style={font=\scriptsize},
    clip=false,
]

\nextgroupplot[title={Original cone $U_{\mathrm{II}}$}]
\addplot[membraneone, very thick]
    {\bcon*pow(max(x,0),\betacon)};
\addplot[membranetwo, very thick, dashed]
    {0};
\addplot[membranethree, very thick, dotted]
    {-\bcon*pow(max(x,0),\betacon)};

\nextgroupplot[title={Competitor $U_\tau$}]
\addplot[membraneone, very thick]
    {(\bcon/3)*(
        2*pow(max(x,0),\betacon)
        +pow(max(x-\taucon,0)/(1-\taucon),\betacon)
    )};
\addplot[membranetwo, very thick, dashed]
    {(\bcon/3)*(
        -pow(max(x,0),\betacon)
        +pow(max(x-\taucon,0)/(1-\taucon),\betacon)
    )};
\addplot[membranethree, very thick, dotted]
    {-(\bcon/3)*(
        pow(max(x,0),\betacon)
        +2*pow(max(x-\taucon,0)/(1-\taucon),\betacon)
    )};

\end{groupplot}
\begin{scope}[
    shift={($(group c1r1.south)!0.5!(group c2r1.south)
             +(0,-0.7cm)$)}
]
    \draw[membraneone, very thick]
        (-2.4,0) -- (-1.7,0);
    \node[anchor=west] at (-1.55,0) {$u_1$};

    \draw[membranetwo, very thick, dashed]
        (-0.7,0) -- (0,0);
    \node[anchor=west] at (0.15,0) {$u_2$};

    \draw[membranethree, very thick, dotted]
        (1,0) -- (1.7,0);
    \node[anchor=west] at (1.85,0) {$u_3$};
\end{scope}
\end{tikzpicture}%
}
\caption{Competitor in \Cref{lem:type-ii-delayed-branching}
for $\gamma=1/2$, $w_{1}=w_{2}=1$, and $\tau=0.45$.}
\label{fig:type-ii-delayed-separation}
\end{figure}

\begin{proof}
W.l.o.g. assume $c_1-c_2\geq c_2-c_3$.
On $I=(-1,1)$, take 
(see \Cref{fig:type-ii-delayed-separation})
\[
 U_\tau
 =(c_1-c_2)x_+^\beta\mathbf e_{3,1}
 +(c_2-c_3)
 \left(\frac{(x-\tau)_+}{1-\tau}\right)^\beta
 \mathbf e_{3,2}, \qquad 0 < \tau < 1.
\]
In view of \eqref{eq:n3-two-gap-coefficients-rewrite}, $U_\tau$ is admissible competitor for
$U_{\mathrm{II}}$ on $I$. Using  $|\mathbf e_{3,1}|^2=|\mathbf e_{3,2}|^2=2/3$
and $\mathbf e_{3,1}\cdot\mathbf e_{3,2}=1/3$, 
\[
\begin{aligned}
 |U_\tau'|^2
 &=\frac{2\beta^2}{3}\left[
 (c_1-c_2)^2x_+^{2\beta-2}
 +(c_2-c_3)^2
 \frac{(x-\tau)_+^{2\beta-2}}{(1-\tau)^{2\beta}}
 \right]+
 \frac{2\beta^2(c_1-c_2)(c_2-c_3)}3
 \frac{x_+^{\beta-1}(x-\tau)_+^{\beta-1}}
 {(1-\tau)^\beta},\\
 W_3^\gamma(U_\tau,\mathbf w_\gamma)
 &=w_{1}(c_1-c_2)^\gamma x_+^{\beta\gamma}
 +w_{2}(c_2-c_3)^\gamma
 \frac{(x-\tau)_+^{\beta\gamma}}{(1-\tau)^{\beta\gamma}}.
\end{aligned}
\]
Integrating all terms by change of
variables, and using $\beta\gamma+1=2\beta-1$, we obtain
\[
\begin{aligned}
 \int_I|U_\tau'|^2\,dx
 &=\frac{2\beta^2}{3(2\beta-1)}
 \left[
 (c_1-c_2)^2+\frac{(c_2-c_3)^2}{1-\tau}
 \right]+
 \frac{2\beta^2(c_1-c_2)(c_2-c_3)}3
 \int_0^1(\tau+(1-\tau)s)^{\beta-1}s^{\beta-1}\,ds,\\
 \int_I W_3^\gamma(U_\tau,\mathbf w_\gamma)\,dx
 &=\frac{
 w_{1}(c_1-c_2)^\gamma
 +(1-\tau)w_{2}(c_2-c_3)^\gamma}
 {2\beta-1}.
\end{aligned}
\]
 Note setting $\tau =0$ recovers those for $U_{\mathrm{II}}.$ Now we compute their energy difference
 \[
\begin{aligned}
 &\J_{3,1}^\gamma(U_\tau,I)
 -\J_{3,1}^\gamma(U_{\mathrm{II}},I)\\
 &=\frac{2\beta^2(c_2-c_3)^2}{3(2\beta-1)}
 \left(\frac1{1-\tau}-1\right)
 -\frac{\tau w_{2}(c_2-c_3)^\gamma}{2\beta-1}+
 \frac{2\beta^2(c_1-c_2)(c_2-c_3)}3
 \left[
 \int_0^1(\tau+(1-\tau)s)^{\beta-1}s^{\beta-1}\,ds
 -\frac1{2\beta-1}
 \right].
\end{aligned}
\]
By \eqref{eq:n3-two-gap-coefficients} and
$c_1+c_2+c_3=0$,
\[
\begin{aligned}
 w_{2}(c_2-c_3)^\gamma
 &=-\beta^2c_3(c_2-c_3)=\frac{\beta^2}{3}
 \left[(c_1-c_2)(c_2-c_3)+2(c_2-c_3)^2\right].
\end{aligned}
\]
Thus the terms containing $(c_2-c_3)^2$ combine.
Collecting the remaining mixed terms, we conclude
\[
\begin{aligned}
 F(\tau)
 &:=
 \frac{\J_{3,1}^\gamma(U_\tau,I)
       -\J_{3,1}^\gamma(U_{\mathrm{II}},I)}{\beta^2}=\frac{2(c_2-c_3)^2\tau^2}
 {3(2\beta-1)(1-\tau)}+
 \frac{(c_1-c_2)(c_2-c_3)}3
 \left[
 2\int_0^1(\tau+(1-\tau)s)^{\beta-1}s^{\beta-1}\,ds
 -\frac{2+\tau}{2\beta-1}
 \right].
\end{aligned}
\]
It remains to identify the range of $\gamma$ for which $F$ is strictly negative near $0.$
Note $F(0)=0$. Now differentiate
\[
\begin{aligned}
 F'(\tau)
 &=\frac{2(c_2-c_3)^2}{3(2\beta-1)}
 \frac{\tau(2-\tau)}{(1-\tau)^2}+
 \frac{(c_1-c_2)(c_2-c_3)}3
 \left[
 2(\beta-1)\int_0^1
 (\tau+(1-\tau)s)^{\beta-2}(1-s)s^{\beta-1}\,ds
 -\frac1{2\beta-1}
 \right].
\end{aligned}
\]
Then dominated convergence and
\(
 \int_0^1(1-s)s^{2\beta-3}\,ds
 =\frac1{(2\beta-2)(2\beta-1)}
\) imply
\(
 F'(0+)=0.
\) Differentiating again,
\[
\begin{aligned}
 F''(\tau)
 &=\frac{4(c_2-c_3)^2}{3(2\beta-1)(1-\tau)^3}+
 \frac{2(c_1-c_2)(c_2-c_3)}3(\beta-1)(\beta-2)
 \int_0^1
 (\tau+(1-\tau)s)^{\beta-3}(1-s)^2s^{\beta-1}\,ds.
\end{aligned}
\]
As $\tau\downarrow0$, monotone convergence gives
\[
\begin{aligned}
 \lim_{\tau\downarrow0}
 \int_0^1
 (\tau+(1-\tau)s)^{\beta-3}(1-s)^2s^{\beta-1}\,ds=\int_0^1(1-s)^2s^{2\beta-4}\,ds
 =
 \begin{cases}
 +\infty,&1<\beta\leq3/2,\\
 \displaystyle\frac{2}{(2\beta-3)(2\beta-2)(2\beta-1)},
 &3/2<\beta<2.
 \end{cases}
\end{aligned}
\]
Since $(\beta-1)(\beta-2)<0$,  $F''(\tau)\to-\infty$ when $1<\beta\leq3/2$.
For $3/2<\beta<8/5$, using $c_1-c_2\geq c_2-c_3$ gives
\[
\begin{aligned}
 \lim_{\tau\downarrow0}F''(\tau)
 &=\frac{2(c_2-c_3)}{3(2\beta-1)}
 \left[
 2(c_2-c_3)
 -(c_1-c_2)\frac{2-\beta}{2\beta-3}
 \right]\leq
 \frac{2(c_2-c_3)^2}{3(2\beta-1)}
 \left[2-\frac{2-\beta}{2\beta-3}\right]<0.
\end{aligned}
\]
Thus $F''(\tau)<0$ for small $\tau>0$
whenever $1<\beta<8/5$, equivalently $0<\gamma<3/4$.
\end{proof}

\subsubsection{Minimality for $\gamma$ near $1$}

We now assume \eqref{eq:weights_notcanonical} and prove
minimality near $\gamma=1$ by a calibration.

\begin{figure}[htbp]
\centering
\vspace{-6pt} 
\setlength{\abovecaptionskip}{1pt}
\includegraphics[width=0.7\textwidth]{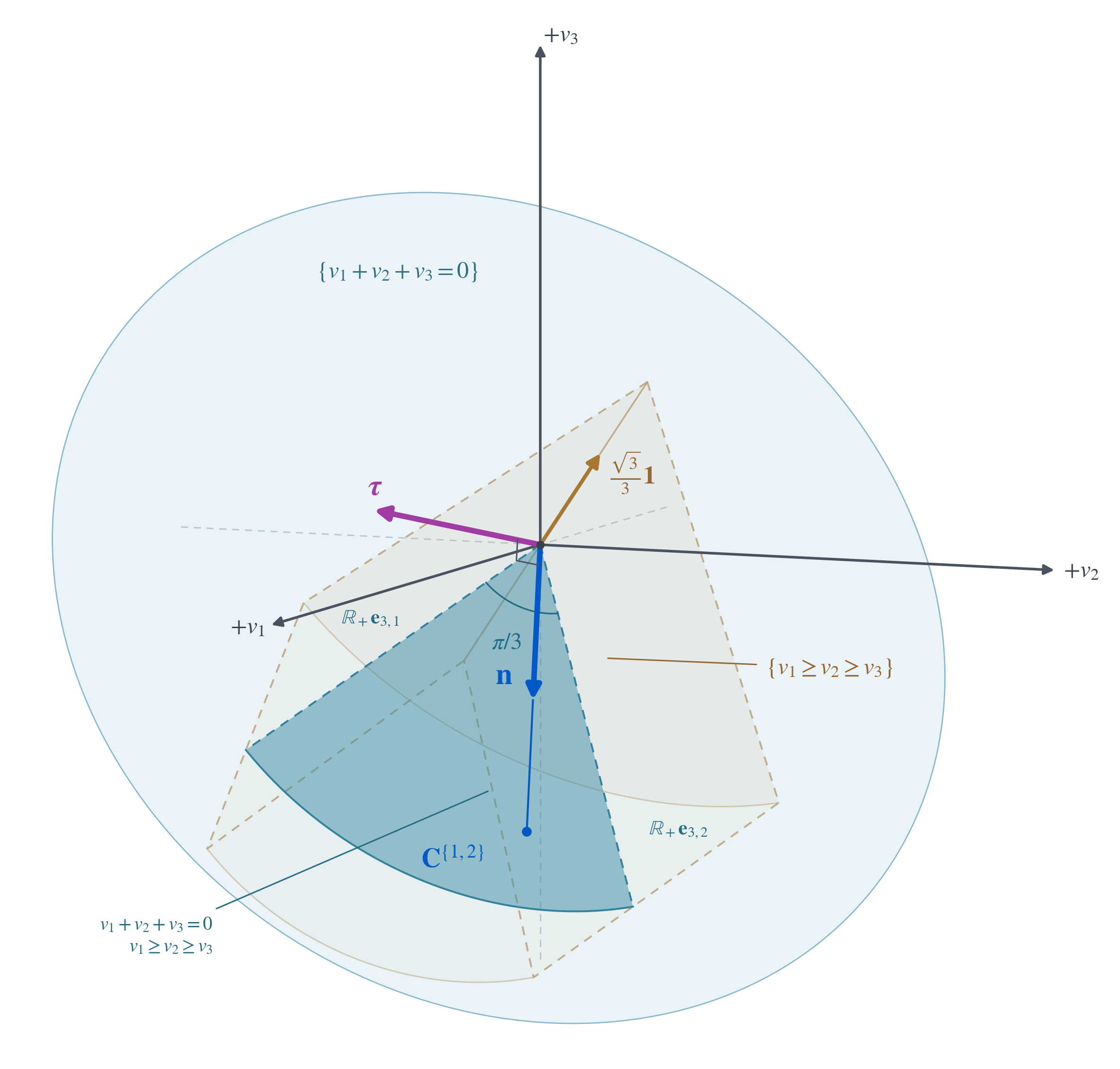}
\caption{The region $\{v_1\geq v_2\geq v_3\}$
and its intersection with the $\{\sum_{i=1}^3 v_i =0\}$.
The directions $\mathbf n,\boldsymbol\tau$ as in \eqref{eq:N3-n-tau} are orthonormal
within this plane, whose unit normal is
$\frac{\sqrt3}{3}\mathbf1$.}
\label{fig:n3-zero-mean-geometry}
\end{figure}

We choose orthonormal directions in the zero-mean plane $\{(v_1, v_2, v_3)\in \R^3 \mid \sum_{i=1}^3 v_i = 0\}$,
one along our Type~II cone and the other perpendicular to it (see \Cref{fig:n3-zero-mean-geometry}),
\begin{equation}\label{eq:N3-n-tau}
     \mathbf n=
 \frac{\mathbf C^{\{1,2\}}}{|\mathbf C^{\{1,2\}}|},
 \qquad
 \boldsymbol\tau=
 \frac{(c_2-c_3,c_3-c_1,c_1-c_2)}
 {\sqrt3|\mathbf C^{\{1,2\}}|}.
\end{equation}

Now for any ordered $V = (v_1, v_2, v_3)$, one decompose 
\begin{equation}\label{eq:N3-typeII-V-pq}
    V=p\mathbf n+q\boldsymbol\tau+h \frac{\sqrt{3}}{3}\mathbf1, \qquad p=\mathbf n\cdot V,\qquad
 q=\boldsymbol\tau\cdot V,\qquad
 h=\frac{v_1+v_2+v_3}{\sqrt{3}}.
\end{equation}
In particular, $U_{\mathrm{II}} $ in $(p, q, h)$-variables takes the simple form
\[p= p_0(x) := |\mathbf{C}^{\{1, 2\}}| x_+^\beta, \qquad q = 0, \qquad h = 0.\]

\begin{figure}[htbp]
\centering
\vspace{-6pt} 
\setlength{\abovecaptionskip}{1pt}
\includegraphics[width=0.45\textwidth]{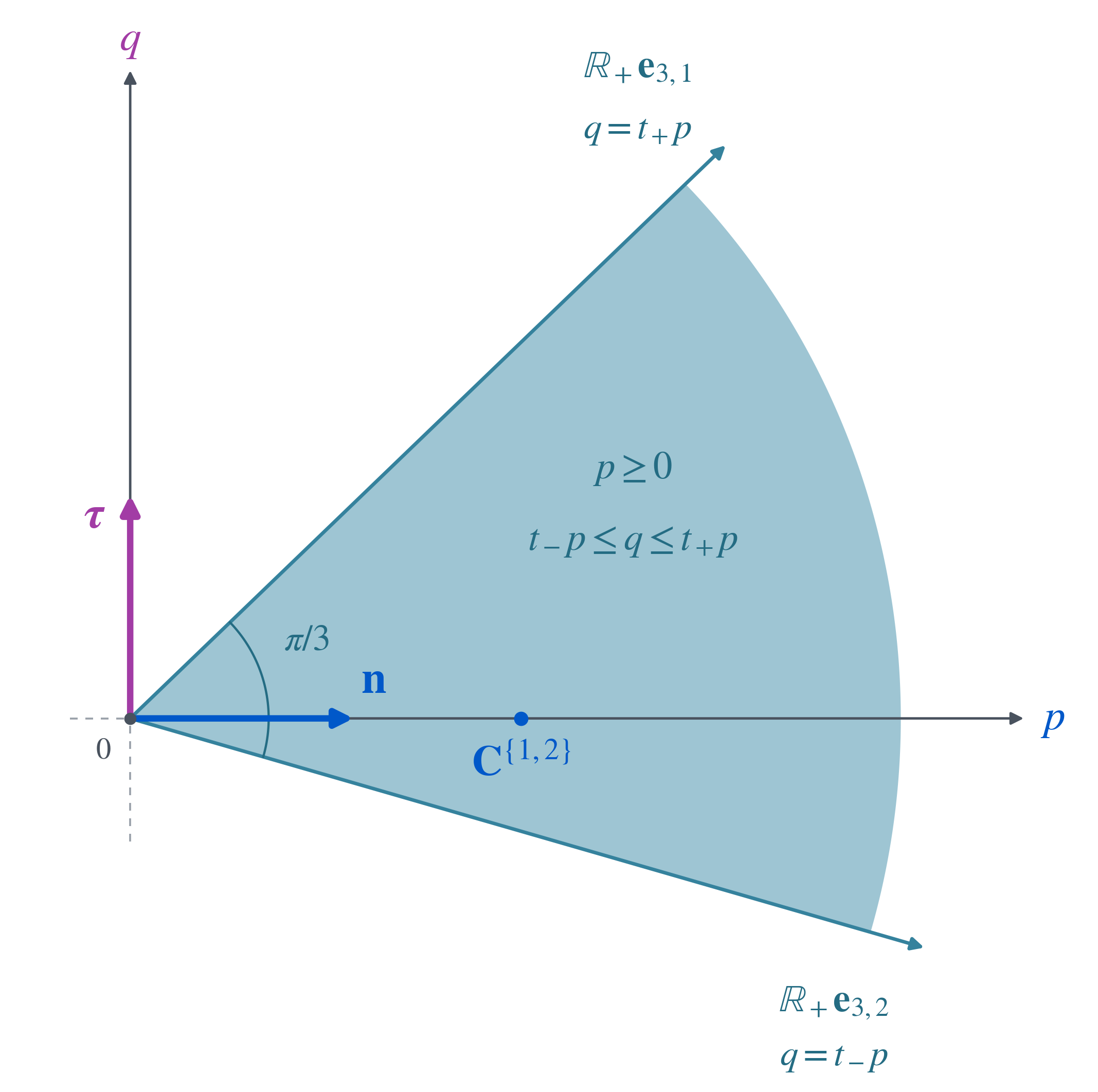}
\caption{\(\{v_1\geq v_2\geq v_3\}\cap \{\sum_{i=1}^3 v_i = 0\}\),
with $(p, q)$-coordinates and boundaries $q=t_\pm p$.}
\label{fig:n3-ordered-sector}
\end{figure}

The energy takes the form
\begin{equation}\label{eq:gap-functional-p-q-variables-3-1}
 \J_{3,1}^\gamma(V,I)
 =\int_I\left(
 |p'|^2+|q'|^2+|h'|^2
 +W_3^\gamma(p\mathbf n+q\boldsymbol\tau,\mathbf w_\gamma)
 \right)dx.
\end{equation}
We first rewrite the admissible region in  \((p,q,h)\)-variables. Using $c_1+c_2+c_3=0$ in \eqref{eq:N3-n-tau} gives
\begin{equation}\label{eq:gap-variables-p-q-variables-3-1}
    v_1-v_2
 =\frac{(c_1-c_2)p-\sqrt3c_3q}
 {|\mathbf C^{\{1,2\}}|}, \qquad v_2-v_3
 =\frac{(c_2-c_3)p-\sqrt3c_1q}
 {|\mathbf C^{\{1,2\}}|}.
\end{equation}
Since $c_1>0>c_3$, ordering is equivalent to
\begin{equation}
    \label{eq:gap-variables-t--t+-variables-3-1}
 p\geq0,\qquad t_-p\leq q\leq t_+p,
 \qquad
 t_-=\frac{c_1-c_2}{\sqrt3c_3}<0,\qquad
 t_+=\frac{c_2-c_3}{\sqrt3c_1}>0.
\end{equation} Therefore
\(\{v_1 \geq v_2 \geq v_3\} = \{p \geq 0 , \ t_-p\leq q\leq t_+p, \ h \in \R\}.\) In particular for $p>0$, the ratio $t=q/p$ determines
the direction in the zero-mean plane, with $t=0$
corresponding to $U_{\mathrm{II}}$; see \Cref{fig:n3-ordered-sector}. 

We seek a lower bound for the energy determined only
by the boundary values and attained by $U_{\mathrm{II}}$. For any $C^1$ function $\Psi=\Psi(p,q)$, using
\eqref{eq:gap-functional-p-q-variables-3-1} and completing
the square gives
\begin{equation}\label{eq:type-ii-explicit-calibration-motivation}
\begin{aligned}
 &\J_{3,1}^\gamma(V,(r,s))
 -\bigl[\Psi(p(s),q(s))-\Psi(p(r),q(r))\bigr]\\
 &=\int_r^s\left[
 \left|p'-\frac12\partial_p\Psi(p,q)\right|^2
 +\left|q'-\frac12\partial_q\Psi(p,q)\right|^2
 +|h'|^2
 +W_3^\gamma(p\mathbf n+q\boldsymbol\tau,\mathbf w_\gamma)
 -\frac14\left(
 |\partial_p\Psi(p,q)|^2+|\partial_q\Psi(p,q)|^2
 \right)
 \right]dx.
\end{aligned}
\end{equation}
The purpose of our calibration argument is to seek $\Psi$ satisfying
\begin{equation}\label{eq:type-ii-explicit-calibration-motivation-1}
    \frac14\left(
 |\partial_p\Psi(p,q)|^2+|\partial_q\Psi(p,q)|^2
 \right)
 \leq W_3^\gamma(p\mathbf n+q\boldsymbol\tau,\mathbf w_\gamma)
 \qquad\text{for }p\geq0,\quad t_-p\leq q\leq t_+p,
\end{equation}
so that \eqref{eq:type-ii-explicit-calibration-motivation}
gives the desired lower bound, and require
\begin{equation}\label{eq:type-ii-explicit-calibration-motivation-2}
     \partial_p\Psi(p_0(x),0)=2p_0'(x),
 \qquad
 \partial_q\Psi(p_0(x),0)=0
\end{equation}
so that the equality is attained by $U_{\mathrm{II}}$ with
$p_0(x)=|\mathbf C^{\{1,2\}}|x_+^\beta$; Indeed, since $q=h=0$ along $U_{\mathrm{II}}$ and
\eqref{eq:n3-two-gap-coefficients} gives
\(
 |p_0'|^2=W_3^\gamma(p_0\mathbf n,\mathbf w_\gamma),
\) \eqref{eq:type-ii-explicit-calibration-motivation-2} suffices for every term in the integrand of
\eqref{eq:type-ii-explicit-calibration-motivation}
to vanish.

Now we motivate the construction of $\Psi$ starting from \eqref{eq:type-ii-explicit-calibration-motivation-2}.
Since $\beta-1=\beta\gamma/2$, we have
\[
 p_0'
 =\beta|\mathbf C^{\{1,2\}}|x_+^{\beta-1}
 =\sqrt{K_\gamma}\,p_0^{\gamma/2}, \qquad  K_\gamma
 :=W_3^\gamma(\mathbf n,\mathbf w_\gamma)
 =\beta^2|\mathbf C^{\{1,2\}}|^{2-\gamma}.
\]
Thus \eqref{eq:type-ii-explicit-calibration-motivation-2} naturally prescribes
\[
\partial_p\Psi(p,0)=2\sqrt{K_\gamma}\,p^{\gamma/2},
 \qquad
 \partial_q\Psi(p,0)=0.
\]
Integrating the first identity and setting $\Psi(0, 0) =0 $ gives
\begin{equation}\label{eq:type-ii-explicit-calibration-motivation-4}
    \Psi(p,0)
 =\frac{2\sqrt{K_\gamma}}{1+\gamma/2}\,p^{1+\gamma/2}.
\end{equation}

 If one naively put $\Psi(p, q) = \frac{2\sqrt{K_\gamma}}{1+\gamma/2}\,p^{1+\gamma/2}$ independent of $q$, then
\begin{equation}\label{eq:type-ii-explicit-calibration-motivation-3}
    \frac14\left(
|\partial_p\Psi(p,q)|^2+|\partial_q\Psi(p,q)|^2
 \right)=K_\gamma p^\gamma.
\end{equation}
However, the potential in $(p, t)$-coordinates write as 
\[W_3^\gamma(V, \mathbf{w}_\gamma) = p^\gamma \sum_{i = 1}^2 w_{i}\left(n_i - n_{i+1} + t(\tau_i -\tau_{i+1})\right)^\gamma = p^\gamma W_3^\gamma(\mathbf{n} + t  \boldsymbol\tau, \mathbf{w}_\gamma) = K_\gamma p^\gamma f_{\gamma, \rho}(t),\]
where, using \eqref{eq:gap-variables-p-q-variables-3-1} and now actually imposing \eqref{eq:weights_notcanonical}
\begin{equation}\label{eq:N3-typeII-f}
\begin{aligned}
 f_{\gamma,\rho}(t)
 &:=\frac{
 W_3^\gamma(\mathbf n+t\boldsymbol\tau,\mathbf w_\gamma)}
 {W_3^\gamma(\mathbf n,\mathbf w_\gamma)}=\frac{
 w_{1}(c_1-c_2-\sqrt3c_3t)^\gamma
 +w_{2}(c_2-c_3-\sqrt3c_1t)^\gamma}
 {w_{1}(c_1-c_2)^\gamma
 +w_{2}(c_2-c_3)^\gamma}\\
 &=\frac{
 (c_1-c_2-\sqrt3c_3t)^\gamma
 +\rho(c_2-c_3-\sqrt3c_1t)^\gamma}
 {(c_1-c_2)^\gamma+\rho(c_2-c_3)^\gamma},
 \qquad t_-\leq t\leq t_+.
\end{aligned}
\end{equation}
It is easy to check that 
\begin{equation}\label{eq:N3-typeII-f-properties}
     f_{\gamma, \rho}(0) = 1, \qquad f_{\gamma, \rho}'(0) = 0, \qquad f_{\gamma, \rho} \in C^0([t_-, t_+]) \cap C^\infty((t_-, t_+)) \quad \text{is strictly concave.}
\end{equation}
Indeed, $f_{\gamma, \rho}(0)=1$ by definition, and using
\eqref{eq:n3-two-gap-coefficients} gives
\[
 f_{\gamma, \rho}'(0)
 =\frac{
 \nabla W_3^\gamma(\mathbf n,\mathbf w_\gamma)
 \cdot\boldsymbol\tau}{K_\gamma}
 =\gamma\,\mathbf n\cdot\boldsymbol\tau=0.
\] Regularity follows by the positive interior gaps, and differentiating twice gives
\begin{equation}\label{eq:N3-typeII-f-second-derivative}
     f_{\gamma, \rho}''(t)
 =-\frac{\gamma(1-\gamma)}{K_\gamma}
 \sum_{i=1}^2w_{i}
 (\tau_i-\tau_{i+1})^2
 \left[n_i-n_{i+1}
       +t(\tau_i-\tau_{i+1})\right]^{\gamma-2}<0.
\end{equation}
Hence $f_{\gamma, \rho}$ is strictly concave and attains its
unique maximum at $t=0$. In particular,
\[
 W_3^\gamma(p\mathbf n+q\boldsymbol\tau,\mathbf w_\gamma)
 =K_\gamma p^\gamma f_{\gamma, \rho}(q/p)
 <K_\gamma p^\gamma
 \qquad\text{for }p>0,\quad q\ne0.
\]
Therefore \eqref{eq:type-ii-explicit-calibration-motivation-3}
cannot satisfy
\eqref{eq:type-ii-explicit-calibration-motivation-1}
away from $q=0$.

We must add a correction to \eqref{eq:type-ii-explicit-calibration-motivation-4} that decreases
the squared gradient of $\Psi$ while preserving
\eqref{eq:type-ii-explicit-calibration-motivation-2}. The simplest correction preserving the value and both
first derivatives along $q=0$ is quadratic in $q$.
We also retain the homogeneity $1+\gamma/2$ in
\eqref{eq:type-ii-explicit-calibration-motivation-4},
so that the squared gradient has degree matching $W_3^\gamma$. We thus try
\[
 \Psi(p,q)
 =\sqrt{K_\gamma}\left[
 \frac{2}{1+\gamma/2}p^{1+\gamma/2}
 +\lambda q^2p^{\gamma/2-1}
 \right],
 \qquad p>0,
\]
where $\lambda>0$ will be chosen below. Writing $t=q/p$, we obtain
\[
\begin{aligned}
 \partial_p\Psi(p,q)
 &=\sqrt{K_\gamma}\left[
 2p^{\gamma/2}
 -\lambda(1-\gamma/2)q^2p^{\gamma/2-2}
 \right]=\sqrt{K_\gamma}\,p^{\gamma/2}
 \left[2-\lambda(1-\gamma/2)t^2\right],\\
 \partial_q\Psi(p,q)
 &=2\lambda\sqrt{K_\gamma}\,qp^{\gamma/2-1}
 =2\lambda\sqrt{K_\gamma}\,p^{\gamma/2}t.
\end{aligned}
\]
Consequently,
\[
\begin{aligned}
 \frac{
|\partial_p\Psi(p,q)|^2+|\partial_q\Psi(p,q)|^2}
 {4K_\gamma p^\gamma}
 &=
 \frac14\left[2-\lambda(1-\gamma/2)t^2\right]^2
 +\lambda^2t^2=
 1-\left[\lambda(1-\gamma/2)-\lambda^2\right]t^2
 +\frac{\lambda^2(1-\gamma/2)^2}{4}t^4.
\end{aligned}
\]
Thus the decrease contributed by $\partial_p\Psi$ is
linear in $\lambda$, whereas the additional contribution
of $\partial_q\Psi$ is quadratic in $\lambda$.
For $0<\lambda<1-\gamma/2$, the resulting $t^2$ term is negative.

At $\gamma=1$, the coefficient for $-t^2$ is
$\lambda/2-\lambda^2$, which is largest at $\lambda=1/4$.
This motivates the fixed choice
\begin{equation}\label{eq:type-ii-explicit-calibration}
 \Psi(p,q)
 =\sqrt{K_\gamma}\left[
 \frac{2}{1+\gamma/2}p^{1+\gamma/2}
 +\frac14q^2p^{\gamma/2-1}
 \right],
 \qquad p>0,
 \qquad \Psi(0,0)=0.
\end{equation}
The equalities
\eqref{eq:type-ii-explicit-calibration-motivation-2}
are preserved. Moreover, since $|q|\leq
\max\{|t_-|,|t_+|\}p \leq \sqrt{3} p$ in the ordered sector,
\[
 \Psi(p,q)=O(p^{1+\gamma/2}),
 \qquad
 |\partial_p\Psi(p,q)|+|\partial_q\Psi(p,q)|
 =O(p^{\gamma/2}),
\]
so $\Psi$ extends in $C^1$ to the
vertex, with both derivatives zero there.

For this choice of \eqref{eq:type-ii-explicit-calibration}, one in fact compute precisely 
\[
 \frac{
|\partial_p\Psi(p,q)|^2+|\partial_q\Psi(p,q)|^2}
 {4K_\gamma p^\gamma}
 =
 1-\left(\frac{1-\gamma/2}{4}-\frac1{16}\right)t^2
 +\frac{(1-\gamma/2)^2}{64}t^4.
\]
Therefore, by \eqref{eq:N3-typeII-f}, verifying
\eqref{eq:type-ii-explicit-calibration-motivation-1}
reduces to showing
\begin{equation}\label{eq:N3-typeII-key-inequality}
    1-\left(\frac{1-\gamma/2}{4}-\frac1{16}\right)t^2
 +\frac{(1-\gamma/2)^2}{64}t^4
 \leq f_{\gamma,\rho}(t),
 \qquad t_-\leq t\leq t_+, \qquad  \frac{1}{2}< \rho < 2.
\end{equation}
We establish this inequality
when $\gamma$ is sufficiently close to $1$.

\begin{lemma}\label{lem:N3-typeII-key-inequality}
For
$\rho=w_{2}/w_{1} \in (1/2, 2)$, set
\[
 \kappa_{\gamma,\rho}
 :=\max\left\{
 \frac{1-f_{\gamma,\rho}(t_-)}{t_-^2},
 \frac{1-f_{\gamma,\rho}(t_+)}{t_+^2}
 \right\}.
\]
Then \eqref{eq:N3-typeII-key-inequality} holds provided
\begin{equation}\label{eq:type-ii-weighted-calibration-condition}
 \kappa_{\gamma,\rho}
 \leq
 \frac{36-20\gamma-3\gamma^2}{256}.
\end{equation}
In particular, 
\eqref{eq:type-ii-weighted-calibration-condition}
holds for $\gamma$
sufficiently close to $1$ depending only on
$\rho_*$ and $\rho(\cdot).$
\end{lemma}
\begin{proof}

\emph{Step 1.}
We first show that the endpoint values defining
$\kappa_{\gamma,\rho}$ control $f_{\gamma,\rho}$ over the
interval \(t_-\leq t\leq t_+\).
By the formula for $f_{\gamma,\rho}''$ as in
\eqref{eq:N3-typeII-f-second-derivative},
the function $-f_{\gamma,\rho}''$ is convex.
Taylor and \eqref{eq:N3-typeII-f-properties} then give
\[
 \frac{1-f_{\gamma,\rho}(t)}{t^2}
 =-\int_0^1(1-s)f_{\gamma,\rho}''(st)\,ds,
 \qquad t\in(t_-,t_+)\setminus\{0\}.
\]
The quotient is therefore convex, with continuous
extension $-f_{\gamma,\rho}''(0)/2$ at $t=0$.
It is also continuous at $t_\pm$, so its maximum
is attained at an endpoint. Consequently,
\begin{equation*}
 1-\kappa_{\gamma,\rho}t^2
 \leq f_{\gamma,\rho}(t)\leq1,
 \qquad t_-\leq t\leq t_+.
\end{equation*}
Since $t^2\leq3$ throughout $[t_-,t_+]$,
\eqref{eq:type-ii-weighted-calibration-condition}
and rearranging prove \eqref{eq:N3-typeII-key-inequality}
\[
\begin{aligned}
 1-\frac{3-2\gamma}{16}t^2
 +\frac{(1-\gamma/2)^2}{64}t^4
 &\leq
 1-\left[
 \frac{3-2\gamma}{16}
 -\frac{3(1-\gamma/2)^2}{64}
 \right]t^2
 \leq1-\kappa_{\gamma,\rho}t^2
 \leq f_{\gamma,\rho}(t).
\end{aligned}
\]

\emph{Step 2.} We show \eqref{eq:type-ii-weighted-calibration-condition} as $\gamma \uparrow 1.$
Since
\(
 \frac{36-20\gamma-3\gamma^2}{256}\to\frac{13}{256},
\)
it suffices to show
$\kappa_{\gamma,\rho}\to0$.

Multiplying $\mathbf{w}_\gamma$ by constants leaves $t_\pm$, $f_{\gamma,\rho}$
and $\kappa_{\gamma,\rho}$ unchanged, so we work with
\(
 w_{1}=1,\ w_{2}=\rho,
\)
and $\mathbf C^{\{1,2\}}$ denotes the corresponding solution to
\eqref{eq:n3-two-gap-coefficients}. We first show 
\begin{equation}\label{eq:gap-variables-3-1-convergence-C*}
    \mathbf C^{\{1,2\}}
 \longrightarrow
 \mathbf C_*=\frac14(1,\rho_*-1,-\rho_*)
 \qquad\text{as }\gamma\uparrow1.
\end{equation}

Indeed, by \Cref{lem:eighted-critical-coefficient-system},
$\mathbf C^{\{1,2\}}$ minimizes
\[
 F_{\{1,2\}}(c)
 =\frac12|c|^2
 -\frac{\beta^{-2}}{\gamma}
 \left[
 (c_1-c_2)^\gamma
 +\rho(c_2-c_3)^\gamma
 \right], \qquad \{c_1\geq c_2\geq c_3\}.
\]
For $C$ independent of $\gamma$ near $1$ changing from line to line,
\begin{align*}
    F_{\{1, 2\}}(c)&\geq \frac{1}{2}|c|^2 - \frac{\beta^{-2}}{\gamma} \left(1 + c_1 - c_2 + \rho(1 + c_2 -c_3)\right)  \geq \frac{1}{2}|c|^2 - \frac{\beta^{-2}}{\gamma}\left(1 + \rho\right)\left(1 + \sqrt{2}|c|\right)\geq \frac{1}{4}|c|^2 - C.
\end{align*}
Combining with minimality, \(
 \frac14|\mathbf C^{\{1,2\}}|^2-C
 \leq F_{\{1,2\}}(\mathbf C^{\{1,2\}})
 \leq F_{\{1,2\}}(0)=0
\)
gives \(|\mathbf C^{\{1,2\}}| \leq C.\) Moreover, the functional $F_{\{1, 2\}}$ converge locally uniformly on compact subsets of $\{c_1 \geq c_2 \ge c_3\}$ to
\[
\begin{aligned}
 F_*(c)
 &=\frac12|c|^2
 -\frac14\left[
 (c_1-c_2)+\rho_*(c_2-c_3)
 \right]=\frac12|c|^2
 -\frac14(1,\rho_*-1,-\rho_*)\cdot c=\frac12|c-\mathbf C_*|^2-\frac12|\mathbf C_*|^2,
\end{aligned}
\]
where $\mathbf C_*$
is strictly ordered and the unique minimizer of $F_*$; In particular,
\[
\begin{aligned}
 \sup_{\substack{|c|\leq C\\c_1\geq c_2\geq c_3}}
 |F_{\{1,2\}}(c)-F_*(c)|
 &\leq C\left[
 (1-\gamma)+|\rho(\gamma)-\rho_*|
 \right]\longrightarrow0.
\end{aligned}
\]
Using minimality $F_{\{1, 2\}}(\mathbf{C}^{\{1, 2\}}) \leq F_{\{1, 2\}}(\mathbf{C}_*),$ one conclude \eqref{eq:gap-variables-3-1-convergence-C*} via
\begin{align*}
 \frac12|\mathbf C^{\{1,2\}}-\mathbf C_*|^2
 &=F_*(\mathbf C^{\{1,2\}})-F_*(\mathbf C_*)\leq
 2\sup_{\substack{|c|\leq C\\c_1\geq c_2\geq c_3}}
 \left|F_{\{1,2\}}(c)-F_*(c)\right|
 \longrightarrow0.
\end{align*}

Now Using \eqref{eq:gap-variables-3-1-convergence-C*} and \eqref{eq:gap-variables-t--t+-variables-3-1},
we obtain
\begin{equation}\label{eq:gap-variables-t--t+-variables-3-1-conbvergence}
\begin{aligned}
 t_-
 &=\frac{c_1-c_2}{\sqrt3c_3}
 \longrightarrow
 \frac{\frac14-\frac14(\rho_*-1)}
 {\sqrt3\left(-\frac14\rho_*\right)}
 =-\frac{2-\rho_*}{\sqrt3\rho_*},&
 t_+
 &=\frac{c_2-c_3}{\sqrt3c_1}
 \longrightarrow
 \frac{\frac14(\rho_*-1)+\frac14\rho_*}
 {\sqrt3\,\frac14}
 =\frac{2\rho_*-1}{\sqrt3}.
\end{aligned}
\end{equation}
Both limits are nonzero. At the limiting exponent and weights, substituting
$\mathbf C_*$ into \eqref{eq:N3-typeII-f} gives
\[
\begin{aligned}
 f_{1,\rho_*}(t)
 &=\frac{
 (2-\rho_*+\sqrt3\rho_*t)
 +\rho_*(2\rho_*-1-\sqrt3t)}
 {(2-\rho_*)+\rho_*(2\rho_*-1)}
 =1,
 \qquad
 -\frac{2-\rho_*}{\sqrt3\rho_*}
 \leq t\leq\frac{2\rho_*-1}{\sqrt3}.
\end{aligned}
\]
The formula \eqref{eq:N3-typeII-f} is continuous in $\gamma, \rho, c, t$ up to $\gamma = 1.$
Thus \eqref{eq:gap-variables-3-1-convergence-C*}
and \eqref{eq:gap-variables-t--t+-variables-3-1-conbvergence}
give
\[
 \kappa_{\gamma,\rho(\gamma)}
 =\max\left\{
 \frac{1-f_{\gamma,\rho}(t_-)}{t_-^2},
 \frac{1-f_{\gamma,\rho}(t_+)}{t_+^2}
 \right\}
 \longrightarrow0.
\]
\end{proof}

Now we conclude the proof for Type~II cones.

\begin{proof}[Proof of \Cref{thm:type-ii-minimality}]
\Cref{lem:type-ii-stability-condition,lem:type-ii-delayed-branching}
give non-minimality for every $0<\gamma<\gamma_+$. For minimality, using \eqref{eq:weights_notcanonical} and \Cref{lem:N3-typeII-key-inequality}, one may choose
$\gamma_{\mathrm{II}} =\gamma_{\mathrm{II}}(\rho_*,\rho(\cdot))\in(\gamma_+,1)$ such that
\eqref{eq:N3-typeII-key-inequality} holds whenever
$\gamma_{\mathrm{II}}<\gamma<1$, and the calibration
\eqref{eq:type-ii-explicit-calibration} satisfies
\eqref{eq:type-ii-explicit-calibration-motivation-1}
and \eqref{eq:type-ii-explicit-calibration-motivation-2}. Let $I=(r,s)$ be bounded interval and $V$ an admissible competitor for $U_{\mathrm{II}}$ in $I$.
Write $V$ as in \eqref{eq:N3-typeII-V-pq}. Since it agrees with $(p_0,0,0)$ at $\{r, s\}$, \eqref{eq:type-ii-explicit-calibration-motivation}
gives
\[
\begin{aligned}
 \J_{3,1}^\gamma(V,I)
 &\geq \Psi(p(s),q(s))-\Psi(p(r),q(r))=\Psi(p_0(s),0)-\Psi(p_0(r),0)=\J_{3,1}^\gamma(U_{\mathrm{II}},I),
\end{aligned}
\]
where the last equality follows from
\eqref{eq:type-ii-explicit-calibration-motivation-2}
and $|p_0'|^2=W_3^\gamma(p_0\mathbf n,\mathbf w_\gamma)$.
Hence $U_{\mathrm{II}}$ is minimizing.
\end{proof}
\begin{remark}\label{rmk:typeIIrho=1}
For weights $\rho\equiv1$, we give an explicit choice
of $\gamma_{\mathrm{II}}$ in \Cref{thm:type-ii-minimality}.
By symmetry, $c_2=0$ and $c_1=-c_3$, so
\eqref{eq:N3-typeII-f} gives
\[
 t_\pm=\pm\frac1{\sqrt3},
 \qquad
 f_{\gamma,1}(t)
 =\frac{(1+\sqrt3\,t)^\gamma+(1-\sqrt3\,t)^\gamma}{2}.
\]
Hence $f_{\gamma,1}(t_\pm)=2^{\gamma-1}$ and
$\kappa_{\gamma,1}=3(1-2^{\gamma-1})$.
By Step~1 of \Cref{lem:N3-typeII-key-inequality},
$(1-f_{\gamma,1}(t))/t^2$ attains its maximum at the endpoints,
where $t^2=1/3$. Thus
\eqref{eq:N3-typeII-key-inequality} holds precisely when
\[
 3(1-2^{\gamma-1})+\frac{(1-\gamma/2)^2}{192}
 \leq\frac{3-2\gamma}{16},
 \qquad\text{equivalently}\qquad
2304\,2^{\gamma-1}\geq2164+92\gamma+\gamma^2.
\]
We may therefore take $\gamma_{\mathrm{II}}\approx0.968$
to be the unique root of
\(
 2304\,2^{\gamma-1}=2164+92\gamma+\gamma^2.
\)
\end{remark}

\subsection{The two-sided cone with different active gaps}
\label{sec:type-iv-minimality}

We next consider
\[
 U_{\mathrm{IV}}(x)
 =\mathbf C^{\{1\}}x_+^\beta
  +\mathbf C^{\{2\}}x_-^\beta,
 \qquad
 \mathbf C^{\{k\}}=a_k\mathbf e_{3,k},
\]
under \eqref{eq:weights_notcanonical}, with
$a_k$ given by \eqref{eq:single-gap-ak}.
\begin{theorem}[Type IV]
\label{thm:type-iv-minimality}
For every $\mathbf w_\gamma\in(0,\infty)^2$,
$U_{\mathrm{IV}}$ is not minimizing for
$0<\gamma\leq2/3$.
Under \eqref{eq:weights_notcanonical}, $U_{\mathrm{IV}}$ is minimizing for
$\gamma_{\mathrm{IV}}<\gamma<1$ with $\gamma_{\mathrm{IV}}$ depending on
$\rho_*, \rho(\cdot) \in (1/2, 2).$
\end{theorem}

\subsubsection{Non-minimality for $\gamma$ near $0$}

We move the two opening points past each other.
On the overlap, the two scalar profiles have
derivatives of opposite signs, whose interaction decreases the Dirichlet energy.

\begin{lemma}
\label{lem:type-iv-overlap-instability}
For every $\mathbf w_\gamma\in(0,\infty)^2$,
$U_{\mathrm{IV}}$ is not minimizing when
$0<\gamma\leq2/3$.
\end{lemma}
\begin{proof}
On $I=(-1,1)$, take (see \Cref{fig:type-iv-overlap})
\[
 U_\tau
 =a_1\left(\frac{(x+\tau_1)_+}{1+\tau_1}\right)^\beta
       \mathbf e_{3,1}
  +a_2\left(\frac{(\tau_2-x)_+}{1+\tau_2}\right)^\beta
       \mathbf e_{3,2}, \qquad \tau_1=\frac{\tau}{a_1},\
 \tau_2=\frac{\tau}{a_2}, \qquad  0<\tau<\min\{a_1,a_2\}.
\]
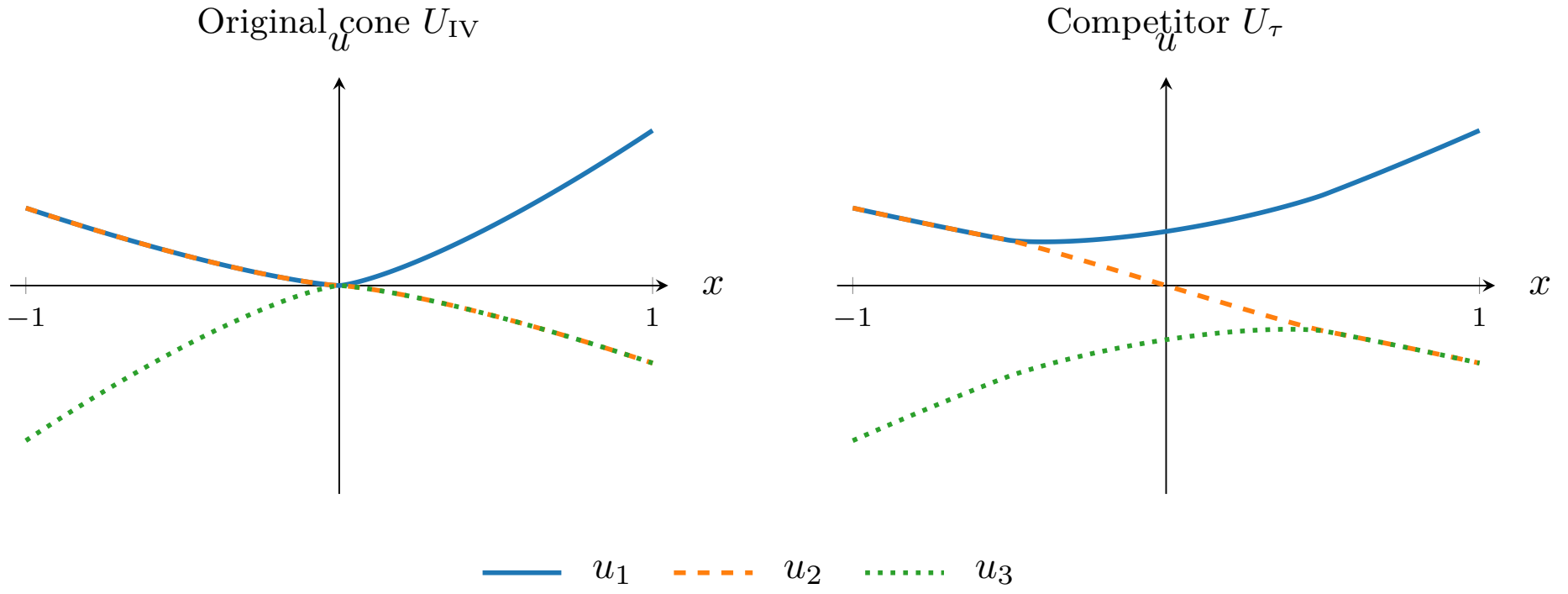
\begin{figure}[htbp]
\centering
\vspace{-6pt} 
\setlength{\abovecaptionskip}{1pt}
\definecolor{membraneone}{RGB}{31,119,180}
\definecolor{membranetwo}{RGB}{255,127,14}
\definecolor{membranethree}{RGB}{44,160,44}
\resizebox{0.75\textwidth}{!}{%
\begin{tikzpicture}
\pgfmathsetmacro{\betacon}{4/3}
\pgfmathsetmacro{\aonecon}{pow(27/32,2/3)}
\pgfmathsetmacro{\atwocon}{pow(27/32,2/3)}
\pgfmathsetmacro{\taucon}{0.45}
\pgfmathsetmacro{\tauonecon}{\taucon/\aonecon}
\pgfmathsetmacro{\tautwocon}{\taucon/\atwocon}

\begin{groupplot}[
    group style={
        group size=2 by 1,
        horizontal sep=1.5cm
    },
    width=0.45\textwidth,
    height=0.32\textwidth,
    domain=-1:1,
    samples=401,
    axis lines=middle,
    xmin=-1.05, xmax=1.05,
    ymin=-0.8, ymax=0.8,
    xtick={-1,1},
    ytick=\empty,
    xlabel={$x$},
    ylabel={$u$},
    every axis x label/.style={
        at={(ticklabel* cs:1.03)}, anchor=west
    },
    every axis y label/.style={
        at={(ticklabel* cs:1.03)}, anchor=south
    },
    title style={font=\small},
    tick label style={font=\scriptsize},
    clip=false,
]

\nextgroupplot[title={Original cone $U_{\mathrm{IV}}$}]
\addplot[membraneone, very thick]
    {(2/3)*\aonecon*pow(max(x,0),\betacon)
     +(1/3)*\atwocon*pow(max(-x,0),\betacon)};
\addplot[membranetwo, very thick, dashed]
    {-(1/3)*\aonecon*pow(max(x,0),\betacon)
     +(1/3)*\atwocon*pow(max(-x,0),\betacon)};
\addplot[membranethree, very thick, dotted]
    {-(1/3)*\aonecon*pow(max(x,0),\betacon)
     -(2/3)*\atwocon*pow(max(-x,0),\betacon)};

\nextgroupplot[title={Competitor $U_\tau$}]
\addplot[membraneone, very thick]
    {(2/3)*\aonecon
        *pow(max(x+\tauonecon,0)/(1+\tauonecon),\betacon)
     +(1/3)*\atwocon
        *pow(max(\tautwocon-x,0)/(1+\tautwocon),\betacon)};
\addplot[membranetwo, very thick, dashed]
    {-(1/3)*\aonecon
        *pow(max(x+\tauonecon,0)/(1+\tauonecon),\betacon)
     +(1/3)*\atwocon
        *pow(max(\tautwocon-x,0)/(1+\tautwocon),\betacon)};
\addplot[membranethree, very thick, dotted]
    {-(1/3)*\aonecon
        *pow(max(x+\tauonecon,0)/(1+\tauonecon),\betacon)
     -(2/3)*\atwocon
        *pow(max(\tautwocon-x,0)/(1+\tautwocon),\betacon)};

\end{groupplot}

\begin{scope}[
    shift={($(group c1r1.south)!0.5!(group c2r1.south)
             +(0,-0.7cm)$)}
]
    \draw[membraneone, very thick]
        (-2.4,0) -- (-1.7,0);
    \node[anchor=west] at (-1.55,0) {$u_1$};

    \draw[membranetwo, very thick, dashed]
        (-0.7,0) -- (0,0);
    \node[anchor=west] at (0.15,0) {$u_2$};

    \draw[membranethree, very thick, dotted]
        (1,0) -- (1.7,0);
    \node[anchor=west] at (1.85,0) {$u_3$};
\end{scope}
\end{tikzpicture}%
}
\caption{Competitor in \Cref{lem:type-iv-overlap-instability}
for $\gamma=1/2$, $w_{1}=w_{2}=1$, and $\tau=0.45$.}
\label{fig:type-iv-overlap}
\end{figure}

$U_\tau$ is ordered and agrees with
$U_{\mathrm{IV}}$ at $\pm1$. 
Using $|\mathbf e_{3,k}|^2=2/3$,
$\mathbf e_{3,1}\cdot\mathbf e_{3,2}=1/3$, we obtain
\[
\begin{aligned}
 |U_\tau'|^2
 &=\frac{2\beta^2}{3}\left[
 a_1^2\frac{(x+\tau_1)_+^{2\beta-2}}{(1+\tau_1)^{2\beta}}
 +a_2^2\frac{(\tau_2-x)_+^{2\beta-2}}{(1+\tau_2)^{2\beta}}
 \right]-
 \frac{2\beta^2a_1a_2}{3}
 \frac{(x+\tau_1)_+^{\beta-1}(\tau_2-x)_+^{\beta-1}}
 {(1+\tau_1)^\beta(1+\tau_2)^\beta},\\
 W_3^\gamma(U_\tau,\mathbf w_\gamma)
 &=w_{1}a_1^\gamma
   \frac{(x+\tau_1)_+^{\beta\gamma}}{(1+\tau_1)^{\beta\gamma}}
 +w_{2}a_2^\gamma
   \frac{(\tau_2-x)_+^{\beta\gamma}}{(1+\tau_2)^{\beta\gamma}}.
\end{aligned}
\]
Integrating all terms by change of
variables, and using $\beta\gamma+1=2\beta-1$, we obtain
\[
\begin{aligned}
 \int_I|U_\tau'|^2\,dx
 &=\frac{2\beta^2}{3(2\beta-1)}
 \left[
 \frac{a_1^2}{1+\tau_1}+\frac{a_2^2}{1+\tau_2}
 \right]-
 \frac{2\beta^2a_1a_2}{3}
 \frac{(\tau_1+\tau_2)^{2\beta-1}}
 {(1+\tau_1)^\beta(1+\tau_2)^\beta}
 \int_0^1s^{\beta-1}(1-s)^{\beta-1}\,ds,\\
 \int_IW_3^\gamma(U_\tau,\mathbf w_\gamma)\,dx
 &=\frac{
 (1+\tau_1)w_{1}a_1^\gamma
 +(1+\tau_2)w_{2}a_2^\gamma}
 {2\beta-1}.
\end{aligned}
\]
Setting $\tau=0$ recovers the corresponding expressions
for $U_{\mathrm{IV}}$. Thus their energy difference is
\[
\begin{aligned}
 &\J_{3,1}^\gamma(U_\tau,I)
 -\J_{3,1}^\gamma(U_{\mathrm{IV}},I)\\
 &=\frac{2\beta^2}{3(2\beta-1)}
\sum_{k=1}^2a_k^2\left(\frac1{1+\tau_k}-1\right)
 +\frac1{2\beta-1}
 \sum_{k=1}^2\tau_k w_{k}a_k^\gamma-
 \frac{2\beta^2a_1a_2}{3}
 \frac{(\tau_1+\tau_2)^{2\beta-1}}
 {(1+\tau_1)^\beta(1+\tau_2)^\beta}
 \int_0^1s^{\beta-1}(1-s)^{\beta-1}\,ds.
\end{aligned}
\]
By \eqref{eq:single-gap-ak},
$w_{k}a_k^\gamma=\frac23\beta^2a_k^2$.
Hence the terms containing $a_k^2$ combine. Collecting these terms, we conclude
\[
\begin{aligned}
 F(\tau)
 &:=
 \frac{3}{2\beta^2}
 \left[
 \J_{3,1}^\gamma(U_\tau,I)
 -\J_{3,1}^\gamma(U_{\mathrm{IV}},I)
 \right]=\frac1{2\beta-1}
 \left[
 \frac{a_1^2\tau_1^2}{1+\tau_1}
 +\frac{a_2^2\tau_2^2}{1+\tau_2}
 \right]-
 \frac{a_1a_2(\tau_1+\tau_2)^{2\beta-1}}
 {(1+\tau_1)^\beta(1+\tau_2)^\beta}
 \int_0^1s^{\beta-1}(1-s)^{\beta-1}\,ds\\
 &=\frac{\tau^2}{2\beta-1}
 \left[
 \frac1{1+\tau/a_1}+\frac1{1+\tau/a_2}
 \right]-
 \frac{a_1a_2\,\tau^{2\beta-1}
       (a_1^{-1}+a_2^{-1})^{2\beta-1}}
 {(1+\tau/a_1)^\beta(1+\tau/a_2)^\beta}
 \int_0^1s^{\beta-1}(1-s)^{\beta-1}\,ds.
\end{aligned}
\]
It remains to identify the range of $\gamma$ for which $F$ is strictly negative near $0.$
Note $F(0)=0$. Now differentiate
\[
\begin{aligned}
 F'(\tau)
 &=\frac1{2\beta-1}
 \sum_{k=1}^2
 \frac{\tau(2+\tau/a_k)}{(1+\tau/a_k)^2}-
 \frac{a_1a_2(a_1^{-1}+a_2^{-1})^{2\beta-1}
       \tau^{2\beta-2}}
 {(1+\tau/a_1)^\beta(1+\tau/a_2)^\beta}
 \left[
 2\beta-1-\beta\tau\sum_{k=1}^2\frac1{a_k+\tau}
 \right]
 \int_0^1s^{\beta-1}(1-s)^{\beta-1}\,ds.
\end{aligned}
\]
Since $\beta>1$, both terms tend to zero as
$\tau\downarrow0$, so $F'(0+)=0$.
Differentiating again,
\[
\begin{aligned}
 &F''(\tau)
 =\frac2{2\beta-1}
 \sum_{k=1}^2\frac1{(1+\tau/a_k)^3}\\
 &\quad-
 \frac{a_1a_2(a_1^{-1}+a_2^{-1})^{2\beta-1}
       \tau^{2\beta-3}}
 {(1+\tau/a_1)^\beta(1+\tau/a_2)^\beta}
 \int_0^1s^{\beta-1}(1-s)^{\beta-1}\,ds\cdot 
 \Bigg[
 \left(
 2\beta-1-\beta\tau\sum_{k=1}^2\frac1{a_k+\tau}
 \right)^2
 -(2\beta-1)
 +\beta\tau^2\sum_{k=1}^2\frac1{(a_k+\tau)^2}
 \Bigg].
\end{aligned}
\]
As $\tau\downarrow0$, the expression in brackets
converges to
\(
 (2\beta-1)^2-(2\beta-1)
 =(2\beta-1)(2\beta-2)>0.
\) For $1<\beta<3/2$, 
$\tau^{2\beta-3}$ diverges so
\(
 \lim_{\tau\downarrow0}F''(\tau)=-\infty.
\)
At $\beta=3/2$, using
$\int_0^1\sqrt{s(1-s)}\,ds=\pi/8$,
\[
\begin{aligned}
 \lim_{\tau\downarrow0}F''(\tau)
 &=2-\frac{\pi}{4}a_1a_2
       \left(\frac1{a_1}+\frac1{a_2}\right)^2=2-\frac{\pi}{4}
       \left(\frac{a_1}{a_2}+\frac{a_2}{a_1}+2\right)
 \leq2-\pi<0.
\end{aligned}
\]
Therefore $F''(\tau)<0$ for small $\tau>0$
whenever $1<\beta\leq3/2$, equivalently
$0<\gamma\leq2/3$.
\end{proof}

\subsubsection{Minimality for $\gamma$ near $1$}
For the next two lemmas, normalize
\[
  w_1=1,\qquad w_2=\rho\in(1/2,2),\qquad
  a_k=\left(\frac{3w_k}{2\beta^2}\right)^{1/(2-\gamma)},
  \qquad k=1,2.
\]

First consider the profile $a_kx^\beta$ on a half-interval.
We quantify the possible energy decrease from lifting
its value at the origin, while keeping the other
endpoint fixed.
\begin{lemma}\label{lem:type-iv-scalar-bridge}
For $2/3<\gamma<1$, there exists $\ep_\gamma\geq0$,
depending only on $\gamma$, with $\ep_\gamma\to0$ as
$\gamma\uparrow1$, such that
\begin{equation}\label{eq:type-iv-scalar-bridge}
\begin{aligned}
\int_0^L\left(\frac23|v'|^2+w_kv^\gamma\right)\,dx
 -\int_0^L\left(
   \frac23|(a_kx^\beta)'|^2+w_k(a_kx^\beta)^\gamma
 \right)\,dx\geq-\ep_\gamma h^{1+\gamma/2}, \qquad k = 1, 2,
\end{aligned}
\end{equation}
for every $L>0$, $h\geq0$, and nonnegative
$v\in H^1(0,L)$ satisfying
$v(0)=h$ and $v(L)=a_kL^\beta$.
\end{lemma}

\begin{figure}[htbp]
\centering
\vspace{-6pt} 
\setlength{\abovecaptionskip}{1pt}
\definecolor{membraneone}{RGB}{31,119,180}
\definecolor{membranetwo}{RGB}{255,127,14}
\definecolor{membranethree}{RGB}{44,160,44}
\resizebox{0.85\textwidth}{!}{%
\begin{tikzpicture}
\pgfmathsetmacro{\gammabridge}{3/4}
\pgfmathsetmacro{\betabridge}{2/(2-\gammabridge)}
\pgfmathsetmacro{\abridge}
    {pow(3/(2*\betabridge^2),1/(2-\gammabridge))}
\pgfmathsetmacro{\rbridge}{3/4}
\pgfmathsetmacro{\Tbridge}{2}
\pgfmathsetmacro{\Lbridge}{\rbridge*\Tbridge}
\pgfmathsetmacro{\Hbridge}
    {\abridge*pow(\rbridge,\betabridge)}
\pgfmathsetmacro{\endbridge}
    {\abridge*pow(\Lbridge,\betabridge)}

\begin{groupplot}[
    group style={
        group size=2 by 1,
        horizontal sep=2cm
    },
    width=0.43\textwidth,
    height=0.33\textwidth,
    samples=201,
    axis lines=left,
    xmin=0, ymin=0,
    title style={font=\small},
    tick label style={font=\scriptsize},
    label style={font=\small},
    every axis x label/.style={
        at={(ticklabel* cs:1.03)}, anchor=west
    },
    legend style={
        at={(0.5,-0.23)}, anchor=north,
        draw=none, font=\small,
        /tikz/every even column/.append style={column sep=8pt}
    },
    legend columns=2,
    clip=false
]

\nextgroupplot[
    title={Original variables},
    domain=0:\Lbridge,
    xmax=1.12*\Lbridge,
    ymax=1.18*\endbridge,
    xlabel={$x$},
    xtick={\rbridge,\Lbridge},
    xticklabels={$\widetilde h$,$L$},
    ytick={\Hbridge,\endbridge},
    yticklabels={$h$,$a_1L^\beta$}
]

\draw[gray!50,densely dotted]
    (axis cs:0,\Hbridge)
    -- (axis cs:\rbridge,\Hbridge)
    -- (axis cs:\rbridge,0);
\draw[gray!50,densely dotted]
    (axis cs:0,\endbridge)
    -- (axis cs:\Lbridge,\endbridge)
    -- (axis cs:\Lbridge,0);

\addplot[membraneone,very thick]
    {\abridge*pow(x,\betabridge)};
\addlegendentry{$a_1x^\beta$}

\addplot[membranetwo,very thick,dashed]
    {\Hbridge*pow(
        1-0.3*(x/\rbridge)+0.4*(x/\rbridge)^2,
        \betabridge
    )};
\addlegendentry{$v(x)$}

\addplot[
    only marks,mark=*,mark size=1.8pt,
    membraneone,forget plot
] coordinates {(0,0)};
\addplot[
    only marks,mark=*,mark size=1.8pt,
    membranetwo,forget plot
] coordinates {(0,\Hbridge)};
\addplot[
    only marks,mark=*,mark size=1.8pt,
    black,forget plot
] coordinates {(\Lbridge,\endbridge)};

\nextgroupplot[
    title={Normalized variables},
    domain=0:\Tbridge,
    xmax=2.24, ymax=2.36,
    xlabel={$t$},
    xtick={2}, xticklabels={$T$},
    ytick={1,2}, yticklabels={$1$,$T$}
]

\addplot[
    name path=bridgeReference,
    draw=none,forget plot
] {x};
\addplot[
    name path=bridgeCompetitor,
    draw=none,forget plot
] {1-0.3*x+0.4*x^2};
\addplot[
    membranethree!12,draw=none,forget plot
] fill between[
    of=bridgeReference and bridgeCompetitor,
    soft clip={domain=0:1.25}
];

\draw[gray!50,densely dotted]
    (axis cs:0,2) -- (axis cs:2,2) -- (axis cs:2,0);

\addplot[membraneone,very thick] {x};
\addlegendentry{$t$}

\addplot[membranetwo,very thick,dashed]
    {1-0.3*x+0.4*x^2};
\addlegendentry{$z(t)$}

\draw[<->,membranethree,thick]
    (axis cs:0.45,0.45) -- (axis cs:0.45,0.946)
    node[midway,left=3pt,font=\small] {$\zeta(t)$};

\addplot[
    only marks,mark=*,mark size=1.8pt,
    membraneone,forget plot
] coordinates {(0,0)};
\addplot[
    only marks,mark=*,mark size=1.8pt,
    membranetwo,forget plot
] coordinates {(0,1)};
\addplot[
    only marks,mark=*,mark size=1.8pt,
    black,forget plot
] coordinates {(2,2)};

\end{groupplot}
\end{tikzpicture}%
}
\caption{The scalar comparison in
\eqref{eq:type-iv-scalar-bridge}, for $k=1$ and $\gamma=3/4$.}
\label{fig:type-iv-scalar-bridge}
\end{figure}
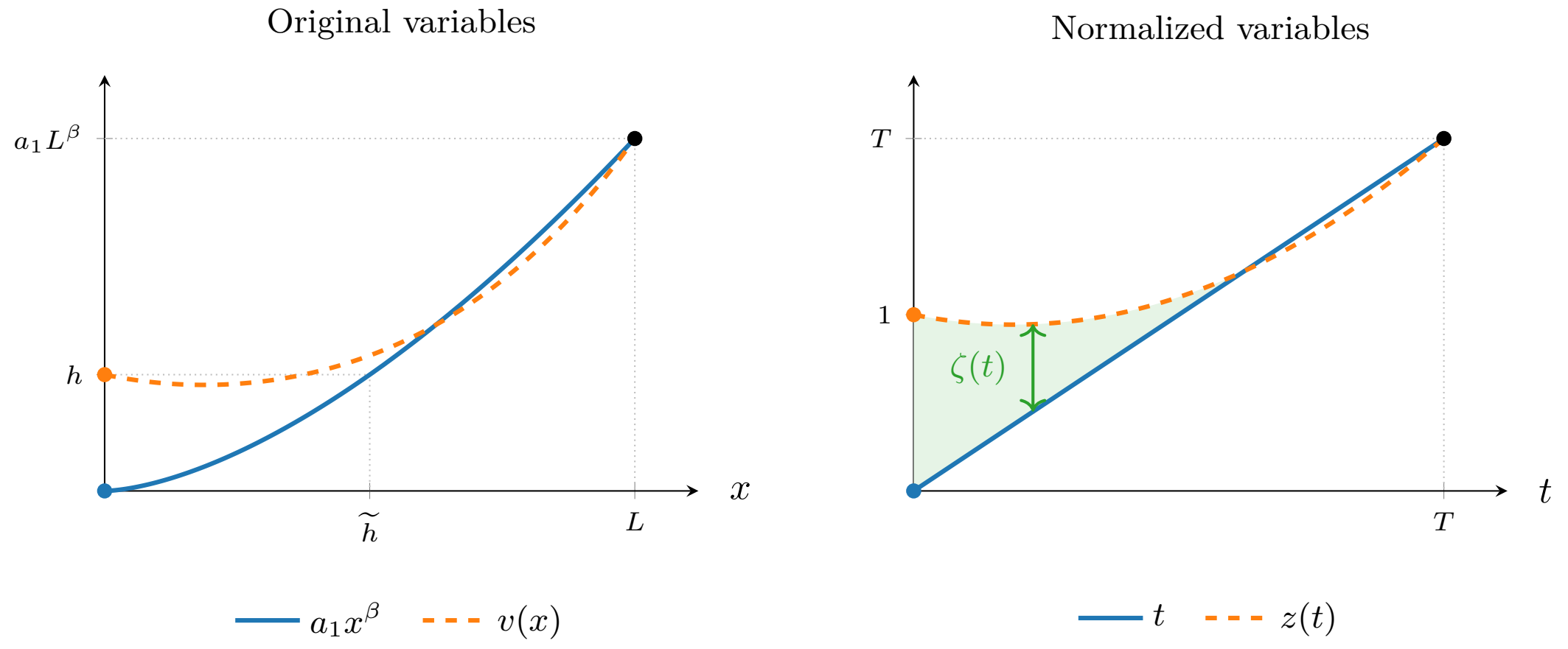

\begin{proof}
\emph{Step 1.}
Fix $k\in\{1,2\}$. $h=0$ follows from
\eqref{eq:single-gap-minimizers-comparison-energy-0}.
For $h>0$, first assume $v>0$ on $[0,L]$. We normalize
\[
 \widetilde h=(h/a_k)^{1/\beta},
 \qquad T=L/\widetilde h,
 \qquad
 z(t)=\left(\frac{v(\widetilde ht)}h\right)^{1/\beta}.
\]
Then $z(0)=1$, $z(T)=T$, and in particular, positivity ensures $z\in H^1(0,T)$.
Using $\beta\gamma=2\beta-2$ and
\eqref{eq:single-gap-ak},
\[
\int_0^L\left(\frac23|v'|^2+w_kv^\gamma\right)\,dx
 =
 \beta\sqrt{\frac{2w_k}{3}}\,h^{1+\gamma/2}
 \int_0^Tz^{\beta\gamma}(|z'|^2+1)\,dt.
\]
Under the same normalization, $a_kx^\beta$ becomes $z(t)=t$.
Thus, denoting the energy difference in
\eqref{eq:type-iv-scalar-bridge} by $\Delta E$,
\begin{equation}
\label{eq:type-iv-scalar-bridge-midstep-estimate-0}
\begin{aligned}
 \frac{\Delta E}
 {\beta\sqrt{2w_k/3}\,h^{1+\gamma/2}}
 &=
 \int_0^Tz^{\beta\gamma}(|z'|^2+1)\,dt
 -\frac{2T^{1+\beta\gamma}}{1+\beta\gamma}=
 \int_0^Tz^{\beta\gamma}|z'-1|^2\,dt
 -\frac{2}{1+\beta\gamma}.
\end{aligned}
\end{equation}

\emph{Step 2.}
Set $\zeta(t)=(z(t)-t)_+$, so that
$\zeta(0)=1$ and $\zeta(T)=0$. Then
\[
 \int_0^Tz^{\beta\gamma}|z'-1|^2\,dt
 \geq
\int_0^T(t+\zeta)^{\beta\gamma}|\zeta'|^2\,dt.
\]
To bound the latter uniformly in $T$, we seek $B(t,s)$
satisfying
\begin{equation}
\label{eq:type-iv-scalar-bridge-midstep-estimate}
 (t+\zeta)^{\beta\gamma}|\zeta'|^2
 +\frac{d}{dt}B(t,\zeta(t))\geq0,
 \qquad B(t,0)=0.
\end{equation}
Integrating would then give
\begin{equation}
\label{eq:type-iv-scalar-bridge-midstep-estimate-2}
\int_0^T(t+\zeta)^{\beta\gamma}|\zeta'|^2\,dt
 \geq B(0,1)-B(T,0)=B(0,1).
\end{equation}
To choose $B$, we complete the square,
\[
\begin{aligned}
 (t+\zeta)^{\beta\gamma}|\zeta'|^2
 +\frac{d}{dt}B(t,\zeta(t))=
 (t+\zeta)^{\beta\gamma}
 \left|\zeta'
 +\frac{\p_sB(t,\zeta)}
 {2(t+\zeta)^{\beta\gamma}}\right|^2
 +\p_tB(t,\zeta)
 -\frac{|\p_sB(t,\zeta)|^2}
 {4(t+\zeta)^{\beta\gamma}}.
\end{aligned}
\]
Thus it suffices to find $B$ with
\[
 \p_tB(t,s)\geq
 \frac{|\p_sB(t,s)|^2}{4(t+s)^{\beta\gamma}},
 \qquad t,s\geq0,\quad t+s>0.
\]
To this end, define
\begin{equation}\label{eq:typeIVcalibration}
 B(t,s)=2(\beta\gamma-1)
 \int_0^s r(t+r)^{\beta\gamma-1}\,dr.
\end{equation}
Then $B(t,0)=0$. Since $1<\beta\gamma<2$,
\[
\begin{aligned}
 \p_sB(t,s)
 &=2(\beta\gamma-1)s(t+s)^{\beta\gamma-1},\\
 \p_tB(t,s)
 &=2(\beta\gamma-1)^2
   \int_0^s r(t+r)^{\beta\gamma-2}\,dr\geq
 (\beta\gamma-1)^2s^2(t+s)^{\beta\gamma-2}
 =\frac{|\p_sB(t,s)|^2}{4(t+s)^{\beta\gamma}}.
\end{aligned}
\]
Therefore \eqref{eq:type-iv-scalar-bridge-midstep-estimate}
holds, with
\(
 B(0,1)=\frac{2(\beta\gamma-1)}{1+\beta\gamma}.
\)
Substituting \eqref{eq:type-iv-scalar-bridge-midstep-estimate-2} into
\eqref{eq:type-iv-scalar-bridge-midstep-estimate-0} gives
\[
 \Delta E\geq
 -2\beta\sqrt{\frac{2w_k}{3}}\,
 \frac{2-\beta\gamma}{1+\beta\gamma}
 h^{1+\gamma/2}.
\]
Since $w_k\leq2$, this proves
\eqref{eq:type-iv-scalar-bridge} with
\begin{equation}
\label{eq:type-iv-scalar-bridge-epsilon}
 \ep_\gamma
 =\frac{4\beta}{\sqrt3}
   \frac{2-\beta\gamma}{1+\beta\gamma}
 =\frac{32}{\sqrt3}\,
   \frac{1-\gamma}{4-\gamma^2}
 \downarrow0
 \qquad\text{as }\gamma\uparrow1.
\end{equation}

\emph{Step 3.}
We remove the positivity assumption. For $h>0$, set
\[
 \ell(x)=h\left(1-\frac{x}{L}\right)
 +a_kL^\beta\frac{x}{L},
 \qquad
 v_\delta=(1-\delta)v+\delta\ell,
 \qquad 0<\delta<1.
\]
Then $v_\delta$ has the same endpoint values as $v$ and
\(
 v_\delta\geq\delta\min\{h,a_kL^\beta\}>0.
\)
Apply \eqref{eq:type-iv-scalar-bridge} to $v_\delta$.
Since $v_\delta\to v$ in $H^1(0,L)$ and uniformly on
$[0,L]$, sending $\delta\downarrow0$ proves \eqref{eq:type-iv-scalar-bridge} for all admissible $v.$
\end{proof}

Now we consider the single-gap cone on a half-interval.
We quantify the possible energy change from opening
both gaps at the origin, while keeping the other
endpoint fixed. 
\begin{lemma}\label{lem:type-iv-one-sided-bridge}
Assume \eqref{eq:weights_notcanonical}. There exists $c=c(\rho_*)>0$ s.t. for $\gamma$
sufficiently close to $1$ depending only on
$\rho_*$ and $\rho(\cdot)$, and $\ep_\gamma$ given by
\Cref{lem:type-iv-scalar-bridge}, one has
\begin{equation}\label{eq:type-iv-one-sided-bridge}
  \J_{3,1}^\gamma(V,(0,L))
  -\J_{3,1}^\gamma(x^\beta\mathbf C^{\{k\}},(0,L))
  \geq
  -\ep_\gamma(h_1+h_2)^{1+\gamma/2}
  +c\,h_{3-k}^{1+\gamma/2},
  \qquad k=1,2,
\end{equation}
for every $L > 0$, $h_1, h_2 \geq 0$, and ordered
$V\in H^1((0,L);\R^3)$ satisfying \((v_i-v_{i+1})(0) = h_i\) and
$V(L)=L^\beta\mathbf C^{\{k\}}$.
\end{lemma}
\begin{figure}[htbp]
\centering
\vspace{-6pt} 
\setlength{\abovecaptionskip}{1pt}
\definecolor{membraneone}{RGB}{31,119,180}
\definecolor{membranetwo}{RGB}{255,127,14}
\definecolor{membranethree}{RGB}{44,160,44}
\resizebox{0.85\textwidth}{!}{%
\begin{tikzpicture}
\pgfmathsetmacro{\gammabridge}{3/4}
\pgfmathsetmacro{\betabridge}{2/(2-\gammabridge)}
\pgfmathsetmacro{\abridge}
    {pow(3/(2*\betabridge^2),1/(2-\gammabridge))}
\pgfmathsetmacro{\honebridge}{0.30*\abridge}
\pgfmathsetmacro{\htwobridge}{0.45*\abridge}
\pgfmathsetmacro{\vonezero}
    {(2*\honebridge+\htwobridge)/3}
\pgfmathsetmacro{\vtwozero}
    {(-\honebridge+\htwobridge)/3}
\pgfmathsetmacro{\vthreezero}
    {(-\honebridge-2*\htwobridge)/3}

\begin{groupplot}[
    group style={
        group size=2 by 1,
        horizontal sep=1.8cm
    },
    width=0.43\textwidth,
    height=0.33\textwidth,
    domain=0:1,
    samples=201,
    axis lines=middle,
    xmin=-0.24, xmax=1.13,
    ymin=-0.36, ymax=0.51,
    xtick={1}, xticklabels={$L$},
    ytick=\empty,
    xlabel={$x$},
    every axis x label/.style={
        at={(ticklabel* cs:1.03)}, anchor=west
    },
    title style={font=\small},
    tick label style={font=\scriptsize},
    legend style={
        at={(0.5,-0.14)}, anchor=north,
        draw=none, font=\small,
        /tikz/every even column/.append style={column sep=5pt}
    },
    legend columns=3,
    clip=false
]

\nextgroupplot[
    title={Single-gap cone $x^\beta\mathbf C^{\{1\}}$}
]
\draw[gray!45,densely dotted]
    (axis cs:1,{-\abridge/3})
    -- (axis cs:1,{2*\abridge/3});

\addplot[membraneone,very thick]
    {(2/3)*\abridge*pow(x,\betabridge)};
\addlegendentry{$u_1$}

\addplot[membranetwo,very thick,dashed]
    {-(1/3)*\abridge*pow(x,\betabridge)};
\addlegendentry{$u_2$}

\addplot[membranethree,very thick,dotted]
    {-(1/3)*\abridge*pow(x,\betabridge)};
\addlegendentry{$u_3$}

\addplot[
    black,only marks,mark=*,mark size=1.5pt,forget plot
] coordinates {
    (0,0)
    (1,{2*\abridge/3})
    (1,{-\abridge/3})
};

\nextgroupplot[title={Ordered competitor $V$}]
\draw[gray!45,densely dotted]
    (axis cs:1,{-\abridge/3})
    -- (axis cs:1,{2*\abridge/3});

% v_1-v_2=a_1 x^beta+h_1(1-x),
% v_2-v_3=h_2(1-x), and v_1+v_2+v_3=0.
\addplot[membraneone,very thick]
    {(2*\abridge*pow(x,\betabridge)
       +(2*\honebridge+\htwobridge)*(1-x))/3};
\addlegendentry{$v_1$}

\addplot[membranetwo,very thick,dashed]
    {(-\abridge*pow(x,\betabridge)
       +(-\honebridge+\htwobridge)*(1-x))/3};
\addlegendentry{$v_2$}

\addplot[membranethree,very thick,dotted]
    {(-\abridge*pow(x,\betabridge)
       +(-\honebridge-2*\htwobridge)*(1-x))/3};
\addlegendentry{$v_3$}

\draw[gray!50,densely dotted]
    (axis cs:-0.12,\vonezero) -- (axis cs:0,\vonezero);
\draw[gray!50,densely dotted]
    (axis cs:-0.12,\vtwozero) -- (axis cs:0,\vtwozero);
\draw[gray!50,densely dotted]
    (axis cs:-0.12,\vthreezero) -- (axis cs:0,\vthreezero);

\draw[<->,thick]
    (axis cs:-0.12,\vtwozero)
    -- (axis cs:-0.12,\vonezero)
    node[midway,left=3pt,font=\small] {$h_1$};
\draw[<->,thick]
    (axis cs:-0.12,\vthreezero)
    -- (axis cs:-0.12,\vtwozero)
    node[midway,left=3pt,font=\small] {$h_2$};

\addplot[
    membraneone,only marks,mark=*,
    mark size=1.5pt,forget plot
] coordinates {(0,\vonezero)};
\addplot[
    membranetwo,only marks,mark=*,
    mark size=1.5pt,forget plot
] coordinates {(0,\vtwozero)};
\addplot[
    membranethree,only marks,mark=*,
    mark size=1.5pt,forget plot
] coordinates {(0,\vthreezero)};
\addplot[
    black,only marks,mark=*,mark size=1.5pt,forget plot
] coordinates {
    (1,{2*\abridge/3})
    (1,{-\abridge/3})
};

\end{groupplot}
\end{tikzpicture}%
}
\caption{The comparison in
\eqref{eq:type-iv-one-sided-bridge}, for $k=1$,
$\ell=2$, $\gamma=3/4$, and $\rho=3/2$.}
\label{fig:type-iv-one-sided-openings}
\end{figure}
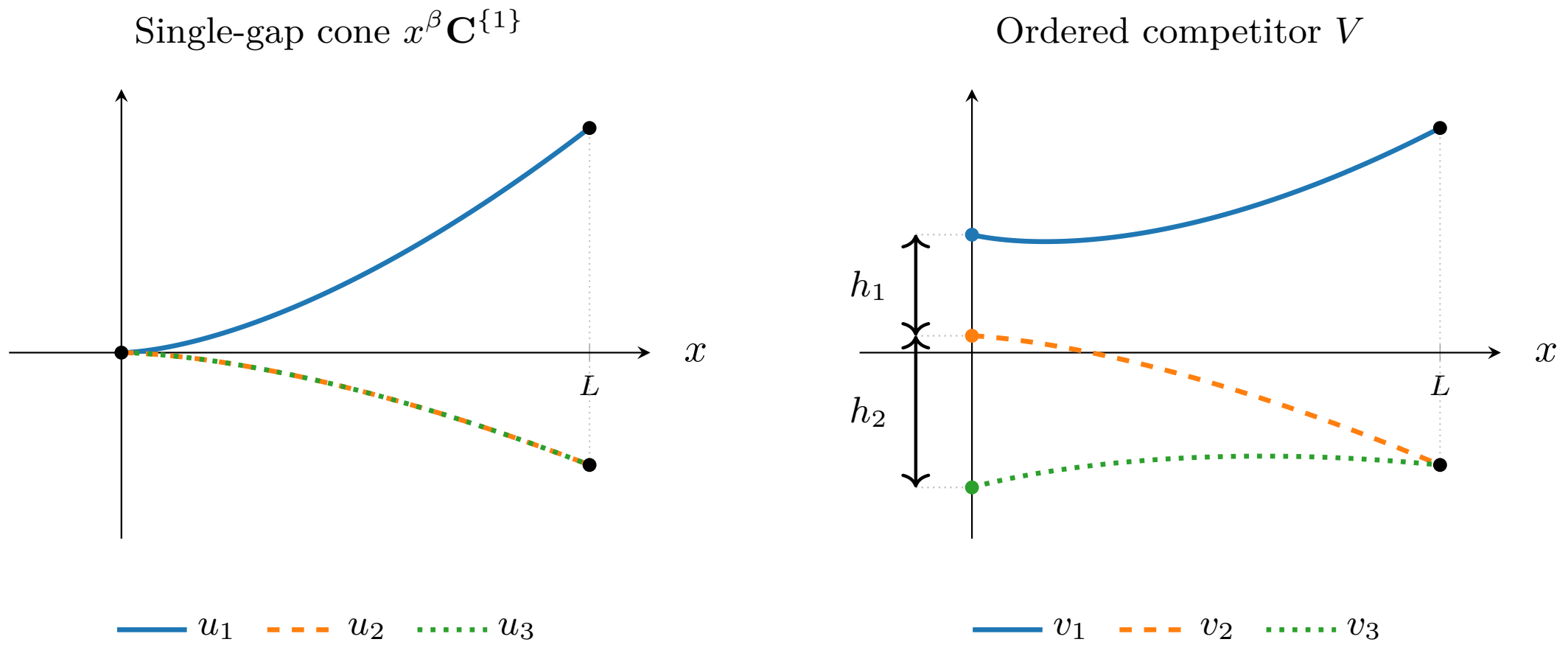
\begin{proof}
Fix $k\in\{1,2\}$ and set $\ell=3-k$.
Decompose as in \Cref{fig:type-iv-one-sided-projection}
\[
 V=p\mathbf e_{3,k}
 +q\left(\mathbf e_{3,\ell}-\frac12\mathbf e_{3,k}\right)
 +\frac h{\sqrt3}\mathbf1, \qquad 
 p=\frac{\mathbf e_{3,k}\cdot V}{|\mathbf e_{3,k}|^2},
 \qquad q=v_\ell-v_{\ell+1},
 \qquad h=\frac{v_1+v_2+v_3}{\sqrt3}.
\]
\begin{figure}[htbp]
\centering
\vspace{-6pt} 
\setlength{\abovecaptionskip}{1pt}
\definecolor{bridgeplane}{RGB}{62,139,164}
\definecolor{bridgeprojection}{RGB}{31,119,180}
\definecolor{bridgeresidual}{RGB}{159,60,171}
\definecolor{bridgemean}{RGB}{178,126,49}
\resizebox{0.65\textwidth}{!}{%
\begin{tikzpicture}[
    x={(-0.80cm,-0.232cm)},
    y={(1.95cm,-0.09518cm)},
    z={(0cm,2.09275cm)},
    >=Latex,
    line cap=round,
    line join=round,
    every node/.style={font=\small}
]
\coordinate (O) at (0,0,0);
\coordinate (P) at (2,-1,-1);
\coordinate (Q) at (2,0.2,-2.2);
\coordinate (V) at
    ({2+0.95/sqrt(3)},
     {0.2+0.95/sqrt(3)},
     {-2.2+0.95/sqrt(3)});

\coordinate (A) at
    ({6.5/sqrt(6)},
     {-3.25/sqrt(6)},
     {-3.25/sqrt(6)});
\coordinate (B) at
    ({3.25/sqrt(6)},
     {3.25/sqrt(6)},
     {-6.5/sqrt(6)});

\path[fill=bridgeplane!9,draw=bridgeplane!55]
    plot[domain=0:360,samples=121,variable=\t]
    ({3.65*2*cos(\t)/sqrt(6)},
     {3.65*(-cos(\t)/sqrt(6)+sin(\t)/sqrt(2))},
     {3.65*(-cos(\t)/sqrt(6)-sin(\t)/sqrt(2))})
    -- cycle;

\path[fill=bridgeplane!32]
    (O) -- (A)
    -- plot[domain=0:60,samples=41,variable=\t]
    ({3.25*2*cos(\t)/sqrt(6)},
     {3.25*(-cos(\t)/sqrt(6)+sin(\t)/sqrt(2))},
     {3.25*(-cos(\t)/sqrt(6)-sin(\t)/sqrt(2))})
    -- cycle;

\draw[bridgeplane,thick]
    plot[domain=0:60,samples=41,variable=\t]
    ({3.25*2*cos(\t)/sqrt(6)},
     {3.25*(-cos(\t)/sqrt(6)+sin(\t)/sqrt(2))},
     {3.25*(-cos(\t)/sqrt(6)-sin(\t)/sqrt(2))});

\draw[gray!45,dashed] (O) -- (-1.5,0,0);
\draw[gray!45,dashed] (O) -- (0,-1.5,0);
\draw[gray!45,dashed] (O) -- (0,0,-2.2);

\draw[->,gray!80,thick]
    (O) -- (4.3,0,0)
    node[left] {$+v_1$};
\draw[->,gray!80,thick]
    (O) -- (0,3.4,0)
    node[right] {$+v_2$};
\draw[->,gray!80,thick]
    (O) -- (0,0,3.2)
    node[above] {$+v_3$};

\draw[->,bridgeplane,thick] (O) -- (A)
    node[left,fill=white,inner sep=2pt]
    {$\mathbb R_+\mathbf e_{3,1}$};

\draw[bridgeplane,thick,dashed] (O) -- (B)
    node[below right]
    {$\mathbb R_+\mathbf e_{3,2}$};

\draw[->,bridgemean,very thick]
    (O) -- ({1.5/sqrt(3)},
            {1.5/sqrt(3)},
            {1.5/sqrt(3)})
    node[right] {$\frac{\sqrt3}{3}\mathbf1$};

\draw[gray!65,densely dashed] (O) -- (Q);
\draw[->,black!75,thick] (O) -- (V);

\draw[->,bridgeprojection,very thick] (O) -- (P);
\draw[->,bridgeresidual,very thick] (P) -- (Q);
\draw[->,bridgemean,very thick] (Q) -- (V);

\draw[bridgeresidual,thin]
    ({2-0.14*2/sqrt(6)},
     {-1+0.14/sqrt(6)},
     {-1+0.14/sqrt(6)})
    --
    ({2-0.14*2/sqrt(6)},
     {-1+0.14/sqrt(6)+0.14/sqrt(2)},
     {-1+0.14/sqrt(6)-0.14/sqrt(2)})
    --
    (2,{-1+0.14/sqrt(2)},{-1-0.14/sqrt(2)});

\fill[black!75] (O) circle (1.2pt);
\fill[bridgeprojection] (P) circle (1.5pt);
\fill[bridgeresidual] (Q) circle (1.5pt);
\fill[black!75] (V) circle (1.5pt);

\node[above left,text=bridgeprojection] at (P)
    {$p\mathbf e_{3,1}$};
\node[below right,text=bridgeresidual] at (Q)
    {$V-\frac h{\sqrt3}\mathbf1$};
\node[above right=2pt] at (V) {$V$};

\draw[bridgeresidual!65,thin]
    ($(P)!0.5!(Q)$)
    -- ([xshift=-4cm,yshift=-4.4cm]O)
    -- ([xshift=-4cm,yshift=-5cm]O)
    node[
        below,
        text=bridgeresidual,
        fill=white,
        inner sep=2pt
    ]
    {$q\bigl(\mathbf e_{3,2}
             -\frac12\mathbf e_{3,1}\bigr)$};

\node[
    right=5pt,
    text=bridgemean,
    fill=white,
    inner sep=2pt
] at ($(Q)!0.55!(V)$)
    {$\frac h{\sqrt3}\mathbf1$};

\node[
    text=bridgeplane,
    fill=white,
    inner sep=2pt
] at ([xshift=-3cm,yshift=3cm]O)
    {$\{v_1+v_2+v_3=0\}$};

\end{tikzpicture}%
}
\caption{The orthogonal decomposition in
\Cref{lem:type-iv-one-sided-bridge}, for $k=1$, $\ell=2$.}
\label{fig:type-iv-one-sided-projection}
\end{figure}
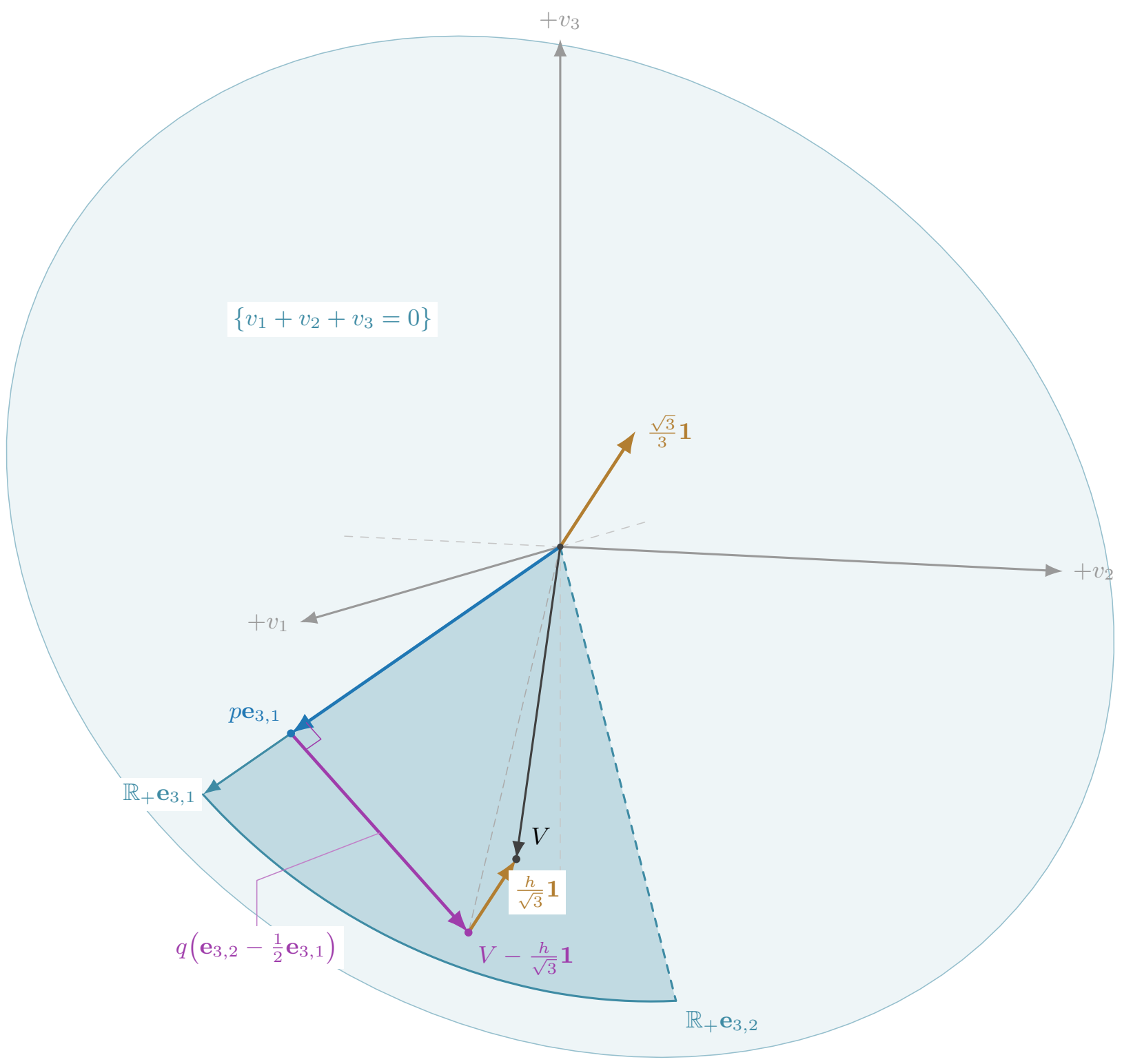
These directions are orthogonal, with
\(
 |\mathbf e_{3,k}|^2=\frac23,\
 \left|\mathbf e_{3,\ell}
       -\frac12\mathbf e_{3,k}\right|^2=\frac12.
\)
Consequently,
\[
 |V'|^2=\frac23|p'|^2+\frac12|q'|^2+|h'|^2.
\]
For the potential, $p=v_k-v_{k+1}+q/2\geq0$ and
subadditivity give
\[
\begin{aligned}
 W_3^\gamma(V,(1,\rho))
 &=w_k\left(p-\frac q2\right)^\gamma+w_\ell q^\gamma\geq w_kp^\gamma
 +(w_\ell-2^{-\gamma}w_k)q^\gamma.
\end{aligned}
\]
Dropping $|h'|^2$, we obtain
\[
\begin{aligned}
 \J_{3,1}^\gamma(V,(0,L))
&\geq\int_0^L\left(\frac23|p'|^2+w_kp^\gamma\right)\,dx+
 \int_0^L\left[
 \frac12|q'|^2+(w_\ell-2^{-\gamma}w_k)q^\gamma
 \right]\,dx.
\end{aligned}
\]
For the first integral, use
\(
 p(0)=h_k+h_\ell/2,\  p(L)=a_kL^\beta
\)
and apply \eqref{eq:type-iv-scalar-bridge}:
\[
\begin{aligned}
\int_0^L\left(\frac23|p'|^2+w_kp^\gamma\right)\,dx
 &\geq
 \J_{3,1}^\gamma(x^\beta\mathbf C^{\{k\}},(0,L))
 -\ep_\gamma(h_k+h_\ell/2)^{1+\gamma/2}\geq
 \J_{3,1}^\gamma(x^\beta\mathbf C^{\{k\}},(0,L))
 -\ep_\gamma(h_1+h_2)^{1+\gamma/2}.
\end{aligned}
\]
For the second integral,
 \eqref{eq:weights_notcanonical} ensures
$w_\ell-2^{-\gamma}w_k>0$ for $\gamma$ 
close to $1$. Using $q(0)=h_\ell$ and $q(L)=0$,
\[
\begin{aligned}
 \int_0^L\left[
 \frac12|q'|^2+(w_\ell-2^{-\gamma}w_k)q^\gamma
 \right]\,dx
 &\geq
 \sqrt{2(w_\ell-2^{-\gamma}w_k)}
 \int_0^L|q'|q^{\gamma/2}\,dx\\
 &\geq
 \sqrt{2(w_\ell-2^{-\gamma}w_k)}
 \left|\int_0^Lq'q^{\gamma/2}\,dx\right|
 =
 \frac{\sqrt{2(w_\ell-2^{-\gamma}w_k)}}{1+\gamma/2}
 h_\ell^{1+\gamma/2}.
\end{aligned}
\]
The coefficient converges to
\[
 \frac23\sqrt{2\rho_*-1}\quad\text{if }k=1,
 \qquad
 \frac23\sqrt{2-\rho_*}\quad\text{if }k=2.
\]
Hence it is bounded below by $c=c(\rho_*)>0$
for $\gamma$ sufficiently close to $1$, depending only
on $\rho_*$ and $\rho(\cdot)$.
Combining the estimates and recalling $\ell=3-k$
proves \eqref{eq:type-iv-one-sided-bridge}.
\end{proof}

\begin{proof}[Proof of \Cref{thm:type-iv-minimality}]
Non-minimality follows from
\Cref{lem:type-iv-overlap-instability}. For minimality, assume \eqref{eq:weights_notcanonical}.
By scaling, it suffices to take $w_1=1$ and
$w_2=\rho$. Fix $R>0$ and an ordered competitor $V$ for
$U_{\mathrm{IV}}$ on $(-R,R)$, so that
\[
  V(R)=R^\beta\mathbf C^{\{1\}},\qquad
  V(-R)=R^\beta\mathbf C^{\{2\}}.
\]
Set $h_i=(v_i-v_{i+1})(0)$ for $i=1,2$.
Apply \Cref{lem:type-iv-one-sided-bridge} on $(0,R)$
to $V(x)$ with $k=1$ and to $V(-x)$ with $k=2$. Adding gives
\begin{align*}
 \J_{3,1}^\gamma(V,(-R,R))
 -\J_{3,1}^\gamma(U_{\mathrm{IV}},(-R,R))&\geq
c\left(h_1^{1+\gamma/2}+h_2^{1+\gamma/2}\right)
 -2\ep_\gamma(h_1+h_2)^{1+\gamma/2}\\
 &\geq
 \left(c\,2^{-\gamma/2}-2\ep_\gamma\right)
 (h_1+h_2)^{1+\gamma/2}.
\end{align*}
Since $c=c(\rho_*)>0$ and $\ep_\gamma\to0$, there exists
$\gamma_{\mathrm{IV}}
=\gamma_{\mathrm{IV}}(\rho_*,\rho(\cdot))\in(2/3,1)$
such that this is nonnegative for all
$\gamma_{\mathrm{IV}}\leq\gamma<1$. Thus $U_{\mathrm{IV}}$ minimizes on every $(-R, R)$.
For any bounded interval, extend the competitor by
$U_{\mathrm{IV}}$ to a larger centered interval.
\end{proof}

\begin{remark}\label{rmk:typeIVrho=1}
For $\rho\equiv1$, we give an explicit choice of
$\gamma_{\mathrm{IV}}$ in \Cref{thm:type-iv-minimality}. By convexity,
\[
 (h_1+h_2/2)^{1+\gamma/2}
 +(h_2+h_1/2)^{1+\gamma/2}
 \leq
 \left(\frac32\right)^{1+\gamma/2}
 \left(h_1^{1+\gamma/2}+h_2^{1+\gamma/2}\right),
\]
with equality when $h_1=h_2$.
Adding the two estimates in
\Cref{lem:type-iv-one-sided-bridge} gives,
\[
\begin{aligned}
 \J_{3,1}^\gamma(V,I)
 -\J_{3,1}^\gamma(U_{\mathrm{IV}},I)\geq
 \left[
 \frac{\sqrt{2(1-2^{-\gamma})}}{1+\gamma/2}
 -\frac{\ep_\gamma}{\sqrt2}
  \left(\frac32\right)^{1+\gamma/2}
 \right]
\left(h_1^{1+\gamma/2}+h_2^{1+\gamma/2}\right).
\end{aligned}
\]
By \eqref{eq:type-iv-scalar-bridge-epsilon}, the coefficient
is nonnegative precisely when
\[
 \sqrt{1-2^{-\gamma}}
 \geq
 4\sqrt3\,\frac{1-\gamma}{2-\gamma}
 \left(\frac32\right)^{\gamma/2}.
\] 
We may take
$\gamma_{\mathrm{IV}}\approx0.911$ to be the unique root of
\(
 (2-\gamma)^2(1-2^{-\gamma})
 =48(1-\gamma)^2\left(\frac32\right)^\gamma,\ \text{for }\gamma\in(2/3,1).
\)
\end{remark}

\appendix
\section{Improvement of flatness for scalar Alt-Phillips with Lipschitz perturbation}
\label{sec:scalar-perturbed-flatness}

Let $0<\gamma<1$ and $\beta=2/(2-\gamma)$.
We consider nonnegative local minimizers $g$ of
\begin{equation}\label{eq:scalar-perturbed-energy}
 \int_\O\left(|\nabla g|^2+g_+^\gamma+G(x,g)\right)dx,
 \qquad G\in\mathbf G_K,
\end{equation}
where $\mathbf G_K$ denotes the scalar version of
\eqref{eq:g-lip}. For $B_r(x_0)\Subset\O$, the rescaling
\begin{equation*}
 g_{x_0,r}(x)=r^{-\beta}g(x_0+rx), \qquad
 G_{x_0,r}(x,t)=r^{2-2\beta}
                 G(x_0+rx,r^\beta t)
                 \in\mathbf G_{Kr^{2-\beta}}
\end{equation*}
preserves minimality in $B_1$. One has flat blowup profile
$\beta^{-\beta}(x\cdot\nu)_+^\beta$. We set
\begin{equation*}
 w=\beta g^{1/\beta},\qquad g=\beta^{-\beta}w^\beta,
\end{equation*}
so the flat profile for $w$ is $(x\cdot\nu)_+$, and
$\{w>0\}=\{g>0\}$. We introduce viscosity solution as in \cite{DSS21}.
\begin{definition}
\label{def:scalar-perturbed-viscosity}
A nonnegative function $w\in C^0(\O)$ is a
\emph{viscosity solution with parameter $K>0$} if
\begin{enumerate}
\item[\textup{(i)}] \emph{Interior inequality.}
$w$ satisfies
\begin{equation}\label{eq:scalar-perturbed-viscosity}
 \left|w\Delta w+(\beta-1)(|\nabla w|^2-1)\right|
 \leq \frac{\beta^{\beta-1}}{2}K\,w^{2-\beta}
 \qquad\text{in }\{w>0\}
\end{equation}
in the viscosity sense, i.e., for any
$B_\rho(x_0)\Subset\{w>0\}$ and
$\varphi\in C^2(B_\rho(x_0))$ with $\varphi(x_0)=w(x_0)$,
\[
\begin{aligned}
 \varphi\geq w \text{ in }B_\rho(x_0)
 &\quad\Longrightarrow\quad
 \varphi(x_0)\Delta\varphi(x_0)+(\beta-1)(|\nabla\varphi(x_0)|^2-1)
 \geq -\frac{\beta^{\beta-1}}{2}K\,\varphi(x_0)^{2-\beta},\\
 \varphi\leq w \text{ in }B_\rho(x_0)
 &\quad\Longrightarrow\quad
 \varphi(x_0)\Delta\varphi(x_0)+(\beta-1)(|\nabla\varphi(x_0)|^2-1)
 \leq \frac{\beta^{\beta-1}}{2}K\,\varphi(x_0)^{2-\beta}.
\end{aligned}
\]
\item[\textup{(ii)}] \emph{Boundary equation.}
$w$ satisfies $|\nabla w|=1$ on $\p\{w>0\}\cap\O$ in the
viscosity sense, i.e, for any
$x_0\in\p\{w>0\}\cap\O$, $B_\rho(x_0)\Subset\O$, and
$\psi\in C^2(B_\rho(x_0))$ with
$\psi(x_0)=0$ and $\nabla\psi(x_0)\neq0$,
\begin{equation}\label{eq:scalar-viscosity-boundary-condition}
\begin{aligned}
 \psi_+\geq w \text{ in }B_\rho(x_0)
 &\quad\Longrightarrow\quad |\nabla\psi(x_0)|\geq1,\\
 \psi_+\leq w \text{ in }B_\rho(x_0)
 &\quad\Longrightarrow\quad |\nabla\psi(x_0)|\leq1.
\end{aligned}
\end{equation}
\end{enumerate}
\end{definition}
If $w \in C^0(B_r(x_0))$ is a viscosity solution with parameter $K$, then
\(
 w_{x_0,r}(x)=r^{-1}w(x_0+rx) \in C^0(B_1)
\) is a viscosity solution with parameter $Kr^{2-\beta}$. We begin by checking $g$ leads to a viscosity solution.
\begin{lemma}
\label{lem:scalar-perturbed-viscosity}
If $g$ locally minimizes \eqref{eq:scalar-perturbed-energy},
then $w=\beta g^{1/\beta}$ is locally Lipschitz and is a
viscosity solution with parameter $K$ in the sense of
\Cref{def:scalar-perturbed-viscosity}.
\end{lemma}
\begin{proof}
\emph{Step 1.} Lipschitz continuity of $w$ follows from $g\in C^{1,\beta-1}$ as in \Cref{thm:optimal-regularity}, and $g \geq 0.$

\emph{Step 2.}
\eqref{eq:vertical-first-variation} applied to \eqref{eq:scalar-perturbed-energy} and using continuity of $g,$
\[
 \gamma g^{\gamma-1}-K
 \le2\Delta g\le\gamma g^{\gamma-1}+K,
 \qquad\{g>0\}\text{ in the viscosity sense}.
\]
Let $\varphi$ touch $w$ from above or below at
$x_0\in\{w>0\}$. On a ball around $x_0$,
$\varphi>0$, and $\beta^{-\beta}\varphi^\beta$ touches $g$
from the same side. Using $\beta(\gamma-1)=\beta-2$ and
$\gamma\beta/2=\beta-1$ gives
\begin{equation}\label{eq:viscosity-solution-chainrule}
    \begin{aligned}
 2\Delta(\beta^{-\beta}\varphi^\beta)
 -\gamma(\beta^{-\beta}\varphi^\beta)^{\gamma-1}
 &=2\beta^{1-\beta}\varphi^{\beta-2}
   \left[\varphi\Delta\varphi
   +(\beta-1)(|\nabla\varphi|^2-1)\right].
\end{aligned}
\end{equation}
Dividing by
$2\beta^{1-\beta}\varphi(x_0)^{\beta-2}>0$ gives \eqref{eq:scalar-perturbed-viscosity}.

\emph{Step 3.} W.l.o.g. let $\psi_+$ touch $w$ at $0\in \p\{w> 0\}$ and $\psi(0)=0$, $\nabla\psi(0)=p e_d$ for $p>0$. By \Cref{thm:gap-blowup-cone}, \(
 w_{0,r_j} = \beta g_{0, r_j}^{1/\beta}\to w_\infty=\beta g_\infty^{1/\beta}\ \text{locally uniformly in }\R^d 
\) to a nontrivial one-homogeneous solution to
\begin{equation}\label{eq:scalar-blowup-interior}
 w_\infty\Delta w_\infty
 +(\beta-1)(|\nabla w_\infty|^2-1)=0
 \qquad\text{in }\{w_\infty>0\}\text{ in the viscosity sense.}
\end{equation}
Moreover,
\(
 r_j^{-1}\psi(r_jx)=p x_d+O(r_j|x|^2),
\)
so the comparison passes to the limit:
\[
 \begin{aligned}
 \psi_+\text{ touches $w$ from below (resp. above)}
 &\quad\Longrightarrow\quad w_\infty\ge (\text{resp.} \leq)\  p(x_d)_+.
 \end{aligned}
\]
We only work with $w_\infty\ge p(x_d)_+$ and prove $p \leq 1$; the other case is similar. For $0<\ep<1$, set
\[
 \varphi_\ep(x)=c_\ep(x_d-\ep|x|), \qquad c_\ep=
 \min_{\substack{\theta\in\S^{d-1}\\\theta_d>\ep}}
 \frac{w_\infty(\theta)}{\theta_d-\ep}
 \ge p.
\]
We compare \( \varphi_\ep\) with $w_\infty$ on $\S^{d-1}.$ For $\theta_d \leq \ep$, $\varphi_\ep \leq0 \leq w_\infty$; for $\theta_d > \ep$,  one has 
\(\varphi_\ep(\theta) = c_\ep(\theta_d-\ep)\leq w_\infty(\theta)\) by definition of $c_\ep$. Therefore $\varphi_\ep \leq w_\infty$ in $\R^d$ by homogeneity. Then we check they touch at a interior point in $\{w_\infty > 0\}$. Indeed, $\frac{w_\infty(\theta)}{\theta_d-\ep}\geq \frac{p \theta_d}{\theta_d - \ep}\to \infty$ as $\theta_d \downarrow \ep$, so the 
minimum \( c_\ep\) is attained at some $\theta_\ep \in \S^{d-1}\cap \{\theta_{\ep,d}>\ep\},$ and thus $ w_\infty(\theta_\ep)= \varphi_\ep(\theta_\ep)= c_\ep(\theta_{\ep, d}-\ep) > 0$.

Now we use our interior inequality at the point $\theta_\ep$. We compute
\[
\varphi_\ep(\theta_\ep)=c_\ep(\theta_{\ep,d}-\ep),\qquad
 \nabla\varphi_\ep(\theta_\ep)=c_\ep(e_d-\ep\theta_\ep),\qquad
 \Delta\varphi_\ep(\theta_\ep)=-c_\ep\ep(d-1).
\]
Using \eqref{eq:scalar-blowup-interior},
\[
\begin{aligned}
 0
&\ge\varphi_\ep(\theta_\ep)\Delta\varphi_\ep(\theta_\ep)
   +(\beta-1)(|\nabla\varphi_\ep(\theta_\ep)|^2-1)\\
 &=-c_\ep^2\ep(d-1)(\theta_{\ep,d}-\ep)
   +(\beta-1)\left[
     c_\ep^2(1-2\ep\theta_{\ep,d}+\ep^2)-1\right]\\
 &=c_\ep^2\left[
   \beta-1-(2\beta+d-3)\ep\theta_{\ep,d}
   +(\beta+d-2)\ep^2\right]-(\beta-1).
\end{aligned}
\]
Choose $\ep >0$ small enough so
$\beta-1-(2\beta+d-3)\ep>0$.
Since $c_\ep\ge p$ and $\theta_{\ep,d}\le1$, it follows that
\[
 p^2\left[\beta-1-(2\beta+d-3)\ep\right]\le\beta-1.
\]
Letting $\ep\downarrow0$ proves $p\le1$. Together with the other direction, this establishes \eqref{eq:scalar-viscosity-boundary-condition}.
\end{proof}

Uniform convergence to \((x\cdot\nu)_+\) does not by itself locate $\p\{w >0\}$.
We need to use $w$ as viscosity solution to trap $\p\{w >0\}$ between
two such flat profiles. This leads to the following;
see \cite[Lemma~4.6]{DSS21}.

\begin{lemma}
\label{lem:scalar-entry-flatness}
There exists $c=c(d,\gamma)\in(0,1/2]$ such that, if $w$ is
a viscosity solution in $B_2$ with $K\le1$ and
\begin{equation}\label{eq:scalar-uniform-closedness}
     \|w-(x\cdot\nu)_+\|_{L^\infty(B_2)}\le c\ep,
 \qquad \nu\in\S^{d-1},\quad 0<\ep<1/4,
\end{equation}
then
\begin{equation}\label{eq:scalar-entry-slab}
 (x\cdot\nu-\ep)_+\le w(x)\le(x\cdot\nu+\ep)_+
 \qquad\text{in }B_1.
\end{equation}
\end{lemma}

\begin{proof}
We first prove that
\begin{equation}\label{eq:scalar-entry-slab-zero-phase-inclusion}
     w=0\qquad\text{in }B_1\cap\{x\cdot\nu\le-\ep/2\}.
\end{equation}
Fix $x$ in this set and take $r=\ep/2$. Then
$\overline B_r(x)\subset B_2\cap\{y\cdot\nu\le0\}$,
so \eqref{eq:scalar-uniform-closedness} gives
\(
 M:=\max_{\overline B_r(x)}w\le c\ep.
\)
Suppose $w(x)>0$. Now $M|y-x|^2/r^2$ is at least $w$ on
$\partial B_r(x)$ and is strictly below $w$ at $x$.
Consequently,
\[
 a:=\max_{\overline B_r(x)}
 \left(w(y)-\frac{M}{r^2}|y-x|^2\right)>0
\]
is attained at some $z\in B_r(x)$. Thus \(
 \varphi(y)=a+\frac{M}{r^2}|y-x|^2
\)
touches $w$ from above at $z\in\{w>0\}$. At $z$,
\[
 0<\varphi(z)\le M,\qquad
 \Delta\varphi(z)=\frac{2dM}{r^2},\qquad
 |\nabla\varphi(z)|^2
 =\frac{4M^2}{r^4}|z-x|^2
 \le\frac{4M^2}{r^2}.
\]
Now \Cref{def:scalar-perturbed-viscosity} gives
\[
\begin{aligned}
 \beta-1
 &\le \varphi(z)\Delta\varphi(z)
    +(\beta-1)|\nabla\varphi(z)|^2
    +\frac{\beta^{\beta-1}}2K\,\varphi(z)^{2-\beta}\\
 &\le \bigl(2d+4(\beta-1)\bigr)\frac{M^2}{r^2}
    +\frac{\beta^{\beta-1}}2M^{2-\beta}\\
 &\le 4\bigl(2d+4(\beta-1)\bigr)c^2
    +\frac{\beta^{\beta-1}}2(c/4)^{2-\beta}.
\end{aligned}
\]We reach contradiction by choosing $c=c(d,\gamma)\in(0,1/2]$ small
so that RHS is strictly less than $\beta-1$.

Since $c\le1/2$, non-negativity and \eqref{eq:scalar-uniform-closedness} now yields
\begin{align*}
 w(x)
 &\ge\bigl((x\cdot\nu)_+-c\ep\bigr)_+
 \ge(x\cdot\nu-\ep)_+
 &&\text{in }B_1,\\
 w(x)
 &\le(x\cdot\nu)_++c\ep
       \le(x\cdot\nu)_++\frac{\ep}{2}
       \le x\cdot\nu+\ep
       =(x\cdot\nu+\ep)_+
 &&\text{in }B_1\cap\{x\cdot\nu>-\ep/2\}.
\end{align*}
Together with \eqref{eq:scalar-entry-slab-zero-phase-inclusion}, these give
\eqref{eq:scalar-entry-slab}.
\end{proof}

The next two lemmas extend
\cite[Theorem~5.1 and Proposition~6.1]{DSS21}
to \Cref{def:scalar-perturbed-viscosity}.

\begin{lemma}[Harnack inequality]
\label{lem:scalar-flat-harnack}
There exist $\ep_H,c_H\in(0,1)$ and $\eta_H\in(0,1/2)$,
depending only on $d,\gamma$ s.t. if $w$ is a viscosity solution in $B_1$ with parameter $K$,
$0\in\partial\{w>0\}$, and
\begin{equation}
\label{eq:scalar-flat-harnack-slab-bound}
 (x_d+a)_+\le w\le(x_d+b)_+
 \quad\text{in }B_1,\qquad a\le0\le b,
\end{equation}
where
\(
 0<b-a\le\ep_H,\ K\le\ep_H(b-a),
\)
then there exist $a',b'$ such that
\begin{equation}\label{eq:scalar-flat-harnack}
\begin{gathered}
 a\le a'\le b'\le b,\qquad
 b'-a'\le(1-c_H)(b-a),\\
 (x_d+a')_+\le w\le(x_d+b')_+
 \qquad\text{in }B_{\eta_H}.
\end{gathered}
\end{equation}
\end{lemma}
\begin{proof}
\emph{Step 1.} Put $\ep=b-a$. For $\ell = \ell(d) > 0$ small to choose, set
\[
 D=\{|x'|<3/4,\ \ell/2<x_d<1/2\},\qquad
 T=\{|x'|\le1/2,\ x_d=\ell\},\qquad
 \bar x=\frac15e_d.
\]
Choosing
$\ep_H\le\ell/4$ gives
\(
 \frac{\ell}{4}\le x_d+a\le w\le x_d+b\le1\ \text{in }D.
\) Also, by chain rule in \eqref{eq:viscosity-solution-chainrule}, $w^\beta$ solves
\begin{equation}\label{eq:viscosity-solution-chainrule-equation}
    \left|
 \Delta(w^\beta)-\beta(\beta-1)w^{\beta-2}
 \right|
 \le\frac{\beta^\beta}{2}K
 \qquad\text{in }D \ \text{ in the viscosity sense.}
\end{equation}

Define the two nonnegative differences
\[
 v_a=w^\beta-(x_d+a)^\beta,\qquad
 v_b=(x_d+b)^\beta-w^\beta.
\]
Apply the
mean value theorem to $t\mapsto t^{(\beta-2)/\beta}$ so that, for  $\xi_a \in [(x_d+a)^\beta, w^\beta] \subseteq [(\ell/4)^\beta, 1] $ 
\[
\begin{aligned}
 0
 &\ge\beta(\beta-1)
       \bigl(w^{\beta-2}-(x_d+a)^{\beta-2}\bigr)=-(\beta-1)(2-\beta)\xi_a^{-2/\beta}v_a\ge-(\beta-1)(2-\beta)\left(\frac4\ell\right)^2v_a.
\end{aligned}
\]
Now with the above and \eqref{eq:viscosity-solution-chainrule-equation} applied to 
\[
\begin{aligned}
 \Delta v_a
 &=\Delta(w^\beta)
   -\beta(\beta-1)(x_d+a)^{\beta-2}
 =\left[\Delta(w^\beta)
         -\beta(\beta-1)w^{\beta-2}\right]
   +\beta(\beta-1)
      \left[w^{\beta-2}-(x_d+a)^{\beta-2}\right],
\end{aligned}
\]and repeating the same for $v_b$, one obtain the following in the viscosity sense in $D$
\[
 -(\beta-1)(2-\beta)\left(\frac4\ell\right)^2v_s
 -\frac{\beta^\beta}{2}K
 \le\Delta v_s
 \le\frac{\beta^\beta}{2}K,
 \qquad s\in\{a,b\}.
\]
The interior Harnack inequality in $D$ therefore gives a constant
$C_\ell=C_\ell(d,\gamma,\ell)\ge1$ such that
\[
 \inf_Tv_s\ge C_\ell^{-1}v_s(\bar x)-C_\ell K,
 \qquad s\in\{a,b\}.
\]

At $\bar x$, at least one of $v_a(\bar x)$, $v_b(\bar x)$ has size comparable
to $\ep$. Indeed, taking $\ep_H\le1/10$, we have
\[
\begin{aligned}
 v_a(\bar x)+v_b(\bar x)
 &=\left(\frac15+b\right)^\beta
   -\left(\frac15+a\right)^\beta=\beta\int_a^b
       \left(\frac15+t\right)^{\beta-1}\,dt
 \ge\beta\,10^{1-\beta}\ep.
\end{aligned}
\]

Consequently, in view of $K \leq \ep_H \ep$, Harnack gives
\[
\begin{aligned}
 \max\{\inf_Tv_a,\inf_Tv_b\}
 &\ge C_\ell^{-1}\max\{v_a(\bar x),v_b(\bar x)\}-C_\ell K\ge\frac{\beta\,10^{1-\beta}}{2C_\ell}\ep-C_\ell K\ge\left(
       \frac{\beta\,10^{1-\beta}}{2C_\ell}-C_\ell \ep_H
      \right)\ep \ge\beta c_0\ep,
\end{aligned}
\]
where we decrease
\(
 \ep_H\le\frac{\beta\,10^{1-\beta}}{4C_\ell^2},
\)
and then fix
\(
 c_0=c_0(d,\gamma, \ell)=\frac{10^{1-\beta}}{4C_\ell}>0.
\)

If $\inf_Tv_a\ge\beta c_0\ep$, then, for every $x\in T$,
using $\ell/4\le x_d+a\le w(x)\le1$, we obtain
\[
\begin{aligned}
 \beta c_0\ep
 &\le v_a(x)
 =w(x)^\beta-(x_d+a)^\beta
 =\beta\int_{x_d+a}^{w(x)}t^{\beta-1}\,dt
 \le\beta\bigl(w(x)-(x_d+a)\bigr).
\end{aligned}
\]
Repeating the same if $\inf_Tv_b\ge\beta c_0\ep$, one has either
\begin{equation}\label{eq:scalar-harnack-top-gap}
 w\ge x_d+a+c_0\ep\quad\text{on }T,
 \qquad\text{or}\qquad
 w\le x_d+b-c_0\ep\quad\text{on }T.
\end{equation}

\emph{Step 2.} We now construct comparison functions that carry  \eqref{eq:scalar-harnack-top-gap} from $T$ toward the origin. Set
\[
 \mathcal C=\{|x'|<1/2,\ |x_d|<\ell\},\qquad
 \zeta(x)=\frac18-|x'|^2+dx_d^2+x_d.
\]
Choosing $\ell = \min\{1/32, 1/(16d)\}$ gives
\[
 \Delta\zeta=2,\quad
 \partial_d\zeta\ge\frac12,\quad
 |\nabla\zeta|\le2,\quad |\zeta|\le\frac14
 \quad\text{in }\overline{\mathcal C}, \qquad \text{and} \qquad
 \zeta\le-\frac1{16}
 \quad\text{on }\{|x'|=1/2,\ |x_d|\le\ell\}.
\]

Choose $0<\tau\le\min\{c_0,1/4\}$ and define
\[
 \psi_-=x_d+a+\tau\ep\zeta,\qquad
 \psi_+=x_d+b-\tau\ep\zeta.
\]
For sufficiently small $\ep_H$,
\begin{align*}
 \Delta\psi_-&=2\tau\ep,
 &
 |\nabla\psi_-|^2-1
 &=2\tau\ep\,\partial_d\zeta
   +\tau^2\ep^2|\nabla\zeta|^2
 \ge\tau\ep,\\
 \Delta\psi_+&=-2\tau\ep,
 &
 |\nabla\psi_+|^2-1
 &=-2\tau\ep\,\partial_d\zeta
   +\tau^2\ep^2|\nabla\zeta|^2
 \le-\frac{\tau\ep}{2}.
\end{align*}
In particular,
$0<|\nabla\psi_+|<1<|\nabla\psi_-|$ in $\overline{\mathcal C}$.
One may decrease $\ep_H$ further so that
\(
 \frac{\beta^{\beta-1}}2\ep_H
 <\frac{(\beta-1)\tau}{4},
\) and use $0\le w\le1$ in $\mathcal C$, $K \leq \ep_H \ep$ to obtain
\begin{align}
 w\Delta\psi_-+(\beta-1)(|\nabla\psi_-|^2-1)
 &\ge(\beta-1)\tau\ep >\frac{\beta^{\beta-1}}2K w^{2-\beta} , \label{eq:scalar-harnack-interior-compare-1}\\
 w\Delta\psi_++(\beta-1)(|\nabla\psi_+|^2-1)
 &\le-\frac{(\beta-1)\tau\ep}{2} < -\frac{\beta^{\beta-1}}2K w^{2-\beta}.\label{eq:scalar-harnack-interior-compare-2}
\end{align}

\emph{Step 3.} 
Suppose the first inequality in
\eqref{eq:scalar-harnack-top-gap} holds. Then
\begin{align*}
 \psi_-(x)
 &\le x_d+a+\frac{\tau\ep}{4}<w(x),
 && |x'|\le\frac12,\quad x_d=\ell,\\
 \psi_-(x)
 &\le x_d+a-\frac{\tau\ep}{16}<x_d+a\le w(x),
 && |x'|=\frac12,\quad |x_d|\le\ell,\\
 \psi_-(x)
 &\le-\ell+\frac{\tau\ep}{4}<0\le w(x),
 && |x'|\le\frac12,\quad x_d=-\ell.
\end{align*}
Thus $\psi_-<w$ on $\partial\mathcal C$. If
\(
 m=\max_{\overline{\mathcal C}}(\psi_--w)>0,
\)
the maximum occurs at $z\in\mathcal C$, and
$\psi_--m$ touches $w$ from below there.
If $w(z)>0$, the strict interior inequality \eqref{eq:scalar-harnack-interior-compare-1} contradicts
\eqref{eq:scalar-perturbed-viscosity}.
If $w(z)=0$, then $(\psi_--m)_+$ touches $w$ from below.
Its nonzero gradient implies $z\in\partial\{w>0\}$, and
$|\nabla\psi_-(z)|>1$ contradicts
\eqref{eq:scalar-viscosity-boundary-condition}.
Thus
\[
 w\ge(\psi_-)_+\quad\text{in }\mathcal C.
\]

Suppose the second inequality in
\eqref{eq:scalar-harnack-top-gap} holds, put
\(
 S=\overline{\mathcal C\cap\{w>0\}}.
\)
\eqref{eq:scalar-flat-harnack-slab-bound} gives $x_d+b\ge0$ on $S$, so 
\begin{align*}
 w(x)
 &\le x_d+b-c_0\ep
 <x_d+b-\frac{\tau\ep}{4}\le\psi_+(x),
 &&\text{on }S\cap\{x_d=\ell\},\\
 w(x)
 &\le x_d+b
 <x_d+b+\frac{\tau\ep}{16}\le\psi_+(x),
 &&\text{on }S\cap\{|x'|=1/2\}.
\end{align*}
Since $S$ does not meet $\{|x^\prime|\leq 1/2, \,x_d= -\ell\}$, this proves
$w<\psi_+$ on $S\cap\partial\mathcal C$.
If
\(
 m=\max_S(w-\psi_+)>0,
\)
the maximum occurs at $z\in S\cap\mathcal C$, and
\[
 w\le(\psi_++m)_+\quad\text{in }\mathcal C,
 \qquad w(z)=\psi_+(z)+m.
\]If $w(z)>0$, $\psi_++m$ touches $w$ from above, and \eqref{eq:scalar-harnack-interior-compare-2} contradicts
\eqref{eq:scalar-perturbed-viscosity}.
If $w(z)=0$, $(\psi_++m)_+$ touches $w$ from above. Since $z\in S\cap\mathcal C$, we have
$z\in\partial\{w>0\}$, and
$0<|\nabla\psi_+(z)|<1$ contradicts
\eqref{eq:scalar-viscosity-boundary-condition}.
Hence
\[
 w\le(\psi_+)_+\quad\text{in }\mathcal C.
\]

\emph{Step 4.}  Finally, choose $\eta_H\in(0,\ell/2)$ so that
$\zeta\ge1/16$ in $B_{\eta_H}$, and put $c_H=\tau/16$.
The two alternatives give, respectively,
\[
 w\ge(x_d+a+c_H\ep)_+,
 \qquad\text{or}\qquad
 w\le(x_d+b-c_H\ep)_+
 \quad\text{in }B_{\eta_H}.
\]
Taking
\(
 (a',b')=(a+c_H\ep,b)
 \ \text{or}\
 (a',b')=(a,b-c_H\ep)
\)
proves \eqref{eq:scalar-flat-harnack}.
\end{proof}

Under the assumptions of \Cref{lem:scalar-flat-harnack}, set \( \ep=b-a\), \(\Omega=B_1\cap\overline{\{w>0\}},\)  $\theta=1-c_H$, and decrease $\eta_H$ so $\eta_H<\theta$.
As in \cite[Corollary~5.2]{DSS21},
repeated application of \Cref{lem:scalar-flat-harnack} gives
nested intervals $[a_m,b_m]$ s.t.
\[
 (x_d+a_m)_+\le w\le(x_d+b_m)_+
 \quad\text{in }B_{\eta_H^m},\qquad
 b_m-a_m=\ep\theta^m,
\] whenever
\begin{equation}\label{eq:scalar-flat-holder-modulus-assumption}
    \ep\theta^m\eta_H^{-m}\le\ep_H, \qquad 
 \frac{K\eta_H^{m(2-\beta)}}{\ep\theta^m\eta_H^{-m}}
 =
 \frac K\ep
 \left(\frac{\eta_H^{3-\beta}}{\theta}\right)^m
 \le\frac K\ep\le\ep_H.
\end{equation}
Note the latter is always satisfied since we start with $K \leq \ep_H \ep$. If so, the above implies 
\[\frac{a_m}{\ep} \leq v(x) \defeq  \frac{w(x) - x_d}{\ep} \leq \frac{b_m}{\ep},\qquad x \in \O \cap B_{\eta_H^m}, \qquad \text{equivalently} \qquad \underset{\O \cap B_{\eta_H^m}}{\mathrm{osc}} v \leq \theta^m.  \]

For $x\in\Omega$ with $|x|\ge\ep/\ep_H$, choose $m\ge0$
such that $\eta_H^{m+1}\le|x|<\eta_H^m$.
All preceding iterations are now admissible, since
$\ep\theta^k\eta_H^{-k}\le\ep/|x|\le\ep_H$ for $k\le m$ so \eqref{eq:scalar-flat-holder-modulus-assumption} holds. Hence, for $\kappa=\log\theta/\log\eta_H\in(0,1)$,
\begin{equation*}
 |v(x)-v(0)|
 \le\theta^m
 =(\eta_H^m)^\kappa
 \le\eta_H^{-\kappa}|x|^\kappa,
 \qquad x\in\Omega,\quad |x|\ge\frac{\ep}{\ep_H}.
\end{equation*}

Applying the iteration after translating and rescaling at free-boundary points, and combining it with
interior elliptic estimates in $\{w>0\}$, gives the two-point
estimate; see \cite[(4.5)]{DS11}.
\begin{equation}\label{eq:scalar-flat-holder-modulus}
 |v(x)-v(y)|\le C|x-y|^\kappa,
 \qquad x,y\in\Omega\cap B_{1/2},\quad
 |x-y|\ge C\ep/\ep_H,
\end{equation}
where $C$ depends only on $d,\gamma$. Now we follow \cite[Proposition~6.1]{DSS21}.
\begin{lemma}[Improvement of flatness]
\label{lem:scalar-improvement-flatness}
There exist $\ep_F\in(0,1)$, $\rho\in(0,1/2)$,
$\alpha\in(0,2-\beta)$, and $C>0$,
depending on $d,\gamma$ s.t. if $w$ is a viscosity solution in $B_1$ with parameter $K$,
$0\in\partial\{w>0\}$, and for $\nu\in\S^{d-1}$,
\begin{equation}\label{eq:scalar-initial-flatness}
 (x\cdot\nu-\ep)_+\le w\le(x\cdot\nu+\ep)_+
 \quad\text{in }B_1,\qquad
 0<\ep\le\ep_F,\quad K\le\ep_F\ep,
\end{equation}
then there exists some $\nu'\in\S^{d-1}$ s.t.
$|\nu'-\nu|\le C\ep$ and
\begin{equation}\label{eq:scalar-improved-flatness}
 (x\cdot\nu'-\ep\rho^{1+\alpha})_+
 \le w\le
 (x\cdot\nu'+\ep\rho^{1+\alpha})_+
 \qquad\text{in }B_\rho.
\end{equation}
\end{lemma}
\begin{proof}
Assume not, then we fix $\rho, \alpha, C$ to be specified w.r.t. $d$, $\gamma$ later, and choose $w_j$ with
\begin{equation}\label{eq:scalar-initial-flatness-assumption}
    \ep_j\downarrow0,\qquad
 \frac{K_j}{\ep_j}\longrightarrow0,\qquad
 (x_d-\ep_j)_+\le w_j\le(x_d+\ep_j)_+,
 \qquad 0\in\partial\{w_j>0\},
\end{equation}
where we rotate each $\nu_j = e_d$ for simplicity, that does not satisfy \eqref{eq:scalar-improved-flatness} for any choice of $\nu'.$ Set
\[
 \Omega_j=B_1\cap\overline{\{w_j>0\}},\qquad
 v_j=\frac{w_j-x_d}{\ep_j}\quad\text{on }\Omega_j, \qquad \text{so} \qquad |v_j|\leq 1, \qquad v_j(0) =0.
\]

\emph{Step 1. Compactness.} For our sequence $\ep_j$, $K_j$, \eqref{eq:scalar-flat-holder-modulus-assumption} are satisfied for $j$ sufficiently large, and since $\ep_j \downarrow0$, on compact subsets of $B_1\cap\{x_d>0\}$, \eqref{eq:scalar-flat-holder-modulus} gives Ascoli-Arzelà. As in
\cite[proof of Proposition~6.1, Step~1]{DSS21}, after passing to a subsequence, graphs of $v_j$  locally
Hausdorff converge to the graph of
\[
 v\in C^{0, \kappa}(B_{1/2}\cap\{x_d\ge0\}),\qquad
 |v|\le1,\qquad v(0)=0, \qquad \kappa = \kappa(d, \gamma)\in (0, 1).
\]

\emph{Step 2. Limiting equation.}
Set $s=2(\beta-1)\in(0,2)$.
We claim that $v$ solves the limiting problem
\begin{equation*}
 \begin{cases}
 \displaystyle
 \Delta v+s\frac{\partial_dv}{x_d}=0
 &\text{in }B_{1/2}\cap\{x_d>0\},\\[1ex]
 v_s:=\lim_{t\downarrow0}
 \frac{v(x',t)-v(x',0)}{t^{1-s}}=0
 &\text{on }B_{1/2}\cap\{x_d=0\},
 \end{cases}
\end{equation*}
with boundary condition understood in the viscosity
sense of \cite[Definition~7.1]{DSS21}.

\emph{Interior equation.} We follow \cite[Proposition~2.9]{CC95}.
In $\{w_j > 0\}$, using \eqref{eq:scalar-perturbed-viscosity},  the inequality for $v_j$ is
\[
 \left|
 (x_d+\ep_jv_j)\Delta v_j
 +2(\beta-1)\partial_dv_j
 +(\beta-1)\ep_j|\nabla v_j|^2
 \right|
 \le
 \frac{\beta^{\beta-1}}2
 \frac{K_j}{\ep_j}(x_d+\ep_jv_j)^{2-\beta}.
\]
Take ball
$\overline B\subset B_{1/2}\cap\{x_d>0\}$ around $x_0$ s.t. $v\le\phi$ on $\overline B$ and
$v(x_0)=\phi(x_0)$. Fix $\delta>0$. Take $x_j\in\overline B$ maximizing
\(
 v_j(x)-\phi(x)-\delta|x-x_0|^2\ \text{over }\overline B.
\)
The choice of $x_j$ and local uniform convergence then gives
\[
 \delta|x_j-x_0|^2
 \le (v_j-\phi)(x_j)-(v_j-\phi)(x_0)
 \le 2\|v_j-v\|_{L^\infty(\overline B)}
 \longrightarrow0,
\]
so $x_j\to x_0$ and
$x_j\in B$ for large $j$. Thus the function
\[
 x\longmapsto
 \phi(x)+\delta|x-x_0|^2
 +v_j(x_j)-\phi(x_j)-\delta|x_j-x_0|^2
\]
touches $v_j$ from above at $x_j$ in $B$.
Moreover, since $x_{0,d}>0$, flatness from \eqref{eq:scalar-initial-flatness-assumption}
ensure that $x_j\in\{w_j>0\}$ for large $j$.
Applying the viscosity inequality now gives
\[
\begin{aligned}
 (x_{j,d}+\ep_jv_j(x_j))
\bigl(\Delta\phi(x_j)+2d\delta\bigr)+s\bigl(\partial_d\phi(x_j)
            &+2\delta(x_{j,d}-x_{0,d})\bigr)
 +(\beta-1)\ep_j
       |\nabla\phi(x_j)+2\delta(x_j-x_0)|^2\\
       &\ge
 -\frac{\beta^{\beta-1}}2
 \frac{K_j}{\ep_j}
 (x_{j,d}+\ep_jv_j(x_j))^{2-\beta}.
\end{aligned}
\]
Using $|v_j|\leq 1$,
$\ep_j\to0$, and $K_j/\ep_j\to0$, passing $j \to \infty$ gives
\(
 x_{0,d}\bigl(\Delta\phi(x_0)+2d\delta\bigr)
 +s\,\partial_d\phi(x_0)\ge0.
\)
Letting $\delta\downarrow0$  then yields
\(
 \Delta\phi(x_0)
 +s\frac{\partial_d\phi(x_0)}{x_{0,d}}\ge0.
\)
Touching from below is similar.

\emph{Boundary condition.}
For $s\ge1$, the boundary condition in
\cite[Definition~7.1]{DSS21} follows from $|v|\le1$. Now suppose $0<s<1$ and that
\[
 A|x'-y'|^2+B+p x_d^{1-s}, \qquad A, B \in \R, \ y^\prime \in \R^{d-1}
\]
touches $v$ at $x_0\in B_{1/2}\cap\{x_d=0\}$,
from above with $p<0$, or from below with $p>0$.
The function
\[
 \phi(x)=A|x'-y'|^2+B
 -\operatorname{sgn}(p)|x'-x_0'|^2+\frac p2x_d
\]
touches strictly from the same side in a sufficiently small
neighborhood of $x_0$ in $\{x_d\ge0\}$.
Graph convergence gives $z_j\in\Omega_j=B_1 \cap \overline{\{w_j > 0\}}$, $z_j\to x_0$,
and constants $c_j\to0$, such that
\[
 \Phi_j=x_d+\ep_j(\phi+c_j)
\]
touches $w_j$ from the corresponding side at $z_j$ in $\Omega_j$. If $z_j \in \{w_j > 0\}$ one is allowed to use viscosity inequality \eqref{eq:scalar-perturbed-viscosity} directly. If $z_j \in \p\{w_j > 0\}$, one needs to ensure touching of $(\Phi_j)_+$ happens in a full neighborhood around $z_j$. In this case, $(\Phi_j)_+\ge w_j$ also happens in the zero phase, so it remains to work with $(\Phi_j)_+ \leq w_j$ at $z_j \in \p\{w_j > 0\}$ in $\Omega_j$. Take cylinder
$Q=B_r'(x_0')\times(-r,r).$
For large $j$, $z_j\in Q$, and \eqref{eq:scalar-initial-flatness-assumption} ensure
that $B_r'(x_0')\times\{r\}$ lies in $\{w_j>0\}$. Now assume $\Phi_j(y)>0=w_j(y)$ at some $y\in Q$, so moving in positive $e_d$ direction from $y$ reaches a point $q\in Q \cap \p\{w_j > 0\}$. But 
\(
 \partial_d\Phi_j=1+\ep_jp/2>0
\)
implies
\(
 \Phi_j(q)\ge\Phi_j(y)>0=w_j(q),
\)
contradicting $\Phi_j \leq w_j$ on $Q\cap\Omega_j$.
Thus $(\Phi_j)_+$ touches $w_j$ from below at $z_j$
in a full neighborhood $Q$ around $z_j$.

Now one may apply the definition of $w_j$ as \Cref{def:scalar-perturbed-viscosity}.
At $z_j \in \{w_j > 0\}$, using $w_j(z_j)\to0$,
\begin{align*}
 \frac{
 w_j(z_j)\Delta\Phi_j(z_j)
 +(\beta-1)(|\nabla\Phi_j(z_j)|^2-1)}
 {\ep_j}&=
 w_j(z_j)\Delta\phi(z_j)
 +(\beta-1)\bigl(p+\ep_j|\nabla\phi(z_j)|^2\bigr)
 \longrightarrow(\beta-1)p,\\
 \frac{\beta^{\beta-1}}2
 \frac{K_j}{\ep_j}w_j(z_j)^{2-\beta}&\longrightarrow0.
\end{align*}
The sign of $p$ contradicts \eqref{eq:scalar-perturbed-viscosity} touching from either side.
At $z_j\in\partial\{w_j>0\}$,
\[
 |\nabla\Phi_j(z_j)|^2-1
 =\ep_j\bigl(p+\ep_j|\nabla\phi(z_j)|^2\bigr),
\]
also contradicts the sign in
\eqref{eq:scalar-viscosity-boundary-condition} for $\ep_j$ sufficiently small.
Thus $v$ satisfies 
\cite[Definition~7.1]{DSS21}.

\emph{Step 3. Improvement of Flatness.}
By \cite[Theorem~7.2]{DSS21}, there exist
$\alpha_0 = \alpha_0(d)>0$ and $a\in\R^{d-1}$ s.t.
\[
 |a|\le C_0 = C_0(d),\qquad
 |v(x)-a\cdot x'|\le C_0|x|^{1+\alpha_0}
 \quad\text{in }B_{1/4}\cap\{x_d\ge0\}.
\]
Choose
\(
 0<\alpha<\min\{\alpha_0,2-\beta\},\) then $\rho$ small so \(
 C_0(2\rho)^{1+\alpha_0}
 \le\frac18\rho^{1+\alpha}.
\)
Graph Hausdorff convergence gives
\[
 |v_j(x)-a\cdot x'|
 \le\frac14\rho^{1+\alpha}
 \quad\text{on }\Omega_j\cap B_{2\rho},
\] for large $j$.
Set
\(
 \nu_j=\frac{(\ep_j a,1)}
 {\sqrt{1+\ep_j^2|a|^2}}.
\)
Then $|\nu_j-e_d|\le C\ep_j$, with $C = C_0(d)$, and
\[
 |w_j(x)-x\cdot\nu_j|
 \le\frac12\ep_j\rho^{1+\alpha}
 \quad\text{on }\Omega_j\cap B_{2\rho}
\]
for large $j$. To extend the lower bound to $x\in B_\rho \cap \{w_j = 0\}$,
move upward in $e_d$ so $y = (x^\prime, 2\ep_j)$ is positive by \eqref{eq:scalar-initial-flatness-assumption}. The segment from $x$ to $y$ then meets
$z = (x^\prime, z_d)\in B_{2\rho} \cap \p\{w_j >0\}$ with $z_d\ge x_d$. Thus
\(
 x\cdot\nu_j\le z\cdot\nu_j
 \le\frac12\ep_j\rho^{1+\alpha}.
\)
Consequently,
\[
 (x\cdot\nu_j-\ep_j\rho^{1+\alpha})_+
 \le w_j\le
 (x\cdot\nu_j+\ep_j\rho^{1+\alpha})_+
 \quad\text{in }B_\rho,
\]
contradicting the choice of $w_j$.
\end{proof}

Together we establish the following theorem; see \cite[Theorem~1.5]{DSS21}.

\begin{theorem}[Regularity at points with a flat blowup]
\label{thm:scalar-flat-blowup-regularity}
Let $g\ge0$ locally minimize \eqref{eq:scalar-perturbed-energy},
with $G\in\mathbf G_K$. Let
$x_0\in\partial\{g>0\}\cap\O$.
Suppose that, for some $r_j\downarrow0$ and $\nu\in\S^{d-1}$,
\begin{equation}\label{eq:scalar-flat-blowup}
 r_j^{-\beta}g(x_0+r_jx)
 \longrightarrow
 \beta^{-\beta}(x\cdot\nu)_+^\beta
 \qquad\text{locally uniformly in }\R^d.
\end{equation}
Then, for some $r_0>0$,
$\partial\{g>0\}\cap B_{r_0}(x_0)$ is a $C^{1,\alpha}$
graph in the $\nu$ direction, where
$\alpha=\alpha(d,\gamma)\in(0,2-\beta)$.
Its unit normal pointing into $\{g>0\}$ at $x_0$ is $\nu$.
\end{theorem}

\begin{proof}
Let $\ep_F,\rho,\alpha,C$ be the constants in
\Cref{lem:scalar-improvement-flatness}. Set
\[
 w=\beta g^{1/\beta},\qquad
 w_{j}(x)=\frac{w(x_0+r_jx)}{r_j},\qquad
 K_j=Kr_j^{2-\beta}.
\]
By \Cref{lem:scalar-perturbed-viscosity}, $w_j$ is a
viscosity solution with parameter $K_j$. Moreover,
\eqref{eq:scalar-flat-blowup} gives
$w_j\to(x\cdot\nu)_+$ uniformly in $B_2$, and $K_j\to0$.
Fix $0<\ep\le\ep_F/4$. For sufficiently large $j$, and $c$ as the constant in \Cref{lem:scalar-entry-flatness},
\[
 \|w_j-(x\cdot\nu)_+\|_{L^\infty(B_2)}\le c\ep,
 \qquad K_j\le\ep_F\ep.
\]
By \Cref{lem:scalar-entry-flatness}, $w_j$ satisfies
\eqref{eq:scalar-entry-slab}. Together with the above
smallness conditions, this verifies
\eqref{eq:scalar-initial-flatness}.

After applying \Cref{lem:scalar-improvement-flatness}
and rescaling by $w_{j, \rho}(x) = w_j(\rho x)/\rho$, 
\[(x\cdot \nu' - \ep_\rho )_+ \leq w_{j, \rho}\leq (x \cdot \nu' + \ep_\rho)_+ \quad \text{in $B_1$}, \qquad 0 < \ep_\rho= \rho^\alpha\ep\leq \ep_F, \quad \frac{K_{j,\rho}}{\ep_\rho}
 =\rho^{2-\beta-\alpha}\frac{K_j}{\ep}\leq \ep_F.\]
Since $\alpha<2-\beta$, the smallness conditions are
preserved. Moreover, for $p\in\partial\{w_j>0\}\cap B_{1/2}$, using
$|p\cdot\nu|\le\ep$, the re-centering and rescaling $2w_j(p+x/2)$ satisfies the
same hypotheses \eqref{eq:scalar-initial-flatness} for \Cref{lem:scalar-improvement-flatness} with flatness $4\ep$ and parameter
$2^{\beta-2}K_j$. Thus we can iterate uniformly in \(\partial\{w_j>0\}\cap B_{1/2}.\) 

In the following we carry out details to show how improvement of flatness gives $C^{1, \alpha}$ as in \cite{AC81}, \cite{CS05}.

\emph{Step 1.}
Fix $j$ and $p \in \partial\{w_j>0\}\cap B_{1/2}.$ Iterate starting with $\nu_{p, 0} = e_d$ gives a sequence $\nu_{p, m} \in \S^{d-1}$ s.t.
\begin{equation}\label{eq:scalar-flat-blowup-regularity-mid1}
    ((x - p)\cdot \nu_{p, m} - 4\ep \rho^{m\alpha } r_m)_+ \leq w_j \leq ((x - p)\cdot \nu_{p, m} + 4\ep \rho^{m\alpha } r_m)_+, \qquad x \in B_{r_m}(p), \quad r_m = \rho^m/2.
\end{equation}
and
\begin{equation}\label{eq:scalar-flat-blowup-regularity-mid2}
    |\nu_{p, m+1} - \nu_{p, m}|\leq4 C \ep \rho^{m\alpha}.
\end{equation}
Using \eqref{eq:scalar-flat-blowup-regularity-mid2}, there exists $\nu_{p, \infty} \in \S^{d-1}$ s.t.
\begin{equation}\label{eq:scalar-flat-blowup-regularity-mid3}
    |\nu_{p, \infty} - \nu_{p, m}|\leq \frac{4C \ep}{1 - \rho^\alpha}\rho^{m\alpha}, \qquad |\nu_{p, \infty} - e_d|\leq \frac{4C \ep}{1 - \rho^\alpha}.
\end{equation}
Now for any $q \in \p\{w_j > 0\} \cap B_{r_m}(p)$,  \eqref{eq:scalar-flat-blowup-regularity-mid1} gives $|(q - p)\cdot \nu_{p, m} |\leq 4\ep \rho^{m\alpha } r_m.$ For any $p,\, q \in \p\{w_j > 0\} \cap B_{1/2}$ with $0 < |p - q|< 1/2$, one may choose $m$ so $r_{m+1} \leq |p - q|< r_m$. Using \eqref{eq:scalar-flat-blowup-regularity-mid3}
\begin{align}
    |(q -p)\cdot \nu_{p, \infty}| &\leq  |(q -p)\cdot \nu_{p, m}| + |(q - p)\cdot (\nu_{p, \infty} - \nu_{p, m})| \notag \\
    &\leq 4\ep \rho^{m\alpha } r_m + \frac{4C \ep}{1 - \rho^\alpha}\rho^{m\alpha}|q - p|\leq M \ep |q- p|^{1 + \alpha}, \qquad M = M(d, \gamma) > 0. \label{eq:scalar-flat-blowup-regularity-mid4}
\end{align}

\emph{Step 2.} We prove the free boundary is a single Lipschitz graph. Choose $\ep < 1/16$ small s.t. \(\frac{4C \ep}{1 - \rho^\alpha}\leq \frac{1}{4}\),  \(M \ep \leq \frac{1}{4},\) so using \eqref{eq:scalar-flat-blowup-regularity-mid3}, $(\nu_{p, \infty})_d \geq 3/4$ and $|\nu_{p, \infty}'| \leq 1/4$. We consider cylinder $Q = B_{1/8}' \times (-1/8, 1/8).$ \eqref{eq:scalar-flat-blowup-regularity-mid1} at $p= 0$, $m = 0$ implies 
\[w_j(x', -1/8) =0, \qquad w_j(x', 1/8)\geq 1/8 - 2\ep > 0,\] thus each vertical line through $Q$ meets $\p\{w_j > 0\}$ at least once. For any $p,q \in \p\{w_j > 0\} \cap Q$, using \eqref{eq:scalar-flat-blowup-regularity-mid4}
\begin{align*}
    \frac{3}{4}|q_d - p_d|&\leq |\nu_{p, \infty}'||q' - p'| + |(q  - p) \cdot \nu_{p, \infty}|\leq \frac{1}{4}|q' - p'| + M \ep |q- p|^{1 + \alpha}\leq \frac{1}{2}|q' - p'| + \frac{1}{4}|q_d - p_d|.
\end{align*}
Rearranging yields $|q_d -p_d|\leq |q' - p'|$, so two free boundary points with same $x'$-coordinate must coincide. Together with existence of each vertical line, we've shown
\[\p\{w_j > 0\} \cap Q = \{(x', f(x')) : x' \in B_{1/8}'\}, \qquad |f(y') - f(x')|\leq |y' - x'|.\] 

\emph{Step 3.} We prove $f$ is differentiable and $f \in C^{1, \alpha}.$ Put $p = (x', f(x'))$, $q = (y', f(y')).$ \eqref{eq:scalar-flat-blowup-regularity-mid4} now gives
\[|\nu_{p, \infty}'\cdot (y' - x') + (\nu_{p, \infty})_d(f(y') - f(x'))|\leq 2M \ep |y' - x'|^{1 + \alpha}.\] Using Lipschitz bound, $|q - p| \leq \sqrt{2}|y' - x'|.$ Dividing by \((\nu_{p, \infty})_d\) we obtain 
\begin{align*}
    \left|f(y') - f(x') + \frac{\nu_{p, \infty}'}{(\nu_{p, \infty})_d}\cdot (y' - x')\right|&\leq M_1 \ep |y' - x'|^{1 + \alpha}, \qquad M_1 = M_1(d, \gamma)> 0.
\end{align*}
Now $f$ is differentiable at $x'$ and $\nabla f(x') = -\frac{\nu_{p, \infty}'}{(\nu_{p, \infty})_d}.$ Then $f \in C^{1, \alpha}(B_{1/32}')$, $\nu_{p, \infty} = \frac{(-\nabla f(x'), 1)}{\sqrt{1 + |\nabla f(x')|^2}}$ and 
\[|\nabla f(y') - \nabla f(x')|\leq M_2(d, \gamma)\ep |y' - x'|^{ \alpha}, \qquad |\nu_{q, \infty} - \nu_{p, \infty}|\leq M_2 \ep |p - q|^{\alpha}.\]

\emph{Step 4.} Scaling back gives
$C^{1,\alpha}$ regularity of $\p\{g > 0\}$ near $x_0$. Let $n(x_0)$ denote
the unit normal pointing into $\{g>0\}$.
Iterating directly at $0$ for $w_j$ with initial flatness $\ep$
gives directions $\nu_{j,m}$ satisfying
\[
 \nu_{j,0}=\nu,\qquad
 |\nu_{j,m+1}-\nu_{j,m}|
 \le C\ep\rho^{m\alpha}.
\]
By Step~3, $\lim \limits_{m \to \infty}\nu_{j,m}$ is the inward
unit normal at the origin of $\partial\{w_j>0\}$, which equals $n(x_0)$ under translation and positive dilation. Hence
\[
 |n(x_0)-\nu|
 \le \sum_{m=0}^{\infty}
       |\nu_{j,m+1}-\nu_{j,m}|
 \le \frac{C\ep}{1-\rho^\alpha}.
\]
For every small $\ep>0$, \eqref{eq:scalar-flat-blowup} allows us to choose $j=j(\ep)$ so that
the preceding argument applies. Since $n(x_0)$
is independent of $j$, letting $\ep\downarrow0$
gives $n(x_0)=\nu$.
\end{proof}

\bibliographystyle{alpha}
\bibliography{references}
\end{document}